\documentclass[10pt]{article}
\usepackage[utf8]{inputenc}
\usepackage{amssymb}
\usepackage{amsmath,bbm,mathabx,cases,booktabs,amsthm}
\usepackage{enumitem}   
\usepackage{mathrsfs}
\usepackage{caption,subcaption}
\usepackage{graphicx,multirow}
\usepackage{cite}
\usepackage{hyphenat,comment}
\usepackage{authblk}

\usepackage[margin=1.4in]{geometry}

\newtheorem{theorem}{Theorem}[section]
\newtheorem{lemma}[theorem]{Lemma}

\newtheorem{proposition}[theorem]{Proposition}

\theoremstyle{definition}

\newtheorem{example}[theorem]{Example}

\theoremstyle{remark}
\newtheorem{remark}[theorem]{Remark}

\numberwithin{equation}{section}

\begin{document}

\author[1]{Iain Henderson \thanks{henderso@insa-toulouse.fr}}
\affil[1]{Institut de Mathématiques de Toulouse, Université Paul Sabatier, 118 Rte de Narbonne, 31400 Toulouse, France}
\title{Sobolev regularity of Gaussian random fields}

\maketitle
\begin{abstract}
    In this article, we fully characterize the measurable Gaussian processes $(U(x))_{x\in\mathcal{D}}$ whose sample paths lie in the Sobolev space of integer order $W^{m,p}(\mathcal{D}),\ m\in\mathbb{N}_0,\ 1 <p<+\infty$, where $\mathcal{D}$ is an arbitrary open set. The result is phrased in terms of a form of Sobolev regularity of the covariance function on the diagonal. This is then linked to the existence of suitable Mercer or otherwise nuclear decompositions of the integral operators associated to the covariance function and its cross-derivatives. In the Hilbert case $p=2$, additional links are made w.r.t. the Mercer decompositions of the said integral operators, their trace and the imbedding of the RKHS in $W^{m,2}(\mathcal{D})$. We provide simple examples and partially recover recent results pertaining to the Sobolev regularity of Gaussian processes.
\end{abstract}
\section{Introduction}
Sobolev spaces $W^{m,p}(\mathcal{D})$ are central tools in modern mathematics, most notably in the study of partial differential equations (PDEs). These spaces are built upon the notion of weak derivative: $v$ is the weak derivative of $u$ in the direction $x_i$ if for all smooth compactly supported function $\varphi\in C_c^{\infty}(\mathcal{D})$,
\begin{align}
     \int_{\mathcal{D}} u(x)\frac{\partial\varphi}{\partial x_i}(x)dx = - \int_{\mathcal{D}} v(x)\varphi(x)dx
\end{align}
Weak derivatives generalize classical, pointwise defined derivatives. In particular, there are cases where weak derivatives are well defined and pointwise differentiation otherwise fails (see e.g. \cite{evans1998}, Examples 3 and 4 p. 260). The popularity of Sobolev spaces is justified by a number of reasons: first, they have good topological and geometrical properties. They are separable reflexive Banach spaces when $1 < p < +\infty$, and separable Hilbert spaces when $p=2$ (\cite{fournier_adams_sobolev}, Theorem 3.6 p. 61). Through duality, this allows for geometrical interpretations of PDEs which in turn lead to numerous quantitative theoretical results in the study of PDEs\cite{evans1998}. Second, as the Sobolev norm is defined through integrals of powers of the function and its weak derivatives, it is easily interpreted as an energy functional of the said function, which complies with physical interpretations of PDEs. This is a desirable feature as PDEs are generally used for describing physical phenomena. Finally, Sobolev spaces are useful for practical purposes as they are the natural mathematical framework for the celebrated finite element method when seeking numerical solutions to PDEs (\cite{brenner2008mathematical}, Chapter 1).

When a function of interest $u :\mathcal{D} \rightarrow \mathbb{R}$ is unknown, it may be modelled as a sample path of a random field $(U(x))_{x\in\mathcal{D}}$, say a Gaussian process, whose sample paths lie in a suitable function space. This is e.g. frequent in Bayesian inference of functions \cite{van2011information}. Such suitable spaces can indeed happen to be Sobolev spaces, e.g. when $u$ describes a physical quantity. The question at hand in this article is thus the following: when do the sample paths of a given Gaussian process lie in some Sobolev space? This question is closely linked to the recent attention that Gaussian processes have drawn for tackling machine learning problems arising from PDE models; see e.g.  \cite{raissi2017,wahlstrom2013,owhadi_bayes_homog}. Notably (see \cite{CHEN_owhadi_2021}), Gaussian processes seem to provide a numerically competitive and mathematically tractable alternative to the now widespread "physics informed neural networks" (PINNs, \cite{pinns}). For the moment though, the machine learning techniques involving Gaussian processes have only been studied within the framework of spaces of functions with classical smoothness : $C^0, C^1$, etc. As argued before, these spaces are often not as well-suited for studying PDEs as Sobolev spaces.

Though weak differentiability is more general, it is less direct to check than classical differentiability. Weak derivatives are defined implicitly and in the most general case, ensuring Sobolev regularity is not usually done by directly verifying that an integral or a series is finite, as would be the case in $L^p$ spaces. One reason for this is that the existence of the weak derivative has to ensured beforehand. To do so, variational or boundedness criteria are used instead (see Proposition \ref{lemma:linear_form_sobolev}).

In many important particular cases however, handy characterizations of such regularity do exist, which have effectively been used to bypass the implicit definition of Sobolev regularity and generate results on the sample path regularity of Gaussian processes. When $\mathcal{D}=\mathbb{R}^d$, the space $W^{m,p}(\mathbb{R}^d)$ can be characterized in terms of a sufficient decay of of the Fourier transform (\cite{stein1970singular}, Theorem 3 p. 135; \cite{evans1998}, Section 5.8.5; \cite{fournier_adams_sobolev}, Section 7.63). Still in the case $\mathcal{D}=\mathbb{R}^d$, Sobolev regularity is equivalent to the convergence of its de la Vallée Poussin expansion in a suitable space (\cite{nikol2012approximation}, Section 8.9). This fact has been the first to be employed for characterizing the Sobolev regularity of stationary Gaussian processes indexed by the unit cube of $\mathbb{R}^d$ in \cite{cirel1976norms, ibragimov1994conditions}, in terms of the spectral measure of its covariance. For some Banach spaces, explicit Schauder bases are known and lying in such spaces can be translated as the convergence of some coordinate series. This has been exploited in \cite{ciesielski1993quelques} for studying the Besov and Besov-Orlicz regularity of one dimensional Gaussian processes (they are natural generalizations of Sobolev regularity, \cite{fournier_adams_sobolev}). Wavelet analysis is also available for describing Sobolev regularity (\cite{fournier_adams_sobolev}, Section 7.70) and has been used for studying the smoothness of the Brownian motion \cite{ciesielski1991modulus,roynette1993mouvement}. More complex notions such as the existence of an underlying Dirichlet structure have aptly been put to use in  \cite{kerkyacharian2018regularity}. The latter work deals with Besov $B^s_{\infty,\infty}$ regularity, $s>0$, on compact metric spaces, and relies on a convergence analysis of suitable spectral coefficients, using the so called Littlewood-Paley decomposition. In
\cite{steinwart2019convergence}, Karhunen-Loève expansions are used to study whether or not the sample paths of a general second order random process lie in interpolation spaces between the reproducing kernel Hilbert space (RKHS, Section \ref{sub:rkhs} below) of the process and $L^2(\nu)$, where $\nu$ is a $\sigma$-finite measure. This is then applied to study $H^s$-regularity properties of the corresponding sample paths when $s>d/2$ (Corollary 4.5 and 5.7 in \cite{steinwart2019convergence}), with applications to Gaussian processes in particular. Using the notion of mean square derivatives, \cite{SCHEUERER2010} shows that the sample paths of a general second order random field lie in $W^{m,2}(\mathcal{D})$ under an integrability condition of the symmetric cross derivatives of the kernel over the diagonal (Theorem 1). This result strongly suggests that a purely spectral criteria for Sobolev regularity of a random process should exist as the integrals appearing in Theorem 1 of \cite{SCHEUERER2010} exactly correspond to the trace of specific integral operators which are naturally linked to the covariance of the process; in fact, we provide such a criteria in Proposition \ref{prop:sobolev_reg_sto}. For the suitable definition and use of the mean square derivatives of the process, \cite{SCHEUERER2010} additionally requires that the covariance function be continuous over the diagonal as well as its symmetric cross derivatives.

The purpose of this article is to uncover necessary and sufficient characterizations of the Sobolev regularity of positive integer order of a given Gaussian process, in terms of its covariance function. In an attempt to make them both as general and concise as possible, we set the following targets and assumptions.
\begin{itemize}[wide,labelwidth=0em,labelindent=0pt]
    \item The covariance function of the Gaussian process will only be assumed measurable, as in \cite{steinwart2019convergence}. This contrasts with some of the previously mentioned works \cite{ciesielski1993quelques,SCHEUERER2010,kerkyacharian2018regularity}, where the covariance function is assumed continuous. It seems though that assuming the continuity of the covariance (and thus more or less that of the sample paths, \cite{azais_level_2009} p. 31) to examine some Sobolev regularity of potentially low order is an unnatural hypothesis. This is especially true as the dimension of $\mathcal{D}$ increases, since $W^{m,p}(\mathcal{D})$ is embedded in $C_B^0(\mathcal{D})$, the Banach space of continuous and bounded functions over $\mathcal{D}$, only when $m>d/p$ (\cite{fournier_adams_sobolev}, Theorems 4.12 and 7.34).
    
    \item We will not make any regularity or shape assumptions on the open set $\mathcal{D}$. Indeed, Sobolev spaces of integer order are easily defined over arbitrary open sets $\mathcal{D}\subset\mathbb{R}^d$, and thus some results should exist within this general setting. As a result though, we will not deal with fractional Sobolev spaces nor Besov spaces. 
    Indeed, minimal imbedding properties of fractional Sobolev and Besov spaces can only be ensured under additional hypotheses on $\mathcal{D}$, namely enjoying a Lipschitz boundary or the cone condition (\cite{fournier_adams_sobolev}, Theorems 4.12 and 7.34; see also Remark 6.47(1)). For example, one may inconveniently have $W^{s,p}(\mathcal{D}) \not\subset W^{1,p}(\mathcal{D})$ for some $0<s<1$ when the boundary of $\Omega$ is not Lipschitz ( \cite{di2012hitchhiker}, Example 9.1). We will see that elementary characterizations of Sobolev regularity (Lemmas \ref{lemma:linear_form_sobolev} and \ref{lemma:E_et_F_countable_sobolev}) will prove to be enough for our purpose.
    
    \item Our results should lie outside of the assumption that $m > d/p$, where $m,p$ and $d$ correspond to the notation $W^{m,p}(\mathcal{D}),\ \mathcal{D} \subset \mathbb{R}^d$. Indeed, many previous results concerning the Sobolev regularity of a given Gaussian process concern the spaces $H^m(\mathcal{D}) =  W^{m,2}(\mathcal{D})$, $\mathcal{D}\subset\mathbb{R}^d$, only in the case $m > d/2$. This is convenient because it ensures that $H^m(\mathcal{D})$ is continuously embedded in $C_B^0(\mathcal{D})$ when $\mathcal{D}$ is smooth enough, which suppresses the ambiguity of choosing a representer of a function in $H^m(\mathcal{D})$. More specifically, the spaces $H^m(\mathcal{D})$ are actually reproducing kernel Hilbert spaces. But $m > d/2$ excludes the spaces $H^1(\mathbb{R}^2)$ and $H^1(\mathbb{R}^3)$, which are central in the study of many important second order PDEs such as the wave equation, the heat equation, Laplace's equation or Schrödinger's equation.
\end{itemize}
Our characterization of measurable Gaussian processes with sample paths in $W^{m,p}(\mathcal{D})$ is phrased in terms of a form of Sobolev regularity of the covariance function on the diagonal. It is then linked to the existence of suitable Mercer or otherwise nuclear decompositions of the integral operators associated to the covariance function and its symmetric weak cross-derivatives. In the Hilbert case $p=2$, additional links are made w.r.t. the Mercer decompositions of the said integral operators, their trace and the Hilbert-Schmidt nature of the imbedding of the RKHS in $W^{m,2}(\mathcal{D})$. Our results are strongly reminiscient of those found in \cite{SCHEUERER2010}. In particular, this shows that contrarily to what is suggested in \cite{steinwart2019convergence}, p. 370, the Sobolev regularity of the sample paths of a given Gaussian process is not about $d/2$ less than that of the functions of its RKHS. This regularity is rather characterized by purely spectral properties of the covariance operator of the associated Gaussian measure. It just happens that in many standard cases such as with the Matérn kernels of order $\nu$ on "nice" bounded domain $\mathcal{D} \subset \mathbb{R}^d$, their RKHS happens to be $H^{\nu+d/2}(\mathcal{D})$ (\cite{steinwart2019convergence}, Example 4.8) and the imbedding of $H^{\nu+d/2}(\mathcal{D})$ in $H^{s}(\mathcal{D})$ is Hilbert-Schmidt when $s < \nu$. See Example \ref{ex:compare_steinwart} for further details.

The article is organized as follow. In section \ref{sec:preliminaries}, we introduce the necessary notions for properly stating our results as well as some useful lemmas directly related to these notions. In sections \ref{sec:gp_sobolev_LP} and \ref{sec:gp_sobolev_L2}, we state and prove the main results of this article, which treat the general case $p\in(1,+\infty)$ and the special case $p=2$ respectively. In section \ref{sec:ccl}, we conclude and provide some further outlooks. In section \ref{sec:proofs}, we prove the intermediary lemmas used in the main proofs.

\paragraph{Notations} Given a Banach space $X$, $X^*$ denotes its topological dual. Given $x\in X$ and $l\in X^*$, we denote the duality bracket as follow: $l(x) = \langle l,x\rangle_{X^*,X}$. $\mathcal{B}(X)$ denotes the Borel $\sigma$-algebra of $X$ for its norm topology. Given two linear operators $A : X_1\rightarrow Y_1$ and $B : X_2\rightarrow Y_2$, $A\otimes B : X_1\otimes X_2 \rightarrow Y_1\otimes Y_2$ denotes their tensor product which verifies $(A\otimes B)(a\otimes b) = (Aa)\otimes(Bb)$. Given two real valued functions $f$ and $g$, $f\otimes g$ denotes their tensor product defined by $(f\otimes g)(x,y) = f(x)g(y)$. Given $h\in \mathbb{R}^d$, $|h|$ denotes its Euclidean norm. Given $p \in (1,+\infty)$, $q$ will always denote its conjugate: $1/p+1/q=1$ i.e. $q=p/(p-1)$. As usual, when $\mathcal{D}$ is an open set of $\mathbb{R}^d$, we identify the dual of $L^p(\mathcal{D})$ with $L^q(\mathcal{D})$. Explicitly, if $f\in L^p(\mathcal{D})$ and $g\in L^q(\mathcal{D})$, we have
\begin{align}
    \langle f,g\rangle_{L^p,L^q} = \int_{\mathcal{D}}f(x)g(x)dx = \langle g,f\rangle_{L^q,L^p}
\end{align}
When there is no risk of confusion, we will write $||f||_p := ||f||_{L^p(\mathcal{D})}$. If $H$ is a Hilbert space, $\langle\cdot,\cdot\rangle_H$ denotes its inner product. We denote $\mathbb{N}:=\{1,2,...\}$ the set of natural numbers and $\mathbb{N}_0:=\mathbb{N}\cup\{0\}$. Given an open set $\mathcal{D} \subset \mathbb{R}^d$, we write $\mathcal{D}_0 \Subset \mathcal{D}$ if $\mathcal{D}_0 \subset \mathcal{D}$ and $\overline{\mathcal{D}_0}$ is compact. $\lambda^d$ denotes the Lebesgue measure over $\mathbb{R}^d$. $L_{loc}^1(\mathcal{D})$ denotes the space of equivalence classes of locally integrable functions over $\mathcal{D}$, i.e. such that $\int_K|f(x)|dx < +\infty$ for all $K\Subset \mathcal{D}$. Elements of $L_{loc}^1(\mathcal{D})$ are identified when they are equal almost everywhere w.r.t. the Lebesgue measure ($L_{loc}^1(\mathcal{D})$ is a large space: $L^p(\mathcal{D})\subset L_{loc}^1(\mathcal{D})$ for all $p\in[1,+\infty]$). Given an equivalence class $f\in L_{loc}^1(\mathcal{D})$, a \textit{representer} of $f$ is a function $\widehat{f} : \mathcal{D} \rightarrow \mathbb{R}$ such that the equivalence class of $\widehat{f}$ in $L_{loc}^1(\mathcal{D})$ is $f$. We will sometimes denote $f$ and $\widehat{f}$ with the same symbol, e.g. $f$.
Given a function $k$ defined over $\mathcal{D}\times\mathcal{D}$, $\mathcal{E}_k$ denotes the associated integral operator (if well defined):
\begin{align}
 (\mathcal{E}_kf)(x) = \int_{\mathcal{D}}k(x,y)f(y)dy   
\end{align}
The input and output spaces of $\mathcal{E}_k$ will be specified on a case-by-case basis. 
\section{Preliminary notions and results}\label{sec:preliminaries}
In Sections \ref{sub:def_sobolev}, \ref{sub:carac_sobolev_1} and \ref{sub:sobolev_vs_distrib}, we define the notion of Sobolev regularity through the prisms of weak derivatives and generalized functions, and provide handy characterizations of this regularity. We present key notions from operator theory in Section \ref{sub:op_th}. In Section \ref{sub:gp_gm}, we define and provide some useful results related to Gaussian processes and Gaussian measures.
\subsection{Definition of weak derivatives and Sobolev spaces}\label{sub:def_sobolev} Let $\alpha = (\alpha_1,...,\alpha_d) \in \mathbb{N}^d$. Denote $\partial^{\alpha} = \partial_{x_1}^{\alpha_1}...\partial_{x_d}^{\alpha_d}$ the $\alpha^{th}$ derivative, and $|\alpha|:=\sum_{i=1}^d|\alpha_i|$. In this article, the statement "let $|\alpha|\leq m$" will mean "let $\alpha = (\alpha_1,...,\alpha_d) \in \mathbb{N}^d$ be such that $|\alpha|\leq m$". Given a function $k$ defined on $\mathcal{D}\times\mathcal{D}$, $\partial^{\alpha,\alpha}k$ denotes its symmetric cross derivative: $\partial^{\alpha,\alpha}k(x,y) := \partial_{x_1}^{\alpha_1}...\partial_{x_d}^{\alpha_d}\partial_{y_1}^{\alpha_1}...\partial_{y_d}^{\alpha_d}k(x,y)$ (formally, $\partial^{\alpha,\alpha} = \partial^{\alpha} \otimes\partial^{\alpha}$). A function $u \in L_{loc}^1(\mathcal{D})$ has $v \in L_{loc}^1(\mathcal{D})$ for its $\alpha^{th}$ weak derivative if (\cite{fournier_adams_sobolev}, section 1.62)
\begin{align}
\forall \varphi \in C_c^{\infty}(\mathcal{D}),\ \ \ \int_{\mathcal{D}}u(x)\partial^{\alpha}\varphi(x)dx = (-1)^{|\alpha|}\int_{\mathcal{D}}v(x)\varphi(x)dx
\end{align}
$v$ is then unique in $L_{loc}^1(\mathcal{D})$ and is denoted $v = \partial^{\alpha}u$.
Let $p\in[1,+\infty]$. The Sobolev space $W^{m,p}(\mathcal{D})$ is defined as (\cite{fournier_adams_sobolev}, section 3.2)
\begin{align}
W^{m,p}(\mathcal{D}) = \{u \in L^p(\mathcal{D}) : \forall \ |\alpha| \leq m, \partial^{\alpha} u \in L^p(\mathcal{D}) \}
\end{align}
Sobolev spaces are Banach spaces for the norm $||u||_{W^{m,p}} := (\sum_{|\alpha|\leq m} ||\partial^{\alpha}u||_{p}^p)^{1/p}$; they are separable when $p \neq +\infty$ (\cite{fournier_adams_sobolev}, Theorem 3.6 p. 61). When $p=2$, $W^{m,p}$ is usually denoted $H^m(\mathcal{D})$ and is a Hilbert space for the following inner product
\begin{align}
    \langle u,v\rangle_{H^m(\mathcal{D})} := \sum_{|\alpha|\leq m} \langle \partial^{\alpha}u,\partial^{\alpha}v\rangle_{L^2(\mathcal{D})}
\end{align}
Note that we made no assumptions on the regularity of the open set $\mathcal{D}$.

\subsection{Characterization of $W^{m,p}$-regularity for locally integrable functions}\label{sub:carac_sobolev_1}
As for pointwise derivatives, finite difference operators can be used for characterizing Sobolev regularity. 
Given $h \in \mathbb{R}^d$, introduce the translation operator $(\tau_h u)(x) = u(x+h)$, which is bounded over $L^p(\mathbb{R}^d)$. 
Introduce the associated "finite difference operator":
\begin{align}
\Delta_h = {\tau_{h} - Id}
\end{align}
The linear subspace of bounded operators over $L^p(\mathbb{R}^d)$ induced by the translation operators is commutative, as $\tau_{h_1}\circ\tau_{h_2} = \tau_{h_1+h_2} = \tau_{h_2}\circ\tau_{h_1}$. Let $h = (h_1,...,h_m)\in (\mathbb{R}^d)^m$, we define the $m^{th}$ order finite difference operator associated to $h$ to be $\Delta_h := \prod_{i=1}^m \Delta_{h_i}$ where the product symbol denotes the composition of operators. 
When $h \in \mathbb{R}^d$, the adjoint of $\Delta_h$ is also a finite difference operator, which is computable using the change of variable formula. If $h \in \mathbb{R}^d$, then
\begin{align}
\Delta_h^* = {\tau_{-h} - Id}
\end{align}
Finally, when $\alpha = (\alpha_1,...,\alpha_d)\in\mathbb{N}^d$ and $h = (h_1,...,h_d)\in(\mathbb{R}_+^*)^{d}$, we denote by $\delta_h^{\alpha}$ the finite difference approximation of $\partial^{\alpha}$ defined by
\begin{align}\label{eq:def_fd_approx_alpha}
    \delta_h^{\alpha} = \prod_{i=1}^{d}\bigg( \frac{\Delta_{h_ie_{i}}}{h_i}\bigg)^{\alpha_i} = \bigg( \frac{\Delta_{h_1e_{1}}}{h_1}\bigg)^{\alpha_1}\dotsb\ \bigg( \frac{\Delta_{h_de_{d}}}{h_d}\bigg)^{\alpha_d}
\end{align}
Above, $(e_1,...,e_d)$ is the canonical basis of $\mathbb{R}^d$. Depending on which one is the most convenient, we will either use $\Delta_h$ or $\delta_h^{\alpha}$.
We shall use the following characterizations of $W^{m,p}$-regularity, which are straightforward generalizations of Proposition 9.3 from \cite{brezis2011functional} to multiple derivatives. 

\begin{lemma}\label{lemma:linear_form_sobolev}
Suppose that $u \in L_{loc}^1(\mathcal{D})$. Let $m\in \mathbb{N}_0$, $p \in (1,+\infty]$ and introduce $q\geq 1$ the conjugate of $p$ : $1/p + 1/q = 1$. Then the following statements are equivalent
\begin{enumerate}[label=(\roman*)] 
\item $u \in W^{m,p}(\mathcal{D})$
\item (Variational control) for all $\alpha$ such that $|\alpha| \leq m$, there exists a constant $C_{\alpha}$ such that
\begin{align}\label{eq:LP_control}
\forall \varphi \in C_c^{\infty}(\mathcal{D}),\ \ \ \Big|\int_{\mathcal{D}}u(x)\partial^{\alpha}\varphi(x)dx \Big| \leq C_{\alpha} ||\varphi||_{L^q(\mathcal{D})}
\end{align}
In this case, the $L^p$ norm of $\partial^{\alpha}u$ is given by
\begin{align}
||\partial^{\alpha}u||_{L^p(\mathcal{D})} = \sup_{\varphi \in C_c^{\infty}(\mathcal{D})\setminus\{0\}} \bigg|\int_{\mathcal{D}}u(x)\frac{\partial^{\alpha}\varphi(x)}{||\varphi||_{L^q}}dx \bigg|
\end{align}
\item \label{point:fd_control} (Finite difference control) there exists a constant $C$ such that for all open set $\mathcal{D}_0\Subset \mathcal{D}$, for all $l\leq m$ and all $h = (h_1,...,h_{l}) \in (\mathbb{R}^d)^{l}$ such that $\sum_i |h_i| < \text{dist}(\mathcal{D}_0,\partial \mathcal{D})$,
\begin{align}\label{eq:big_delta_control_sobolev}
||\Delta_{h}u||_{L^p(\mathcal{D}_0)} \leq C |h_1|\times ...\times |h_{l}|
\end{align}
Moreover, one can take $C = ||u||_{W^{m,p}(\mathcal{D})}$.
\end{enumerate}
In addition, $u\in W^{m,p}(\mathcal{D})$ then for all $|\alpha|\leq m$ and $\mathcal{D}_0\Subset\mathcal{D}$,
\begin{align}\label{eq:fd_cv_LP}
    ||\delta_h^{\alpha}u - \partial^{\alpha}u||_{L^p(\mathcal{D}_0)} \rightarrow 0
\end{align}
when $h\rightarrow 0$, with $h=(h_1,...,h_d)\in(\mathbb{R}_+)^d$ suitably chosen so that $\delta_h^{\alpha}u$ makes sense.
\end{lemma}
In Point $\ref{point:fd_control}$ above, the assumption that $\sum_i |h_i| < \text{dist}(\mathcal{D}_0,\partial \mathcal{D})$ is only there to ensure that the quantity $\Delta_{h}u(x)$ makes sense when $x\in\mathcal{D}_0$.
\subsection{Sobolev regularity and generalized functions}\label{sub:sobolev_vs_distrib}
The theory of generalized functions (or distributions) provides a flexible way of characterizing Sobolev regularity, by building a larger space in which weak derivatives are always defined. Given an open set $\mathcal{D}$, denote $C_c^{\infty}(\mathcal{D})$ the space of smooth functions with compact support in $\mathcal{D}$. 
Endow it with its usual LF topology, defined e.g. in \cite{treves2006topological}, Chapter 13. This topology is such that the sequence $(\varphi_n)$ converges to $\varphi$ in $C_c^{\infty}(\mathcal{D})$ if and only if there exists a compact set $K \subset \mathcal{D}$ such that $\text{Supp}(\varphi_n) \subset K$ for all $n$ and
\begin{align}
    \forall \alpha = (\alpha_1,...,\alpha_d) \in \mathbb{N}^d, \ \ \ \sup_{x\in K}|\partial^{\alpha}\varphi_n(x) - \partial^{\alpha}\varphi(x)| \longrightarrow 0
\end{align}
Here, $\partial^{\alpha} := \partial_{x_1}^{\alpha_1}\dotsc\partial_{x_d}^{\alpha_d}$. With $C_c^{\infty}(\mathcal{D})$ endowed with this topology, the space of generalized functions, or distributions, is then defined as the topological dual of $C_c^{\infty}(\mathcal{D})$ i.e. the vector space of all continuous linear forms over $C_c^{\infty}(\mathcal{D})$. It is traditionally denoted as follow: $\mathscr{D}'(\mathcal{D}) := C_c^{\infty}(\mathcal{D})'$ (\cite{treves2006topological}, Notation 21.1). A generalized function $T\in\mathscr{D}'(\mathcal{D})$ is said to be regular (\cite{treves2006topological}, p. 224) if it is of the form
\begin{align}
    \forall \varphi\in C_c^{\infty }(\mathcal{D}), \ \ \ T(\varphi) = \int_{\mathcal{D}}u(x)\varphi(x)dx
\end{align}
for some $u\in L_{loc}^1(\mathcal{D})$, in which case one writes $T=T_u$.
Given any function $u\in L_{loc}^1(\mathcal{D})$ and $\alpha\in \mathbb{N}^d$, its distributional derivative $D^{\alpha}u$ is \textit{defined} by the following formula (\cite{treves2006topological}, pp. 248-250):
\begin{align}
    D^{\alpha}u : \varphi \longmapsto (-1)^{|\alpha|} \int_{\mathcal{D}}\partial^{\alpha}\varphi(x)u(x)dx
\end{align}
$D^{\alpha}u$ then also lies in $\mathscr{D}'(\mathcal{D})$.
Sobolev regularity can now be rephrased as follow : $u$ lies in $ W^{m,p}(\mathcal{D})$ iff for all $|\alpha|\leq m$, the distributional derivative $D^{\alpha}u$ is in fact a regular generalized function represented by some $v_{\alpha}\in L^p(\mathcal{D})$ i.e. $D^{\alpha}u = T_{v_{\alpha}}$. Then $v_{\alpha}$ is unique in $L^p(\mathcal{D})$ and $\partial^{\alpha}u = v_{\alpha}$ in $L^p(\mathcal{D})$, where $\partial^{\alpha}u$ is the $\alpha^{th}$ weak derivative of $u$. 

Moreover, the control equation \eqref{eq:LP_control} shows that $\partial^{\alpha}u$ exists and lies in $L^p(\mathcal{D})$ if and only if $D^{\alpha}u : C_c^{\infty}(\mathcal{D}) \rightarrow \mathbb{R}$ can be extended as a continuous linear form over $L^q(\mathcal{D})$. 
Ensuring the existence of such extensions will thus be of prime interest for us, and is the topic of the next lemma.
Specifically, the next result states that given continuous linear or bilinear forms over $C_c^{\infty}(\mathcal{D})$, the existence of extensions of these maps to $L^q(\mathcal{D})$ can be ensured by obtaining suitable estimates on a well chosen \textit{countable} set $E_q \subset C_c^{\infty}(\mathcal{D})$. Restricting ourselves to $E_q$ will allow us to eliminate any measurability issues when introducing the supremum of certain random variables indexed by $E_q$, as a countable supremum of random variables remains a random variable.
Below, we write $||\cdot||_q := ||\cdot||_{L^q(\mathcal{D})}$ for short.
\begin{lemma}[Extending continuous linear and bilinear forms over $C_c^{\infty}(\mathcal{D})$ to $L^p(\mathcal{D})$]\label{lemma:E_et_F_countable_line_bilin}
Let $p \in (1,+\infty)$.
There exists a countable $\mathbb{Q}-$vector space $E_q = \{\Phi_n^q, n \in \mathbb{N}_0\} \subset C_c^{\infty}(\mathcal{D})$ with the following property. 

\begin{enumerate}[label=(\roman*),wide,labelwidth=0em,labelindent=0pt]
\item A distribution $T \in \mathscr{D}'(\mathcal{D})$ is a regular distribution, $T = T_v$, for some $v\in L^p(\mathcal{D})$ iff it verifies the countable estimate for some constant $C > 0$
\begin{align}\label{eq:LP_control_distrib}
\forall \varphi \in E_q, \ \ |T(\varphi)| \leq C ||\varphi||_{q}
\end{align}
or equivalently, $\sup_{n \in \mathbb{N}} |T(\Phi_n^q)|/||\Phi_n^q||_{q} < + \infty$ (here, setting $\Phi_0^p = 0$ without loss of generality).
This is equivalent to $T$ admitting an extension over $L^q(\mathcal{D})$ which is then uniquely given by $T(f) = \int_{\mathcal{D}} f(x)v(x)dx$. Moreover,
\begin{align}
\sup_{n \in \mathbb{N}} \frac{|T(\Phi_n^q)|}{||\Phi_n^q||_{q}} = \sup_{\varphi \in C_c^{\infty}(\mathcal{D})} \frac{|T(\varphi)|}{||\varphi||_{q}}
\end{align}
whether these quantities are finite or not.
\item Let $b$ be a continuous bilinear form over $C_c^{\infty}(\mathcal{D})$. Then  $b$ can be extended to a continuous bilinear form over $L^q(\mathcal{D})$ iff it verifies the countable estimate 
\begin{align}\label{eq:bilin_Lq}
\forall \varphi,\psi \in E_q, \ \ |b(\varphi,\psi)| \leq C ||\varphi||_{q}||\psi||_{q}  
\end{align}
In this case, such an extension is unique and there will exist a unique bounded operator $B : L^q(\mathcal{D}) \rightarrow L^p(\mathcal{D})$ verifying the following identity
\begin{align}
\forall \varphi,\psi \in C_c^{\infty}(\mathcal{D}), \ \ b(\varphi,\psi) = \langle B\varphi,\psi\rangle_{L^p,L^q}
\end{align}
\end{enumerate}
\end{lemma}
The proof of this result can be found in the appendix. It is based on Lemma \ref{lemma:d_separable} below, which is interesting in itself. Recall that a topological space $X$ is \textit{separable} if there exists a countable subset $Y\subset X$ which is dense in $X$ for the topology of $X$. Then the following holds.
\begin{lemma}\label{lemma:d_separable}
$C_c^{\infty}(\mathcal{D})$ endowed with its LF-topology is separable. 
\end{lemma}
A short proof of this result can be found in \cite{henderson_arxiv}, p. 16. See also \cite{gelfand1964generalized}, p. 73, (3) for a statement of this result, or \cite{gapaillard8processus}, Corollaire (1).2, p. 78 for a reference in French.
Given the set $E_q$ provided by Lemma \ref{lemma:E_et_F_countable_line_bilin}, we next define the countable set $F_q$ to be
\begin{align}
F_q := \{\varphi/||\varphi||_{q}, \varphi \in E_q, \varphi \neq 0\} = \{f_n^q, n \in\mathbb{N}\} \subset S_{q}(0,1)
\end{align}
Above, $(f_n^q)_{n\in\mathbb{N}}$ is an enumeration of $F_q$ and $S_{q}(0,1)$ is the unit sphere of $L^q(\mathcal{D})$. The next lemma is then a direct consequence of Lemmas \ref{lemma:linear_form_sobolev} and \ref{lemma:E_et_F_countable_line_bilin}.

\begin{lemma}[Countable characterization of Sobolev regularity]\label{lemma:E_et_F_countable_sobolev} Let $p\in(1,+\infty)$. For any $u \in L_{loc}^1(\mathcal{D})$, $u$ lies in $W^{m,p}(\mathcal{D})$ iff for all multi index $\alpha$ such that $|\alpha| \leq m$, there exists a constant $C_{\alpha}$ such that
\begin{align}\label{eq:LP_control_den}
\forall \varphi \in E_q,\ \ \ \bigg|\int_{\mathcal{D}}u(x)\partial^{\alpha}\varphi(x)dx \bigg| \leq C_{\alpha} ||\varphi||_{q}
\end{align}
Or equivalently,
\begin{align}\label{eq:LP_control_sup}
\sup_{\varphi \in F_q} \Bigg|\int_{\mathcal{D}}u(x)\partial^{\alpha}\varphi(x)dx\Bigg| = \sup_{n \in \mathbb{N}} \Bigg|\int_{\mathcal{D}}u(x)\partial^{\alpha}f_n^q(x)dx \Bigg| < + \infty
\end{align}
Moreover,
\begin{align}
\sup_{\varphi \in F_q} \Bigg|\int_{\mathcal{D}}u(x)\partial^{\alpha}\varphi(x)dx\Bigg|  = \sup_{\varphi \in C_c^{\infty}(\mathcal{D})\setminus\{0\}} \Bigg|\int_{\mathcal{D}}u(x)\frac{\partial^{\alpha}\varphi(x)}{||\varphi||_{q}}dx \Bigg| 
\end{align}
whether these quantities are finite or not. If one of them is finite, then it is equal to $||\partial^{\alpha}u||_{L^p(\mathcal{D})}$.
\end{lemma}
This lemma provides us with a somewhat explicit \textit{countable} criteria for Sobolev regularity, which is valid whatever the open set $\mathcal{D}$. This result is not surprising because the spaces $ W^{m,p}(\mathcal{D}),\ p\in(1,+\infty)$, are separable.


\subsection{Tools from operator theory}\label{sub:op_th}
The following reminders may be found in \cite{bogachev1998gaussian}, Section A.2. Let $H_1$ and $H_2$ be two Hilbert spaces, and $X$ and $Y$ two Banach spaces. 
\begin{enumerate}[label=(\roman*),wide,labelwidth=0em,labelindent=0pt]
\item A linear operator $T : X \rightarrow Y$ is bounded if $||T|| := 
\sup_{||x||_X=1} ||Tx||_{Y} < +\infty$. A bounded operator $T : X \rightarrow Y$ is compact if $\overline{T(B)}$ is a compact set of $Y$, where $B$ is the closed unit ball of $X$. When $X=Y$, the spectrum of a compact operator is purely discrete, and can be reordered as a sequence $(\lambda_n)_{n\in\mathbb{N}}$ which converges to $0$.
\item If $T : H_1 \rightarrow H_2$ is compact, then $T^*T : H_1 \rightarrow H_1$ is compact, self-adjoint and positive $(\forall x\in H_1, \langle x,T^*Tx\rangle_{H_1} \geq 0)$. If $H_1$ is separable, $T^*T$ can be diagonalized in an orthonormal basis $(e_n)$ of $H_1$. The positive eigenvalues of $T^*T$, $(s_n^2)$, are called the singular values of $T$. If $H_1$ is separable, $T$ is said to be Hilbert-Schmidt if $\sum_{n\in\mathbb{N}} ||Te_n||_{H_2}^2 < + \infty$ for one (equivalently, all) orthonormal basis $(e_n)$ of $H_1$. 
 Every Hilbert-Schmidt operator is compact, and every Hilbert-Schmidt operator $T$ acting on $L^2(\mathcal{D})$ can be written in integral form (\cite{bogachev1998gaussian}, Lemma A.2.13): there exists a "kernel" $k\in L^2(\mathcal{D}\times \mathcal{D})$ such that for all $f\in L^2(\mathcal{D})$,
\begin{align}
(Tf)(x) = \int_{\mathcal{D}}k(x,y)f(y)dy = (\mathcal{E}_kf)(x)
\end{align}
If $T$ is symmetric, positive and Hilbert-Schmidt, there exists an orthonormal basis $(\phi_n)$ of $L^2(\mathcal{D})$ comprised of eigenvectors of $T$ with positive eigenvalues $(\lambda_n)$,  such that in $L^2(\mathcal{D}\times \mathcal{D})$, we have
\begin{align}\label{eq:pos_sym_decomp}
k(x,y) = \sum_{n\in\mathbb{N}} \lambda_n\phi_n(x)\phi_n(y)
\end{align}
We will refer to decompositions of $f$ of the form of equation \eqref{eq:pos_sym_decomp} as \textit{Mercer decompositions}, in reference to the celebrated Mercer's theorem (\cite{brislawn1988kernels}, Theorem 1.2).
If $H_1$ is separable, $T$ is said to be trace-class (or nuclear) if
\begin{align}
\sum_{n\in\mathbb{N}}s_n < +\infty
\end{align}
One can then define its trace as the following linear functional, which is independent of the choice of basis $(e_n)$, and equal to the series of the eigenvalues of $T$ (Lidskii's theorem)
\begin{align}\label{eq:def_trace}
\text{Tr}(T) :=\sum_{n\in\mathbb{N}} \langle Te_n,e_n\rangle = \sum_{n\in\mathbb{N}} \lambda_n
\end{align}
Any trace-class operator is Hilbert-Schmidt, and $T$ is Hilbert-Schmidt if and only if $T^*T$ is trace-class. If $H_1 = H_2 = L^2(\mathcal{D})$, if $T$ is trace class with kernel $k$ and if $k$ is sufficiently smooth (say continuous), then the trace of $T =\mathcal{E}_k$ is given by
\begin{align}\label{eq:formula_trace_L2}
    \text{Tr}(T) = \int_{\mathcal{D}}k(x,x)dx
\end{align}
Extensions of the formula \eqref{eq:formula_trace_L2} to general Hilbert-Schmidt kernels $k\in L^2(\mathcal{D}\times\mathcal{D})$ of trace class operators is studied in \cite{brislawn1988kernels}; see also Proposition \ref{prop_bogachev_LP} below.
If $T : H_1 \rightarrow H_1$ is bounded, self-adjoint and positive, then we define its trace as the possibly infinite series of positive scalars $\text{Tr}(T) :=\sum_{n\in\mathbb{N}} \langle Te_n,e_n\rangle$.

\item (\cite{linde1980characterization}, p. 160) A bounded operator $T : X \rightarrow Y$ is nuclear if there exists sequences $(x_n) \subset X^*$ and $(y_n) \subset Y$ with $\sum_{n=1}^{+\infty}||x_n||_{X^*}||y_n||_Y <+\infty$ such that
\begin{align}
    \forall x \in X, \ \ \ Tx = \sum_{n=1}^{+\infty} \langle x_n,x\rangle_{X^*,X} y_n
\end{align}
In this case, we write abusively $T = \sum_{n=1}^{+\infty}x_n\otimes y_n$.
The nuclear norm of $T$ is then defined as
\begin{align}\label{eq:def_nuclear_norm}
    \nu(T) := \inf \bigg\{\sum_{n=1}^{+\infty}||x_n||_{X^*}||y_n||_Y \text{ such that } T = \sum_{n=1}^{+\infty}x_n\otimes y_n\bigg\}
\end{align}
A bounded operator $K : X^* \rightarrow X$ is said to be symmetric if for all $x,y \in X^*, \ \langle x,Ry\rangle = \langle y,Rx\rangle$ and positive if $\langle x,Rx\rangle \geq 0$. 
When $X = Y = H$, where $H$ is a separable Hilbert space, the sets of trace class and nuclear operators coincide; moreover, the same can be said for the trace functional \eqref{eq:def_trace} and the nuclear norm \eqref{eq:def_nuclear_norm} if $T$ has a positive spectrum : $\nu(T) = \text{Tr}(T)$. 
\end{enumerate}

\subsection{Gaussian processes and Gaussian measures over Banach spaces}\label{sub:gp_gm}
Throughout this article, $(\Omega,\mathcal{F},\mathbb{P})$ denotes the same probability space.
\begin{enumerate}[label=(\roman*),wide,labelwidth=0em,labelindent=0pt]
\item The \textit{law} $\mathbb{P}_X$ of a random variable $X:\Omega \rightarrow \mathbb{R}$, is the pushforward measure of $\mathbb{P}$ through $X$, which is defined by $\mathbb{P}_X(B) := \mathbb{P}(X^{-1}(B))$ for all Borel set $B\in\mathcal{B}(\mathbb{R})$ (\cite{bogachev2007measure}, Section 3.7).
\item A \textit{Gaussian process} (\cite{adler2007}, Section 1.2)  $(U(x))_{x\in\mathcal{D}}$ is a family of Gaussian random variables defined over $(\Omega,\mathcal{F},\mathbb{P})$ such that for all $n\in\mathbb{N}$, $(a_1,...,a_n)\in\mathbb{R}^n$ and $(x_1,...,x_n)\in\mathcal{D}^n, \sum_{i=1}^na_iU(x_i)$ is a Gaussian random variable. The law it induces over the function space $\mathbb{R}^{\mathcal{D}}$ endowed with its product $\sigma$-algebra is uniquely determined by its mean and covariance functions, $m(x) = \mathbb{E}[U(x)]$ and $k(x,x') = \text{Cov}(U(x),U(x'))$ (\cite{hoffman_jorgensen_book}, Section 9.8). We then write $(U(x))_{x\in\mathcal{D}}\sim GP(m,k)$. The covariance function $k$ is positive definite over $\mathcal{D}$, which means that for all non negative integer $n$ and $(x_1,...x_n)\in \mathcal{D}^n$, the matrix $(k(x_i,x_j))_{1\leq i,j\leq n}$ is non negative definite.
Conversely, given a positive definite function over an arbitrary set $\mathcal{D}$, there exists a centered Gaussian process indexed by $\mathcal{D}$ with the this function as its covariance function (\cite{adler2007}, p. 11). We shall often denote $\sigma(x) = k(x,x)^{1/2}$. Given $\omega\in\Omega$, the corresponding sample path of $(U(x))_{x\in\mathcal{D}}$ is the following deterministic function $U_{\omega}: \mathcal{D}\rightarrow \mathbb{R}$ defined by $U_{\omega}(x) := U(x)(\omega)$. A Gaussian process is said to be measurable if the map $(\Omega\times\mathcal{D},\mathcal{F}\otimes\mathcal{B}(\mathcal{D})) \rightarrow (\mathbb{R},\mathcal{B}(\mathbb{R})), (\omega,x) \mapsto U(x)(\omega)$ is measurable. If $(U(x))_{x\in\mathcal{D}}$ is measurable, then from Fubini's theorem the maps of the form $x\mapsto k(x,x'), x\mapsto k(x,x)$, etc, are measurable. We further discuss this assumption in Remark \ref{rk:GP_mes}.

We shall need the following lemma pertaining to the sample path-wise integration of Gaussian processes.
\begin{lemma}\label{lemma:u_alpha_phi_GP}
Let $\mathcal{D} \subset \mathbb{R}^d$ be an open set. Let $(U(x))_{x\in \mathcal{D}} \sim GP(0,k)$ be a measurable centered Gaussian process such that its standard deviation function $\sigma$ lies in $L_{loc}^1(\mathcal{D})$. Then the sample paths of $U$ lie in $L_{loc}^1(\mathcal{D})$ almost surely and given $\varphi \in C_c^{\infty}(\mathcal{D})$, the map defined by
\begin{align}
U_{\varphi}^{\alpha} : \Omega \ni \omega \longmapsto (-1)^{|\alpha|}\int_{\mathcal{D}}U(x)(\omega)\partial^{\alpha}\varphi(x)dx
\end{align}
is a Gaussian random variable. Moreover, for all $p\in(1,+\infty)$,  $(U_{\varphi}^{\alpha})_{\varphi\in F_q}$ is a centered Gaussian sequence (i.e. a Gaussian process indexed by $\mathbb{N})$.
\end{lemma}
We will also use the following fact about bounded Gaussian sequences, which can be seen as a weak form of Fernique's theorem (\cite{bogachev1998gaussian}, Theorem 2.8.5, p. 75).
\begin{lemma}[\cite{adler2007}, Theorem 2.1.2]\label{lemma:expo_moments}
Let $(U_n)_{n\in\mathbb{N}}$ be a Gaussian sequence and set $|U| := \sup_{n}|U_n|$. Suppose that $\mathbb{P}(|U| <+\infty) = 1$. Then there exists $\varepsilon > 0$ such that
\begin{align}
\mathbb{E}[\exp(\varepsilon |U|^2)] < + \infty
\end{align}
In particular, $\mathbb{E}[|U|^p] < + \infty$ for all $p \in \mathbb{N}$.
\end{lemma}

\item A \textit{Gaussian measure} $\mu$ (\cite{bogachev1998gaussian}, Definition 2.2.1) over a Banach space $X$ is a measure over its Borel $\sigma$-algebra such that given any $x \in X^*$, the pushforward measure of $\mu$ through the functional $x$ is a Gaussian measure over $\mathbb{R}$ (see Section \ref{sub:gp_gm}$(i)$ for a definition of the pushforward). Gaussian measures are equipped with a mean vector $a_{\mu} \in X^{**}$ and a covariance operator $K_{\mu} : X^* \rightarrow X^{**}$, defined in \cite{bogachev1998gaussian}, Definition 2.2.7. When $X$ is separable, $\mu$ is Radon (\cite{bogachev1998gaussian}, p. 125). This implies that $a_{\mu}$ lies in $X$ and that the covariance operator $K_{\mu}$ maps $X^*$ to $X$ (\cite{bogachev1998gaussian}, Theorem 3.2.3). The vector $a_{\mu}$ and the covariance operator $K_{\mu}$ are defined by the following formulas
\begin{align}
    \forall x \in X^*, \ \langle a_{\mu},x \rangle &=  \int_X \langle x, z\rangle\mu(dz) \\
    \forall x,y \in X^{*}, \langle y, K_{\mu}x\rangle &= \int_X \langle x - a_{\mu},z\rangle \ \langle y - a_{\mu},z\rangle\mu(dz)
\end{align}
In Propositions \ref{prop_type2} and \ref{prop_cotype2}, we present useful characterizations of Gaussian measures $\mu$ over two important classes of Banach spaces: spaces of type $2$ and cotype $2$ respectively. For a definition of spaces of type $2$ and cotype $2$, see e.g. \cite{chobanjan1977gaussian}. In this article, we will only use the fact that $L^p(\mathcal{D})$ is of type $2$ and cotype $p$ when $p \geq 2$, and cotype $2$ and type $p$ when $1\leq p\leq 2$ (see \cite{bogachev1998gaussian}, p. 152).
Moreover we shall restrict ourselves to the case where $X$ is separable. This implies that $\mu$ is Radon, which removes extension problems otherwise considered in \cite{linde1980characterization} and \cite{chobanjan1977gaussian}. 
\begin{proposition}[\cite{linde1980characterization}, Theorem 4]\label{prop_type2}
Let $X$ be a separable Banach space of type $2$, and let $\mu$ be a Gaussian measure over $X$. Then its covariance operator is nuclear. Conversely, given any $a\in X$ and any nuclear, symmetric, positive operator $K : X^*\rightarrow X$, there exists a Gaussian measure over $X$ with mean vector $a$ and covariance operator $K$.
\end{proposition}
Denote $l^2$ the Hilbert space of square summable sequences.
\begin{proposition}[\cite{chobanjan1977gaussian}, Theorem 4.1 and Corollary 4.1]\label{prop_cotype2}
Let $X$ be a separable Banach space of cotype $2$, and let $\mu$ be a Gaussian measure over $X$. Then there exists a continuous linear map $A : l^2 \rightarrow X$ and a trace-class operator $S$ over $l^2$ such that covariance operator of $\mu$ is given by $ASA^*$ (in particular, the covariance operator of $\mu$ is nuclear). In other words, $\mu$ is the pushforward measure of a Gaussian measure $\mu_0$ over $l^2$ through a bounded linear map $A$. Conversely, given any $a\in X$ and any operator of the form $ASA^*$ where $A : l^2 \rightarrow X$ is a bounded linear map and $S$ a trace class operator over $l^2$, there exists a Gaussian measure over $X$ with mean vector $a$ and covariance operator $K$.
\end{proposition}

These results generalize the case where $X$ is a separable Hilbert space, which can be found in \cite{bogachev1998gaussian}, Theorem 2.3.1. We finish with the following handy result describing centered Gaussian measures over $L^p$-spaces.
\begin{proposition}[\cite{bogachev1998gaussian}, Proposition 3.11.15]\label{prop_bogachev_LP}
Let $\mu$ be a centered Gaussian measure over $L^p(\mathcal{D})$ where $1 \leq p < +\infty$ and $\mathcal{D} \subset \mathbb{R}^d$ is an open set. Then there exists a positive definite function $k \in L^p(\mathcal{D}\times\mathcal{D})$ such that the covariance operator of $\mu$ is $\mathcal{E}_k : L^q(\mathcal{D}) \rightarrow L^p(\mathcal{D})$, the integral operator associated to $k$. Moreover, there exists a representer $\tilde{k}$ of $k$ in $L^p(\mathcal{D}\times\mathcal{D})$ which is the covariance function of a measurable Gaussian process $(U(x))_{x\in\mathcal{D}}$ whose sample paths lie in $L^p(\mathcal{D})$ a.s.. Additionally, $\tilde{k}$ verifies
\begin{align}\label{eq:sigma_LP_criteria}
    \int_{\mathcal{D}} \tilde{k}(x,x)^{p/2}dx = \int_{\mathcal{D}} \sigma(x)^{p}dx < +\infty
\end{align}
Finally, $\mathbb{P}_U = \mu$, where $\mathbb{P}_U$ is the pushforward of $\mathbb{P}$ through the Borel-measurable map $\omega \mapsto U_{\omega} \in L^p(\mathcal{D})$. Conversely, given any measurable positive definite function $k$ verifying \eqref{eq:sigma_LP_criteria}, the corresponding integral operator $\mathcal{E}_k : L^q(\mathcal{D}) \rightarrow L^p(\mathcal{D})$ is the covariance operator of a centered Gaussian measure $\mu$ over $L^p(\mathcal{D})$.
\end{proposition}
This result is quite strong, as it ensures the existence of a representer in $L^p(\mathcal{D}\times\mathcal{D})$ of the kernel of any Gaussian covariance operator, which is the covariance function of a measurable Gaussian process. This will enable us to remove awkward measurability issues w.r.t. $\sigma$ and equation \eqref{eq:sigma_LP_criteria}. Without the use of an underlying measurable Gaussian process, these issues are not trivial to deal with, see e.g. \cite{brislawn1988kernels} for an analysis of the Hilbert case $p=2$.
\begin{remark}\label{rk:GP_mes}
Proposition \ref{prop_bogachev_LP} shows that the assumption that a given Gaussian process is measurable is slightly less demanding that it might seem. Ensuring the existence of a measurable modification of a general random process is difficult outside of it being continuous in probability (\cite{doob_sto_pro}, Theorem 2.6 p. 61). Tedious extensions of this result exist (\cite{doob1937stochastic}, Theorem 2.3). For a Gaussian process $(U(x))_{x\in\mathcal{D}} \sim GP(0,k_u)$ however, Proposition \ref{prop_bogachev_LP} shows that the measurability of its covariance function over $\mathcal{D}\times\mathcal{D}$ and the integrability of its standard deviation in $L^p(\mathcal{D})$ (or equivalently, suitable nuclear decompositions of its associated integral operator $\mathcal{E}_k$, see Propositions \ref{prop_type2} and \ref{prop_cotype2}) ensure the existence of a measurable Gaussian process $(V(x))_{x\in\mathcal{D}}\sim GP(0,k_v)$ with the same covariance function in $L_{loc}^1(\mathcal{D}\times\mathcal{D})$. Consequently, $k_u=k_v$ a.e. on $\mathcal{D}\times\mathcal{D}$. Note though that the process $V$ need not be a modification of $U$. Since $k_u=k_v$ a.e., we only have that $U$ and $V$ have "$\lambda^d$-almost" the same finite dimensional marginals : for all $n\in\mathbb{N}$ and almost every $(x_1,...,x_n)\in\mathcal{D}^n$ (in the sense of the Lebesgue measure), $(U(x_1),...,U(x_n))$ and $(V(x_1),...,V(x_n))$ have the same law.
\end{remark}
Throughout this article, we will only consider centered Gaussian processes $(\mathbb{E}[U(x)]\equiv 0)$ and Gaussian measures $(a_{\mu}=0)$. Generalizations of the results of this article to non centered Gaussian processes are straightforward.
\end{enumerate}

\section{Sobolev regularity for Gaussian processes : the general case, $1 < p <+\infty$}\label{sec:gp_sobolev_LP}
We can now state our first result, which deals with $W^{m,p}(\mathcal{D})$-regularity of Gaussian processes, given any $p\in(1,+\infty)$ and any open set $\mathcal{D}\subset \mathbb{R}^d$.
\begin{proposition}[Sample path Banach-Sobolev regularity for Gaussian processes]\label{prop:sobolev_reg_sto_LP}
Let $\mathcal{D} \subset \mathbb{R}^d$ be an open set. Let $(U(x))_{x\in \mathcal{D}} \sim GP(0,k)$ be a measurable centered Gaussian process, defined on a probability set $(\Omega, \mathcal{F},\mathbb{P})$, such that its standard deviation function $\sigma$ lies in $L_{loc}^1(\mathcal{D})$. Let $p\in(1,+\infty)$. The following statements are equivalent :
\begin{enumerate}[label=(\roman*),wide,labelwidth=0em,labelindent=0pt]
\item (Sample path regularity) The sample paths of $(U(x))_{x\in \mathcal{D}}$ lie in $W^{m,p}(\mathcal{D})$ almost surely.
\item \label{pt:integral_crit_LP} (Integral criteria) For all $|\alpha|\leq m$, the distributional derivative $\partial^{\alpha,\alpha}k$  lies in $L^p(\mathcal{D} \times \mathcal{D})$ and admits a representer $k_{\alpha}$ in $L^p(\mathcal{D} \times \mathcal{D})$ which is the covariance function of a measurable Gaussian process. Note $\sigma_{\alpha}(x) := k_{\alpha}(x,x)^{1/2}$, then additionally
\begin{align}
    \int_{\mathcal{D}}\sigma_{\alpha}(x)^pdx < +\infty
\end{align}
\item \label{pt:cov_strct}(Covariance structure) For all $|\alpha|\leq m$, the distributional derivative $\partial^{\alpha,\alpha}k$  lies in $L^p(\mathcal{D} \times \mathcal{D})$ and the associated integral operator $\mathcal{E}_k^{\alpha} : L^q(\mathcal{D}) \rightarrow L^p(\mathcal{D})$ defined by
\begin{align}\label{eq:def_e_k_alpha_LP}
\mathcal{E}_k^{\alpha}f(x) = \int_{\mathcal{D}} \partial^{\alpha,\alpha}k(x,y)f(y)dy
\end{align}
is symmetric, positive and nuclear: there exists $(\lambda_n^{\alpha})_{n\in\mathbb{N}} \geq 0$ and $(\psi_n^{\alpha})_{n \in\mathbb{N}} \subset L^p(\mathcal{D})$ such that
\begin{numcases}{}
    \sum_{n=0}^{+\infty}\lambda_n^{\alpha}||\psi_n^{\alpha}||_{L^p(\mathcal{D})}^2 < + \infty \\
\partial^{\alpha,\alpha}k(x,y) = \sum_{n=0}^{+\infty}\lambda_n^{\alpha}\psi_n^{\alpha}(x)\psi_n^{\alpha}(y) \ \ \text{ in } \ \ L^p(\mathcal{D}\times\mathcal{D}) \label{eq:nuclear_rep}
\end{numcases}
If $1 \leq p \leq 2$, then one can choose $(\lambda_n^{\alpha})$ such that $\sum_n \lambda_n^{\alpha} < + \infty$, and there exists a bounded operator $A_{\alpha} : L^2(\mathcal{D})\rightarrow L^p(\mathcal{D})$ and an orthonormal basis $(\phi_n^{\alpha})$ of $L^2(\mathcal{D})$ such that $\psi_n^{\alpha} = A_{\alpha}\phi_n^{\alpha}$ for all $n\geq 0$ (in particular, we have the uniform bound $||\psi_n^{\alpha}||_p \leq ||A_{\alpha}||$).
\end{enumerate}
\end{proposition}
The proposition above shows that a suitable $L^p$ control of the function $\partial^{\alpha,\alpha}k$ over the diagonal is necessary and sufficient for ensuring the Sobolev regularity of the sample paths of the Gaussian process with covariance function $k$. Formally speaking, the function $(x,y) \mapsto\partial^{\alpha,\alpha}k(x,y)$ is the covariance function of the differentiated process, $(\omega,x) \mapsto \partial^{\alpha}U_{\omega}(x)$. This is formal only, as the weak derivative of the sample paths are only defined up to a set of Lebesgue measure zero, and thus there is no obvious way of defining the joint map $(\omega,x) \mapsto \partial^{\alpha}U_{\omega}(x)$. Note also that the idea of ensuring a suitable control of this covariance function near its diagonal is not with reminding more standard results pertaining to the differentiability in the \textit{mean square} sense of a random process (see e.g. \cite{adler2007}, Section 1.4.2). See \cite{SCHEUERER2010} for similar remarks on the Sobolev regularity of random fields.
\begin{remark}\label{rk:asymmetry}
In Proposition \ref{prop:sobolev_reg_sto_LP}, there is an asymmetry between Point $\ref{pt:integral_crit_LP}$ and Point $\ref{pt:cov_strct}$: one depends on whether $p$ is lower or greater than $2$ while the other does not. Moreover, Points $(ii)$ and $(iii)$ rely on the finiteness of some quantity, so explicit bounds should be sought so that Point $\ref{pt:integral_crit_LP}$ controls Point $\ref{pt:cov_strct}$ and conversely. This is the content of Proposition \ref{prop_LP_bounds}.
\end{remark}

\begin{remark}
Under the assumption that $(U(x))_{x\in\mathcal{D}}$ is measurable, the statement that its sample paths lie in some Sobolev space is \textit{not} up to a modification of the process. This is a consequence of Lemmas \ref{lemma:E_et_F_countable_sobolev}, \ref{lemma:u_alpha_phi_GP} and \ref{lemma:expo_moments}, which show that the Sobolev regularity of its paths is fully determined by the finite dimensional marginals of the process (see equation \eqref{eq:sup_u_al_phi}). This contrasts with more classical results, e.g. pertaining to the continuity of the process (\cite{azais_level_2009}, Section 1.4.1). Still, ensuring the measurability of the process is not really straightforward (see Remark \ref{rk:GP_mes}).
\end{remark}

\begin{example}[Finite rank covariance functions]\label{ex:finite_rank} Let $p\in(1,+\infty)$, $m\in\mathbb{N}_0$ and $N\in\mathbb{N}$. Consider $f_1,...,f_N\in W^{m,p}(\mathcal{D})$ and choose once and for all representers of those functions in $L^p(\mathcal{D})$, also denoted by $f_1,...,f_N$, so that they may be understood as functions in the classical sense. Consider the covariance function $k(x,x') := \sum_{i=1}^N f_i(x)f_i(x')$.
Then obviously, for all $|\alpha|\leq m$, the weak derivative $\partial^{\alpha,\alpha}k$ is given by
\begin{align}
\partial^{\alpha,\alpha}k(x,x') = \sum_{i=1}^N \partial^{\alpha}f_i(x)\partial^{\alpha}f_i(x') \ \ \ \text{in}\  L^p(\mathcal{D}\times\mathcal{D})
\end{align}
and the associated integral operators fulfill the criteria $\ref{pt:cov_strct}$ of Proposition \ref{prop:sobolev_reg_sto_LP}. Thus the corresponding measurable Gaussian process has its sample paths in $W^{m,p}(\mathcal{D})$ almost surely. Note that this was obvious in the first place, since this Gaussian process can be written as $U(x) = \sum_{i=1}^N\xi_if_i(x)$ where $(\xi_1,...,\xi_N)$ are independent standard Gaussian random variables (checking that the covariance function is the right one is trivial). Still, this example fell out of the scope of the previous results pertaining to the Sobolev regularity of Gaussian processes.
\end{example}
\begin{proof}(Proposition \ref{prop:sobolev_reg_sto_LP}) We show $(i) \implies (ii)\ \& \ (iii)$, $(ii) \implies (i)$ and $(iii) \implies (ii)$.\\
\underline{$(i) \implies (ii)\ \& \ (iii):$} Suppose $(i)$ and let $|\alpha|\leq m$. We first prove that the map $N_{\alpha} : (\Omega,\mathcal{F},\mathbb{P}) \rightarrow (\mathbb{R},\mathcal{B}(\mathbb{R})),\  \omega \mapsto ||\partial^{\alpha}U_{\omega}||_{L^p(\mathcal{D})}$ is measurable. Indeed, given $\varphi \in F_q$, the map
\begin{align}
    U_{\varphi}^{\alpha} : \omega \longmapsto\int_{\mathcal{D}} \partial^{\alpha}U_{\omega}(x)\varphi(x)dx  = (-1)^{|\alpha|}\int_{\mathcal{D}} U_{\omega}(x)\partial^{\alpha}\varphi(x)dx 
\end{align}
is a real valued random variable. This follows from the fact that $U_{\omega} \in L_{loc}^1(\mathcal{D})$, the measurability of $U$ and Lemma \ref{lemma:u_alpha_phi_GP}. Note now that from Lemma \ref{lemma:E_et_F_countable_sobolev},
\begin{align}\label{eq:sup_u_al_phi}
    \bigg(\omega \mapsto ||\partial^{\alpha}U_{\omega}||_{L^p(\mathcal{D})}\bigg) = \sup_{\varphi \in F_q} |U_{\varphi}^{\alpha}|
\end{align}
The supremum being taken over a countable set, $N_{\alpha}$ is indeed measurable. Given any $f \in L^p(\mathcal{D})$, a slight modification of this proof shows that $\omega \mapsto ||\partial^{\alpha}U_{\omega} - f||_{L^p(\mathcal{D})}$ is also measurable.
We can now show the map $T_{\alpha} : (\Omega,\mathcal{F},\mathbb{P}) \rightarrow (L^p(\mathcal{D}),\mathcal{B}(L^p(\mathcal{D}))),\  \omega \mapsto \partial^{\alpha}U_{\omega}$ is measurable. Let $f\in L^p(\mathcal{D}), r > 0$ and $B = B(f,r)$ be an open ball in $L^p(\mathcal{D})$. Then from the measurability of $\omega \mapsto ||\partial^{\alpha}U_{\omega} - f||_{L^p(\mathcal{D})}$,
\begin{align}
    T_{\alpha}^{-1}(B) = \{\omega\in \Omega : ||\partial^{\alpha}u-f||_{L^p(\mathcal{D})} < r \} \in \mathcal{F}
\end{align}
Since $L^p(\mathcal{D})$ is a separable metric space, its Borel $\sigma$-algebra is generated by the open balls of $L^p(\mathcal{D})$ (see e.g. \cite{bogachev2007measure}, Exercise 6.10.28). Thus $T_{\alpha}$ is Borel-measurable and the pushforward of $\mathbb{P}$ through $T_{\alpha}$ induces a (centered) probability measure $\mu_{\alpha}$ over the Banach space $L^p(\mathcal{D})$. We show that it is Gaussian. Let $f \in L^q(\mathcal{D})$ and denote $T_f$ the associated linear form over $L^p(\mathcal{D})$. Let $(\phi_n) \subset C_c^{\infty}(\mathcal{D})$ be such that $\phi_n \rightarrow f$ in $L^q(\mathcal{D})$ and $\omega \in\Omega$ be such that $U_{\omega}$ lies in $L_{loc}^1(\mathcal{D})$:
\begin{align}
    T_f(\partial^{\alpha}U_{\omega}) = \int_{\mathcal{D}}\partial^{\alpha}U_{\omega}(x)f(x)dx &= \lim_{n\rightarrow \infty} \int_{\mathcal{D}}\partial^{\alpha}U_{\omega}(x)\phi_n(x)dx \label{eq:weak_GRV}\\
    &= \lim_{n\rightarrow \infty}(-1)^{|\alpha|} \int_{\mathcal{D}}U_{\omega}(x)\partial^{\alpha}\phi_n(x)dx
\end{align}
For each value of $n$, Lemma \ref{lemma:u_alpha_phi_GP} shows that the map $\omega \mapsto (-1)^{|\alpha|} \int_{\mathcal{D}}U_{\omega}(x)\partial^{\alpha}\phi_n(x)dx$ is a Gaussian random variable. Thus $\omega \mapsto T_f(\partial^{\alpha}U_{\omega})$ is a Gaussian random variable as an a.s. limit of Gaussian random variables. This shows that the pushforward of $\mu_{\alpha}$ through $T_f$ is Gaussian, since for all Borel set $B\in\mathcal{B}(\mathbb{R})$,
\begin{align}
    \mu_{\alpha}(T_f^{-1}(B)) &= \mu_{\alpha}(\{g\in L^p(\mathcal{D}):  T_f(g)\in B \}) = \mathbb{P}(\{\omega \in \Omega: T_f(\partial^{\alpha}U_{\omega})\in B\})
\end{align}
Hence, $\mu_{\alpha}$ is Gaussian.
We next show that $\partial^{\alpha,\alpha}k \in L^p(\mathcal{D}\times\mathcal{D})$ and that the covariance operator of $\mu_{\alpha}$ is the integral operator $\mathcal{E}_k^{\alpha} : L^q(\mathcal{D})\rightarrow L^p(\mathcal{D})$ with kernel $\partial^{\alpha,\alpha}k$. Let $\mathcal{D}_0 \Subset \mathcal{D}\times \mathcal{D}$ and $K_0 \Subset \mathcal{D}$ such that $\mathcal{D}_0 \subset K_0\times K_0$ (for example, set $K_1 := \overline{\{x\in \mathcal{D} : \exists y \in \mathcal{D}, (x,y) \in K \}}$, $K_2 := \overline{\{y\in \mathcal{D} : \exists x \in \mathcal{D}, (x,y) \in K \}}$ which are both compact subsets of $\mathcal{D}$ and $K_0 := K_1 \cup K_2$). Let $h=(h_1,...,h_d) \in(\mathbb{R}_+^*)^d$ be such that $\sum_i|h_i| < \text{dist}(K_0,\mathcal{D}_0)$. Use then the bilinearity of the covariance operator:
\begin{align}
    \int_{\mathcal{D}_0}|(\delta_h^{\alpha}\otimes \delta_h^{\alpha}) k(x,y)|^pdxdy &= \int_{\mathcal{D}_0} |\mathbb{E}[\delta^{\alpha}_hU(x)\delta^{\alpha}_hU(y)]|^p dxdy\label{eq:K_0_D_0} \\ 
    &\leq \int_{K_0\times K_0} |\mathbb{E}[\delta^{\alpha}_hU(x)\delta^{\alpha}_hU(y)]|^p dxdy \\
    &\leq \int_{K_0\times K_0} \mathbb{E}[|\delta^{\alpha}_hU(x)\delta^{\alpha}_hU(y)|^p] dxdy \\
    &\leq \mathbb{E}\bigg[\bigg(\int_{K_0}|\delta^{\alpha}_hU(x)|^pdx\bigg)^2\bigg] = \mathbb{E}[||\delta^{\alpha}_hU||_p^{p/2}]\\
    &\leq \mathbb{E}[||U||_{W^{m,p}(\mathcal{D})}^{p/2}] =: C^p < +\infty \label{eq:expectation_sobolev_norm}
\end{align}
The expectation in equation \eqref{eq:expectation_sobolev_norm} is indeed finite because of the following. Given $|\alpha|\leq m$, equation \eqref{eq:sup_u_al_phi} shows that the map $\omega \mapsto ||\partial^{\alpha}U_{\omega}||_p$ is the supremum of a Gaussian sequence which is finite a.s. by assumption; Lemma \ref{lemma:expo_moments} then implies that all the moments of this supremum are finite. Writing then $||U||_{W^{m,p}}$ in terms of these $L^p$ norms yields equation \eqref{eq:expectation_sobolev_norm}. To see that the control \eqref{eq:expectation_sobolev_norm} implies that $\partial^{\alpha,\alpha}k \in L^p(\mathcal{D}\times\mathcal{D})$, we copy the steps of equations \eqref{eq:delta_alpha_control}-\eqref{eq:discrete_ibp}-\eqref{eq:discrete_ibp_control} in the proof of Lemma \ref{lemma:linear_form_sobolev}. Let $\varphi\in C_c^{\infty}(\mathcal{D}\times\mathcal{D})$. Since it is compactly supported in $\mathcal{D}\times\mathcal{D}$, find an open set $\mathcal{D}_0 \Subset \mathcal{D}$ such that $\text{Supp}(\varphi) \subset\mathcal{D}_0$. Use Hölder's inequality and equation \eqref{eq:expectation_sobolev_norm}:
\begin{align}\label{eq:delta_alpha_otimes_k}
    \bigg|\int_{\mathcal{D}\times\mathcal{D}}(\delta_h^{\alpha}\otimes \delta_h^{\alpha}) k(x,y)\varphi(x,y)dxdy\bigg| \leq  ||(\delta_h^{\alpha}\otimes \delta_h^{\alpha}) k||_p||\varphi||_q \leq C||\varphi||_q
\end{align}
Next, use the discrete integration by parts formula:
\begin{align}
    \int_{\mathcal{D}\times\mathcal{D}}(\delta_h^{\alpha}\otimes \delta_h^{\alpha}) k(x,y)\varphi(x,y)dxdy = \int_{\mathcal{D}} k(x,y)(\delta_h^{\alpha}\otimes \delta_h^{\alpha})^*\varphi(x,y)dxdy
\end{align}
When $h\rightarrow 0$, observe that $(\delta_h^{\alpha}\otimes \delta_h^{\alpha})^*\varphi(x,y) \rightarrow \partial^{\alpha,\alpha}\varphi(x,y)$ pointwise. Use Lebesgue's dominated convergence theorem and equation \eqref{eq:delta_alpha_otimes_k} to obtain
\begin{align}
    \bigg|\int_{\mathcal{D}\times\mathcal{D}}k(x,y)\partial^{\alpha,\alpha}\varphi(x,y)dxdy\bigg| \leq C||\varphi||_q
\end{align}
which indeed shows that $\partial^{\alpha,\alpha}k\in L^p(\mathcal{D}\times\mathcal{D})$, from Riesz' lemma. We now identify $K_{\alpha}$, the covariance operator of $\mu_{\alpha}$, in terms of $\partial^{\alpha,\alpha}k$. Let $f,g\in L^q(\mathcal{D})$ and using the density of $ C_c^{\infty}(\mathcal{D})$ in $L^q(\mathcal{D})$ (\cite{fournier_adams_sobolev}, Corollary 2.30), let $(f_n),(g_n) \subset C_c^{\infty}(\mathcal{D})$ be two sequences such that $f_n \rightarrow f$ in $L^q(\mathcal{D})$ and likewise for $g_n$ and $g$. Then (explanation below),
\begin{align}
    \langle f,K_{\alpha}g\rangle_{L^q,L^p} &= \lim_{n\rightarrow \infty}\langle f_n,K_{\alpha}g_n\rangle_{L^q,L^p} \label{eq:seq_cont}\\
    &=\lim_{n\rightarrow \infty}\int_{L^p(\mathcal{D})}\langle f_n,h\rangle_{L^q,L^p} \langle g_n,h\rangle_{L^q,L^p} d\mu_{\alpha}(h) \nonumber\\
    &= \lim_{n\rightarrow \infty}\int_{\Omega} \langle f_n,\partial^{\alpha}U_{\omega}\rangle_{L^q,L^p} \langle g_n,\partial^{\alpha}U_{\omega}\rangle_{L^q,L^p}d\mathbb{P}(\omega)\label{eq:pushf_K_alpha} \\
    &= \lim_{n\rightarrow \infty}\int_{\Omega} \langle \partial^{\alpha}f_n,U_{\omega}\rangle_{L^q,L^p} \langle \partial^{\alpha}g_n,U_{\omega}\rangle_{L^q,L^p}d\mathbb{P}(\omega) \nonumber\\
    &= \lim_{n\rightarrow \infty}\int_{\mathcal{D}\times\mathcal{D}}\partial^{\alpha}f_n(x)\partial^{\alpha}g_n(y)k(x,y)dxdy\label{eq:fubi_K_alpha}\\
    &= \lim_{n\rightarrow \infty}\int_{\mathcal{D}\times\mathcal{D}}f_n(x)g_n(y)\partial^{\alpha,\alpha}k(x,y)dxdy\nonumber\\
    &= \int_{\mathcal{D}\times\mathcal{D}}f(x)g(y)\partial^{\alpha,\alpha}k(x,y)dxdy = \langle f,\mathcal{E}_k^{\alpha}g\rangle_{L^q,L^p}
\end{align}
We used the sequential continuity of $K_{\alpha}$ in equation \eqref{eq:seq_cont}, the transfer theorem for pushforward measure integration (\cite{bogachev2007measure}, Theorem 3.6.1) in equation \eqref{eq:pushf_K_alpha} and Fubini's theorem in equation \eqref{eq:fubi_K_alpha}.
According to Proposition \ref{prop_bogachev_LP}, since $\mu_{\alpha}$ is a Gaussian measure over $L^p(\mathcal{D})$, there exists a representer $k_{\alpha}$ of $\partial^{\alpha,\alpha}k$ in $L^p(\mathcal{D}\times\mathcal{D})$ which is the covariance function of a measurable Gaussian process. Note $\sigma_{\alpha}(x) = k_{\alpha}(x,x)^{1/2}$. Then the same proposition shows that 
\begin{align}
\int_{\mathcal{D}} \sigma_{\alpha}(x)^pdx < + \infty
\end{align}
which shows $(ii)$. By Proposition 3.5.11 from \cite{bogachev1998gaussian}, $\mathcal{E}_k^{\alpha}$ is nuclear and admits a symmetric nonnegative representation as the one in equation \eqref{eq:nuclear_rep}. if $1\leq p \leq 2$,  then $L^p(\mathcal{D})$ is of cotype $2$ and since $\mathcal{E}_k^{\alpha}$ is a Gaussian covariance operator, from Proposition \ref{prop_cotype2} there exists a bounded operator $A_{\alpha} : L^2(\mathcal{D}) \rightarrow L^p(\mathcal{D})$ and a trace class operator $S_{\alpha} : L^2(\mathcal{D}) \rightarrow L^2(\mathcal{D})$ such that $\mathcal{E}_k^{\alpha} = A_{\alpha}S_{\alpha}A_{\alpha}^*$ ($l^2$ and $L^2(\mathcal{D})$ are isomorphic Hilbert spaces). Introduce a Mercer decomposition of $S_{\alpha}$ (equation \eqref{eq:pos_sym_decomp}): $S_{\alpha} = \sum_n \lambda_n^{\alpha}\phi_n^{\alpha} \otimes \phi_n^{\alpha}$. Use the continuity of $A_{\alpha}$ and $A_{\alpha}^*$ to obtain that $\partial^{\alpha,\alpha}k(x,y) = \sum_n \lambda_n^{\alpha}(A_{\alpha}\phi_n^{\alpha})(x) (A_{\alpha}\phi_n^{\alpha})(y)$ in  $L^p(\mathcal{D}\times\mathcal{D})$, which yields $(iii)$. \\
\underline{$(ii) \implies (i) :$} from Proposition \ref{prop_bogachev_LP}, let $(V^{\alpha})$ be a centered measurable Gaussian process with covariance function $k_{\alpha}$. Then its sample paths lie in $L^p(\mathcal{D})$ a.s. and the Gaussian measure it induces over $L^p(\mathcal{D})$ through the map $\omega \mapsto V_{\omega}^{\alpha} \in L^p(\mathcal{D})$ is the centered Gaussian measure with covariance operator $\mathcal{E}_k^{\alpha}$. Given $\varphi\in C_c^{\infty}(\mathcal{D})$, denote $V_{\varphi}^{\alpha}$ the following random variable
\begin{align}
\omega \mapsto \int_{\mathcal{D}} V^{\alpha}_{\omega}(x)\varphi(x)dx
\end{align}
From Lemma \ref{lemma:u_alpha_phi_GP}, $(V_{\varphi}^{\alpha})_{\varphi \in F_q}$ is a Gaussian sequence. It is also centered and using Fubini's theorem to permute $\mathbb{E}$ and $\int$, we have that
\begin{align}
    \mathbb{E}[V_{\varphi}^{\alpha}V_{\psi}^{\alpha}] &= \int_{\mathcal{D}\times\mathcal{D}} \varphi(y)\psi(x)k_{\alpha}(x,y)dxdy\nonumber = \int_{\mathcal{D}\times\mathcal{D}} \varphi(y)\psi(x)\partial^{\alpha,\alpha} k(x,y)dxdy\nonumber\allowdisplaybreaks \\
    &= \int_{\mathcal{D}\times\mathcal{D}} \partial^{\alpha}\varphi(y)\partial^{\alpha}\psi(x) k(x,y)dxdy \\
    \mathbb{E}[U_{\varphi}^{\alpha}U_{\psi}^{\alpha}] &= \int_{\mathcal{D}\times\mathcal{D}} \partial^{\alpha}\varphi(y)\partial^{\alpha}\psi(x) k(x,y)dxdy
\end{align}
Having the same mean and covariance, the two Gaussian sequences $(V_{\varphi}^{\alpha})_{\varphi \in F_q}$ and $(U_{\varphi}^{\alpha})_{\varphi \in F_q}$ have the same finite dimensional marginals. One checks in an elementary fashion that their countable supremums over $F_q$ then have the same law (e.g. by showing that they have the same cumulative distribution function). Recalling from Lemma \ref{lemma:E_et_F_countable_sobolev} that $||V^{\alpha}_{\omega}||_p = \sup_{\varphi \in F_q} |V^{\alpha}_{\varphi}(\omega)|$,
we obtain that
\begin{align}
    1 = \mathbb{P}(||V^{\alpha}_{\omega}||_p < + \infty) = \mathbb{P}(\sup_{\varphi \in F_q} |V^{\alpha}_{\varphi}| < + \infty) = \mathbb{P}(\sup_{\varphi \in F_q} |U^{\alpha}_{\varphi}| < + \infty)
\end{align}
which shows that $\partial^{\alpha}U \in L^p(\mathcal{D})$ almost surely. This is true for all $|\alpha| \leq m$, which shows $(i)$. \\
\underline{$(iii) \implies (ii) :$ if $(iii)$}, then from either Proposition \ref{prop_type2} or \ref{prop_cotype2} depending on whether $p\leq 2$ or $p\geq 2$, there exists a Gaussian measure over $L^p(\mathcal{D})$ whose covariance operator is $\mathcal{E}_k^{\alpha}$ as defined in equation \eqref{eq:def_e_k_alpha_LP}. Proposition \ref{prop_bogachev_LP} yields $(ii)$.
\end{proof}
The following proposition deals with the issues raised in Remark \ref{rk:asymmetry} (asymmetry between Points $(ii)$ and $(iii)$ of Proposition \ref{prop:sobolev_reg_sto_LP}). We recall that the nuclear norm $\nu(T)$ of a nuclear operator $T$ is defined in equation \eqref{eq:def_nuclear_norm}. Contrarily to Proposition \ref{prop:sobolev_reg_sto_LP}, we do not exclude $p=1$.
\begin{proposition}\label{prop_LP_bounds}
Let $\mu$ be a centered Gaussian measure over $L^p(\mathcal{D})$, where $1 \leq p < +\infty$. Let $k\in L^p(\mathcal{D}\times\mathcal{D})$ be the kernel of its covariance operator $(K_{\mu}=\mathcal{E}_k)$, chosen such that $k$ is also the covariance function of the measurable Gaussian process provided by Proposition \ref{prop_bogachev_LP}. Define $\sigma(x) = k(x,x)^{1/2}$ and set $C_p = 2^{p/2}\Gamma(\frac{p+1}{2})/\sqrt{\pi}\ (= \mathbb{E}[|X|^p]$ where $X\sim\mathcal{N}(0,1) )$. Then the following holds.
\begin{itemize}[wide,labelwidth=0em,labelindent=0pt]
    \item if $2 \leq p < +\infty$, then $\mathcal{E}_k$ is nuclear and
    \begin{align}
     C_p^{-\frac{2}{p}} \nu(\mathcal{E}_k) \leq ||\sigma||_p^2 \leq \nu(\mathcal{E}_k)
    \end{align}
    \item if $1 \leq p \leq 2$, there exists a nuclear operator $S$ over $L^2(\mathcal{D})$ and a bounded operator $A : L^2(\mathcal{D})\rightarrow L^p(\mathcal{D})$ such that $\mathcal{E}_k = ASA^*$. Moreover,
    \begin{align}
    \nu(\mathcal{E}_k) \leq  \inf_{\substack{A, S \ s.t. \\ \mathcal{E}_k = ASA^*}}||A||^2\nu(S) \leq ||\sigma||_p^2 \leq C_p^{-\frac{2}{p}} \inf_{\substack{A, S \ s.t. \\ \mathcal{E}_k = ASA^*}}||A||^2\nu(S) \label{eq:control_p_leq_222}
\end{align}
\item if $p=2$ then $||\sigma||_2^2 = \nu(\mathcal{E}_k) = \text{Tr}(\mathcal{E}_k)$.
\end{itemize}
\end{proposition}
It is expected that the nuclear norm of $\mathcal{E}_k$ cannot directly appear on the right hand side of equation \eqref{eq:control_p_leq_222}, as not all nuclear operators are Gaussian covariance operators when $1\leq p < 2$. This proposition in fact suggests that for general Banach spaces $X$ of cotype $2$, the following map defined over the set of Gaussian covariance operators $B : X^* \rightarrow X$,
\begin{align}
    B \mapsto \inf_{\substack{A, S \ s.t. \\ B = ASA^*}}||A||^2\nu(S)
\end{align}
is the natural measurement of the "size" of such operators. When $X$ is of type $2$, this would be the case for the nuclear norm $B \mapsto \nu(B)$. 
\begin{remark}\label{rk:control_sobolev_norm}
Proposition \ref{prop_LP_bounds} is interesting from an application point of view because it strongly suggests that the operator norms appearing in this proposition, as well as the $L^p$ norm of the standard deviation function $\sigma$, are the \textit{correct} quantities for quantitatively controlling the $L^p$ norm of the sample paths of the underlying Gaussian process. For instance, we have the following $L^p$ control in expectation: $\mathbb{E}[||U||_p^p] = C_p||\sigma||_p^p$ (see equation \eqref{eq:LP_norm_GP_sigma}). Applying this fact recursively, we obtain that the $W^{m,p}$-Sobolev norm of the sample paths of the Gaussian process in question is controlled as follow, denoting $\sigma_{\alpha}(x) = \partial^{\alpha,\alpha}k(x,x)^{1/2}$ (temporarily discarding definition and measurability issues w.r.t. $\sigma_{\alpha}$)
\begin{align}\label{eq:control_sobolev_norm}
    \mathbb{E}\big[||U||_{W^{m,p}}^p\big] = C_p\sum_{|\alpha|\leq m} ||\sigma_{\alpha}||_p^p
\end{align}
If such a control cannot be obtained, then it means that the sample paths of $U$ do not lie in $W^{m,p}(\mathcal{D})$ in the first place.
Additional "sharp" controls can then be obtained from equation \eqref{eq:control_sobolev_norm} using Proposition \ref{prop_LP_bounds} and other elementary inequalities involving the expectation. Finally, we have the following asymptotic behaviour of the constant when $p \rightarrow +\infty: C_p^{-2/p} \sim \exp(1)/({p-1})$.
\end{remark}

\begin{proof}(Proposition \ref{prop_LP_bounds})
Suppose first that $p\geq 2$. Let $k = \sum_n \mu_n \psi_n\otimes\phi_n, \mu_n \geq 0$, be a nuclear representation of $k$ (rather, $\mathcal{E}_k$), with $||\psi_n||_p = ||\phi_n||_p = 1$ and $S := \sum_n |\mu_n| < + \infty$. Then, using the discrete Jensen's inequality on the weights $|\mu_n|/S$ and the function $x \mapsto |x|^{p/2}$ ($p/2\geq1$)
\begin{align}
    ||\sigma||_p^p &= \int_{\mathcal{D}} \sigma(x)^pdx = \int_{\mathcal{D}} \bigg(\sum_{n=0}^{+\infty}\mu_n \psi_n(x)\phi_n(x)\bigg)^{p/2}dx \\
    &= S^{p/2}\int_{\mathcal{D}} \bigg(\sum_{n=0}^{+\infty}\frac{|\mu_n|}{S}\frac{\mu_n}{|\mu_n|} \psi_n(x)\phi_n(x)\bigg)^{p/2}dx \\
    &\leq S^{p/2}\int_{\mathcal{D}} \sum_{n=0}^{+\infty}\frac{|\mu_n|}{S}\bigg|\frac{\mu_n}{|\mu_n|}\bigg|^{p/2} |\psi_n(x)|^{p/2}|\phi_n(x)|^{p/2}dx = S^{p/2-1} \sum_{n=0}^{+\infty}|\mu_n|\times ||\psi_n\phi_n||_{p/2}^{p/2}\\
    &\leq S^{p/2-1} \sum_{n=0}^{+\infty}|\mu_n|\times||\psi_n||_p^{p/2}||\phi_n||_p^{p/2} = S^{p/2-1}\sum_{n=0}^{+\infty}|\mu_n| = S^{p/2}\label{eq:jensen1}
\end{align}
We used the Cauchy-Schwarz inequality on $||\psi_n\phi_n||_{p/2}^{p/2}$ in equation \eqref{eq:jensen1}. Since equation \eqref{eq:jensen1} holds whatever the nuclear decomposition of $k$, taking the infimum over $S$ in equation \eqref{eq:jensen1} yields $||\sigma||_p \leq \sqrt{\nu(\mathcal{E}_k)}$.
Conversely, consider the measurable Gaussian process $(U(x))_{x\in\mathcal{D}}$ provided by Proposition \ref{prop_bogachev_LP} with covariance function $k$. Fubini's theorem yields
\begin{align}\label{eq:LP_norm_GP_sigma}
    \mathbb{E}[||U||_p^p] = \mathbb{E}\bigg[\int_{\mathcal{D}}|U(x)|^pdx\bigg] = \int_{\mathcal{D}}\mathbb{E}[|U(x)|^p]dx = \int_{\mathcal{D}}C_p\sigma(x)^pdx = C_p||\sigma||_p^p
\end{align}
where $C_p = 2^{p/2}\Gamma(\frac{p+1}{2})/\sqrt{\pi}$. Indeed, given $X \sim \mathcal{N}(0,\sigma^2)$, then $\mathbb{E}[|X|^p] = C_p\sigma^p$. Moreover, introduce $\mu = \mathbb{P}_U$ the Gaussian measure over $L^p(\mathcal{D})$ induced by $U$, whose covariance operator is $\mathcal{E}_k$ (see Proposition \ref{prop_bogachev_LP}). We successively use the transfer theorem for pushforward measure integration, Jensen's inequality for probability measures $(p/2>1)$ and the nuclear norm estimate from \cite{linde1980characterization}, Theorem 3:
\begin{align}
    \mathbb{E}[||U||_p^p] &= \int_{\Omega} ||U_{\omega}||_p^p\mathbb{P}(d\omega) = \int_{L^p(\mathcal{D})}||f||_p^p\mu(df) = \int_{L^p(\mathcal{D})}||f||_p^{2\times p/2}\mu(df)\\
    &\geq \bigg(\int_{L^p(\mathcal{D})}||f||_p^{2}\mu(df)\bigg)^{p/2} \geq \nu(\mathcal{E}_k)^{p/2}
\end{align}
To conclude, when $2\leq p< + \infty$,
\begin{align}
    C_p^{-\frac{1}{p}} \sqrt{\nu(\mathcal{E}_k)} \leq ||\sigma||_p \leq \sqrt{\nu(\mathcal{E}_k)} \label{eq:control_p_geq_2}
\end{align}
Suppose now that $1 \leq p < 2$. Let $\mu_0$ be a Gaussian measure on $L^2(\mathcal{D})$ and $A : L^2(\mathcal{D}) \rightarrow L^p(\mathcal{D})$ a bounded operator such that $\mu = {\mu_0}_A$ (pushforward of $\mu_0$ through $A$) and $S$ the trace class covariance operator associated to $\mu_0$ (see Proposition \ref{prop_cotype2}). $(U(x))_{x\in\mathcal{D}}$ remains the Gaussian process of Proposition \ref{prop_bogachev_LP} and we have $\mu = \mathbb{P}_U$. Then (explanation below),
\begin{align}
    C_p||\sigma||_p^p &= \mathbb{E}[||U||_p^p] =  \int_{\Omega} ||U_{\omega}||_p^p\mathbb{P}(d\omega) = \int_{L^p(\mathcal{D})}||f||_p^p\mu(df) \label{eq:fubi_pushf}\allowdisplaybreaks \\
    &= \int_{L^2(\mathcal{D})}||Ag||_p^p\mu_0(dg) \leq ||A||^p\int_{L^2(\mathcal{D})}||g||_2^p\mu_0(dg) \label{eq:push_mu_mu0} \\
    &\leq ||A||^p\int_{L^2(\mathcal{D})}\langle g,g\rangle_{L^2}^{p/2}\mu_0(dg) \leq ||A||^p\bigg(\int_{L^2(\mathcal{D})}\langle g,g\rangle_{L^2}\mu_0(dg)\bigg)^{p/2}\label{eq:jensen_mu0} \\
    &\leq ||A||^p\ \text{Tr}(S)^{p/2} = ||A||^p \nu(S)^{p/2} \label{eq:nu_S_L2}
\end{align}
In equation \eqref{eq:fubi_pushf}, we used equation \eqref{eq:LP_norm_GP_sigma} and pushforward integration to write the integral w.r.t. $\mathbb{P}$ as an integral w.r.t. $\mu =\mathbb{P}_U$.
Likewise in equation \eqref{eq:push_mu_mu0} where we write the integral w.r.t. $\mu$ as an integral w.r.t. $\mu_0$ using the pushforward identity $\mu = {\mu_0}_A$. In equation \eqref{eq:jensen_mu0}, we used Jensen's inequality for concave functions ($0 < p/2 < 1$). In equation \eqref{eq:nu_S_L2}, we used the trace identity from \cite{bogachev1998gaussian}, equation 2.3.2 and the one following p. 49. Equation \eqref{eq:nu_S_L2} then yields $||\sigma||_p \leq C_p^{-1/{p}} ||A||\sqrt{\nu(S)}$. In the last equation, taking the infimum over all representations $\mathcal{E}_k = ASA^*$ yields:
\begin{align}
    ||\sigma||_p \leq C_p^{-\frac{1}{p}} \inf_{\substack{A, S \ s.t. \\ \mathcal{E}_k = ASA^*}}||A||\sqrt{\nu(S)}
\end{align}
To prove the other inequality, we use an explicit decomposition $\mathcal{E}_k = ASA^*$ by first setting
\begin{align}
    Af(x) = f(x)\sigma(x)^{1-p/2}
\end{align}
Use Hölder's inequality with $a = 2/p, 1/a + 1/b = 1$ (notice that $a > 1$)
\begin{align}\label{eq:A_bound}
    ||Af||_{p}^p &= \int_{\mathcal{D}} |f(x)|^p\sigma(x)^{p(1-p/2)}dx  \\
    &\leq \bigg(\int_{\mathcal{D}}|f(x)|^2dx\bigg)^{p/2}\bigg(\int_{\mathcal{D}}\sigma(x)^{bp(1-p/2)}dx\bigg)^{1/b}
\end{align}
But $b = \frac{a}{a-1} = \frac{2/p}{2/p-1} = \frac{1}{1-p/2}$ and $b(1-p/2) =  1$, which together with equation \eqref{eq:A_bound} yields
\begin{align}
    ||Af||_{p}^p \leq ||f||_{2}^p||\sigma||_{p}^{p(1-p/2)}
\end{align}
Thus $A : L^2(\mathcal{D})\rightarrow L^p(\mathcal{D})$ is bounded and $||A||\leq ||\sigma||_{p}^{1-p/2}$. One also verifies that $A^* : L^q(\mathcal{D})\rightarrow L^2(\mathcal{D})$ is given by $A^*f(x) = f(x)\sigma(x)^{1-p/2}$, with $||A|| = ||A^*||$.
Define the measurable function $k_0(x,y) := k(x,y)\sigma(x)^{p/2-1}\sigma(y)^{p/2-1}$, and $\sigma_0(x) = k_0(x,x)^{1/2}$. It is positive definite and measurable positive on the diagonal. It verifies
\begin{align}
     ||\sigma_0||_2^2 = \int_{\mathcal{D}}\sigma_0(x)^2dx = \int_{\mathcal{D}}k_0(x,x)dx = \int_{\mathcal{D}}\sigma(x)^pdx= ||\sigma||_p^p < + \infty 
\end{align}
Therefore $\mathcal{E}_{k_0}$, the integral operator over $L^2(\mathcal{D})$ associated to $k_0$, is trace class (Proposition \ref{prop:sobolev_reg_sto_LP}$(ii)$). Observing that $k = (A\otimes A)k_0$ also yields that $\mathcal{E}_k = A \mathcal{E}_{k_0} A^*$.
Thus, from the nuclear norm estimate of \cite{treves2006topological}, Proposition 47.1 pp. 479-480,
\begin{align}
    \nu(\mathcal{E}_k) = \nu(A\mathcal{E}_{k_0}A^*) \leq ||A||\nu(\mathcal{E}_{k_0})||A^*|| \leq ||A||^2 \nu(\mathcal{E}_{k_0}) \leq ||\sigma||_p^{2-p}||\sigma||_p^p = ||\sigma||_p^2
\end{align}
Therefore,
\begin{align}
    \nu(\mathcal{E}_k) \leq \inf_{\substack{A, S \ s.t. \\ \mathcal{E}_k = ASA^*}}||A||^2\nu(S) \leq ||\sigma||_p^2
\end{align}
To finish, note that $C_2 = 1$: therefore, when $p=2$ in equation \eqref{eq:control_p_geq_2}, we recover the fact that $||\sigma||_2^2 = \int_{\mathcal{D}} k(x,x)dx = \nu(\mathcal{E}_k) = \text{Tr}(\mathcal{E}_k)$.
\end{proof}

\section{Sobolev regularity for Gaussian processes : the Hilbert space case, $p=2$}\label{sec:gp_sobolev_L2}
In the case $p=2$, we provide an alternative proof of the integral and spectral criteria of Proposition \ref{prop:sobolev_reg_sto_LP}, based on the study of the "ellipsoids" of Hilbert spaces (see Section \ref{subsub:gaussian_bounded_hilbert}). These geometrical objects are well understood in relation with Gaussian processes (see \cite{dudley1967sizes} or \cite{talagrand2014upper}, Section 2.5). Compared with the general case $p\in (1,+\infty)$, we draw additional links between the different Mercer decompositions of the kernels $\partial^{\alpha,\alpha}k$, the trace of $\mathcal{E}_k^{\alpha}$ and the Hilbert-Schmidt nature of the imbedding of the reproducing kernel Hilbert space (see Section \ref{sub:rkhs} below) associated to $k$ in $H^m(\mathcal{D})$. 

\subsection{Reproducing Kernel Hilbert Spaces (RKHS, \cite{agnan2004})}\label{sub:rkhs}
Consider a general set $\mathcal{D}$ and a positive definite function $k: \mathcal{D}\times\mathcal{D}\rightarrow \mathbb{R}$, i.e. such that given any $n\in\mathbb{N}$ and $(x_1,...,x_n)\in\mathcal{D}^n$, the matrix $(k(x_i,x_j))_{1\leq i,j \leq n}$ is non negative definite. One can then build a Hilbert space $H_k$ of functions defined over $\mathcal{D}$ which contains the functions $k(x,\cdot), x\in\mathcal{D}$ and verifies the \textit{reproducing} identities
\begin{align}
    \langle k(x,\cdot),k(x',\cdot)\rangle_{H_k} &= k(x,x') &&\forall x,x' \in \mathcal{D} \\
    \langle k(x,\cdot),f\rangle_{H_k} &= f(x) &&\forall x \in \mathcal{D}, \ \forall f \in H_k \label{eq:rkhs_f}
\end{align}
$H_k$ is the RKHS of $k$. This space is exactly the set of functions of the form $f(x) = \sum_{i=1}^{+\infty}a_ik(x_i,x)$ such that $||f||_{H_k}^2 = \sum_{i,j=1}^{+\infty}a_i a_jk(x_i,x_j)<+\infty$. If for all $x\in\mathcal{D}$, $k(x,\cdot)$ is measurable, then $H_k$ only contains measurable functions. One may then consider imbedding $H_k$ in some Sobolev space $H^m(\mathcal{D})$. Recall that in $H^m(\mathcal{D})$, functions are equal up to a set of Lebesgue measure zero. If such an imbedding $i:H_k\rightarrow H^m(\mathcal{D})$ is well-defined (i.e. if $f\in H_k$ then its weak derivatives $\partial^{\alpha}f$ exist and lie in $L^2(\mathcal{D})$ for all $|\alpha|\leq m$), we will sometimes use the same notation for $f\in H_k$ and its equivalence class $f\in H^m(\mathcal{D})$; strictly speaking, the latter should be denoted $i(f)$. It may then happen that $i$ is not injective, as with the RKHS associated to the Kronecker delta $k(x,x') = \delta_{x,x'}$ (in this case, we even have $i(H_k) = \{0\})$.
\begin{remark}
In Proposition \ref{prop:sobolev_reg_sto}, we will be interested in the Hilbert-Schmidt nature of the imbedding $i$. However, it may happen that $H_k$ is not separable, such as with the RKHS associated to the Kronecker delta $\delta_{x,x'}$. This results in additional care required for defining the notion of Hilbert Schmidt operators, as the definition from Section \ref{sub:op_th} cannot hold. Still, this case is dealt with in Proposition \ref{prop:sobolev_reg_sto}$(iv)$. See \cite{owhadi2017separability} and \cite{bogachev1998gaussian}, Remark 3.2.9 p. 103 for discussions on non separable RKHS.
\end{remark}

\subsection{Ellipsoids of Hilbert spaces and canonical Gaussian processes \cite{dudley1967sizes}}\label{subsub:gaussian_bounded_hilbert}
Let $(H;\langle,\rangle_H)$ be a separable Hilbert space. We introduce $(V_x)_{x\in H}$ the canonical Gaussian process of $H$, defined as the centered Gaussian process whose covariance function is the inner product of $H$ :
\begin{align}\label{eq:def_canonical}
\mathbb{E}[V_xV_y] = \langle x,y\rangle_H
\end{align}
A subset $K$ of $H$ is said to be Gaussian bounded (GB) if
\begin{align}\label{eq:def_GB}
\mathbb{P}(\sup_{x \in K}|V_x| < + \infty) = 1
\end{align}
The GB property was first introduced for studying the compact sets of Hilbert spaces, see \cite{dudley1967sizes} on that topic.
In equation \eqref{eq:def_GB}, the random variable is defined as $\sup_{x \in K}|V_x| := \sup_{x \in A}|V_x|$ where $A$ is any countable subset of $K$, dense in $K$. Different choices of $A$ only modify $\sup_{x \in K}|V_x|$ on  a set of probability $0$ (\cite{dudley1967sizes}, p. 291), which leaves equation \eqref{eq:def_canonical} unchanged. 
We will use the two following results below, taken from \cite{dudley1967sizes}.
\begin{proposition}[\cite{dudley1967sizes}, p. 293 and \cite{dudley1967sizes}, Proposition 3.4]\label{prop:gb_closure}
We have the two following facts.
\begin{enumerate}[label=(\roman*),wide,labelwidth=0em,labelindent=0pt]
\item If $K$ is a GB-set, then its closed, convex, symmetric hull is a GB-set.
\item The closure of a GB-set is compact.
\end{enumerate}

\end{proposition} 
Given a self-adjoint compact operator $T : H \rightarrow H$, introduce a basis of eigenvectors $x_n$ and its real, positive eigenvalues $\lambda_n$, $\lambda_n \rightarrow 0$. The image of the closed unit ball of $H$, $B = B_H(0,1)$ is the following "ellipsoid" (\cite{dudley1967sizes}, p. 312)
\begin{align}
T(B) = \bigg\{\sum_{\lambda_n > 0}a_nx_n \ s.t. \  \sum_{\lambda_n > 0}{a_n^2}/{\lambda_n^2} \leq 1\bigg\}
\end{align}
The main result we will use is the following.
\begin{proposition}[\cite{dudley1967sizes}, Proposition 6.3]\label{prop:schmidt_ellips} Suppose that $T$ is compact and self-adjoint. Then $T(B)$ is a GB-set if and only if $\sum_{n\in\mathbb{N}} \lambda_n^2 < \infty$, i.e. $T(B)$ is a "Schmidt ellipsoid".
\end{proposition}
We can now state our result pertaining to the $H^m(\mathcal{D})$-regularity of Gaussian processes, given an arbitrary open set $\mathcal{D}\subset\mathbb{R}^d$.
\begin{proposition}[Sample path Hilbert-Sobolev regularity for Gaussian processes]\label{prop:sobolev_reg_sto}
Let $\mathcal{D} \subset \mathbb{R}^d$ be an open set. Let $(U(x))_{x\in \mathcal{D}} \sim GP(0,k)$ be a measurable centered Gaussian process, defined on a probability set $(\Omega, \mathcal{F},\mathbb{P})$, such that its standard deviation function $\sigma$ lies in $L_{loc}^1(\mathcal{D})$. The following statements are equivalent:
\begin{enumerate}[label=(\roman*),wide,labelwidth=0em,labelindent=0pt]
\item (Sample path regularity) The sample paths of $U$ lie in $H^m(\mathcal{D})$ almost surely.
\item (Spectral structure) For all $|\alpha|\leq m$, the distributional derivative $\partial^{\alpha,\alpha}k$  lies in $L^2(\mathcal{D} \times \mathcal{D})$ and the associated integral operator
\begin{align}
\mathcal{E}_k^{\alpha}f(x) = \int_{\mathcal{D}} \partial^{\alpha,\alpha}k(x,y)f(y)dy
\end{align}
is trace class. Equivalently, there exists a representer $k_{\alpha}$ of $\partial^{\alpha,\alpha}k$ in $L^2(\mathcal{D}\times\mathcal{D})$ which is the covariance function of a measurable Gaussian process. Note $\sigma_{\alpha}(x) := k_{\alpha}(x,x)^{1/2}$, then additionally
\begin{align}\label{eq:integ_diag_hilbert}
    \int_{\mathcal{D}}k_{\alpha}(x,x)dx < +\infty
\end{align}
\item (Mercer decomposition) The kernel $k$ has the following Mercer decomposition
\begin{align}
    k(x,y) = \sum_{n=0}^{+\infty}\lambda_n\phi_n(x)\phi_n(y) \ \ \ \text{ in } L^2(\mathcal{D}\times\mathcal{D})
\end{align}
where $(\lambda_n)$ is a non negative sequence and $(\phi_n)$ is an orthonormal basis of $L^2(\mathcal{D})$. Moreover, for all $n\in\mathbb{N}$ such that $\lambda_n \neq 0$, $\partial^{\alpha}\phi_n \in L^2(\mathcal{D})$, $\partial^{\alpha,\alpha}k\in L^2(\mathcal{D}\times\mathcal{D})$, the following equalities hold
\begin{numcases}{}
\label{eq:trace_mercer_hilbert}
    \text{Tr}(\mathcal{E}_k^{\alpha}) = \sum_{n=0}^{+\infty} \lambda_n ||\partial^{\alpha}\phi_n||_2^2 < + \infty\\
\partial^{\alpha,\alpha}k(x,y) = \sum_{n=0}^{+\infty} \lambda_n \partial^{\alpha}\phi_n(x)\partial^{\alpha}\phi_n(y) \ \ \ \text{ in }  L^2(\mathcal{D}\times\mathcal{D})
\end{numcases}

\item \label{pt:imbedding_rkhs} (imbedding of the RKHS) $H_k \subset H^m(\mathcal{D})$, the corresponding natural imbedding $i : H_k \rightarrow H^m(\mathcal{D})$ is continuous and $ii^* : H^m(\mathcal{D}) \rightarrow H^m(\mathcal{D})$ is trace class. Equivalently, $\ker (i)^{\perp}$ endowed with the topology of $H_k$ is a separable Hilbert space and $i : \ker (i)^{\perp} \rightarrow  H^m(\mathcal{D})$ is Hilbert-Schmidt. Moreover,
\begin{align}\label{eq:trace_iistar_enonce}
    \text{Tr}(ii^*) = \sum_{|\alpha|\leq m}\text{Tr}(\mathcal{E}_k^{\alpha})
\end{align}
\end{enumerate}
\end{proposition}
Before proving this result, we discuss Proposition \ref{prop:sobolev_reg_sto} in relation with previous results from the literature.
First, point $\ref{pt:imbedding_rkhs}$ is not without reminding Driscoll's theorem (\cite{kanagawa2018gaussian}, Theorem 4.9) which is widely spread in the machine learning/RKHS community; this theorem states the following. Let $k$ and $r$ be two positive definite functions defined over $\mathcal{D}$, and let $U\sim GP(0,k)$. Suppose that $H_k\subset H_r$ with a Hilbert-Schmidt imbedding, then the sample paths of $U$ lie in $H_r$ almost surely.

Second, Proposition \ref{prop:sobolev_reg_sto} and equation \eqref{eq:integ_diag_hilbert} in particular, is a generalization of Theorem 1 from \cite{SCHEUERER2010} in the case of Gaussian processes; By removing the assumption in \cite{SCHEUERER2010} that the covariance function be continuous on its diagonal as well as its symmetric cross derivatives, the sufficient condition found in \cite{SCHEUERER2010} becomes also necessary.
Finally, Proposition \ref{prop:sobolev_reg_sto} shows that if $p=2$, then in the nuclear decomposition of $\mathcal{E}_k^{\alpha}$ (see Proposition \ref{prop:sobolev_reg_sto_LP}$(iii))$ one can choose $\lambda_n^{\alpha} = \lambda_n$ and $\psi_n^{\alpha} = \partial^{\alpha}\psi_n$. It is not obvious that this should hold when $p\neq2$. 
\begin{example}[Hilbert-Schmidt imbeddings of Sobolev spaces]\label{ex:compare_steinwart}
Proposition \ref{prop:sobolev_reg_sto} can be compared with the results found in \cite{steinwart2019convergence} and its Corollary 4.5 in particular. This corollary states that if $\mathcal{D}\subset\mathbb{R}^d$ is sufficiently smooth, if $H_k \subset H^t(\mathcal{D})$ with a continuous imbedding and if $t > d/2$, then the sample paths of the centered Gaussian process with covariance function $k$ lie in $H^m(\mathcal{D})$ for all real number $m \in [0, t-d/2)$. For example, this is the case if $k$ is a Matérn covariance function of order $t-d/2$; its RKHs is then exactly $H^t(\mathcal{D})$ (\cite{steinwart2019convergence}, Example 4.8).

In the particular case where in addition $m$ is an integer, we recover this result from Proposition \ref{prop:sobolev_reg_sto}. Indeed, it is known that when $m\in (0, t-d/2)$, the imbedding of $H^t(\mathcal{D})$ in $H^m(\mathcal{D})$ is Hilbert-Schmidt. When the involved indexes are non negative integers, this is known as Maurin's theorem (\cite{fournier_adams_sobolev}, Theorem 6.61, p. 202). Maurin's theorem is generalized to fractional indices in \cite{triebel1968approximationszahlen}, Folgerung 1 p. 310 (in German) or \cite{PhamTheLai1973-1974}, Proposition 7.1 (in French). if $H_k \subset H^t(\mathcal{D})$ with a continuous imbedding, then the inclusion map of $H_k$ in $H^m(\mathcal{D})$ is Hilbert-Schmidt for all $m\in[0,t-d/2)\cap\mathbb{N}_0$. From Proposition \ref{prop:sobolev_reg_sto}, we obtain that the sample paths of the corresponding Gaussian process indeed lie in $H^m(\mathcal{D})$.

However, not all RKHS that are subspaces of $H^m(\mathcal{D})$ with a Hilbert-Schmidt imbedding are contained in some $H^t(\mathcal{D})$ with $t > m + d/2$, as the following trivial example shows. Fix any $\varepsilon>0$ and consider the rank one kernel $k(x,x') = f(x)f(x')$ where $f$ is chosen such that $f \in H^m(\mathcal{D})$ and $f\notin H^{m+\varepsilon}(\mathcal{D})$ (choose once and for all a representer of $f$ in $L^2(\mathcal{D})$ so that $f$ is a function in the classical sense). Then $H_k = \text{Span}(f)$ and the imbedding of $H_k$ in $H^m(\mathcal{D})$ is Hilbert-Schmidt since it is rank one; but $H_k \not\subset H^{m+\varepsilon}(\mathcal{D})$. Proposition \ref{prop:sobolev_reg_sto} yields that the associated trivial Gaussian process $U(x)(\omega) = \xi(\omega)f(x)$ where $\xi\sim\mathcal{N}(0,1)$ has its sample paths in $H^m(\mathcal{D})$ (it was obvious in the first place).
\end{example}
\begin{example}[One dimensional case] We build a covariance function which is not pointwise differentiable at any $(q,q')\in\mathbb{Q}\times\mathbb{Q}$, and such that the corresponding Gaussian process has its sample paths in $H^1(\mathbb{R})$.
Let $h_a(x) := \max(0,1-|x-a|)$ be the hat function centered around $a\in\mathbb{R}$. It lies in $H^1(\mathbb{R})$ but it is not differentiable at $x=a,a-1$ and $a+1$. Let $(q_n)$ be an enumeration of $\mathbb{Q}$. Then the following positive definite function over $\mathbb{R}$
\begin{align}
    k(x,x') := \sum_{n=0}^{+\infty} \frac{1}{2^n}h_{q_n}(x)h_{q_n}(x')
\end{align}
is not differentiable in the classical sense at each point $(x,x')$ of the form $(q_n,q_m)$, but the map $ii^*$, with $i : H_k \rightarrow H^1(\mathbb{R})$ the canonical imbedding, is trace-class (use equations \eqref{eq:trace_mercer_hilbert} and \eqref{eq:trace_iistar_enonce}):
\begin{align}
    \text{Tr}(ii^*) &= \text{Tr}(\mathcal{E}_k) + \text{Tr}(\mathcal{E}_k^1) \\
    &\leq \sum_{n=0}^{+\infty} \frac{1}{2^n}||h_{q_n}||_2^2 + \sum_{n=0}^{+\infty} \frac{1}{2^n}||h_{q_n}'||_2^2 \\
    &\leq \sum_{n=0}^{+\infty} \frac{1}{2^n} + \sum_{n=0}^{+\infty} \frac{1}{2^n} \times 2^2 = 10
\end{align}
\end{example}

Before proving Proposition \ref{prop:sobolev_reg_sto}, we shall require a number of lemmas concerning the Mercer decomposition of Hilbert-Schmidt operators over $L^2(\mathcal{D})$. They are proved in Section \ref{sec:proofs}.
\begin{lemma}\label{lemma:extend_balpha}
Let $k$ be a measurable positive definite function defined on an open set $\mathcal{D}$. Suppose that $\sigma \in L_{loc}^1(\mathcal{D})$. Then $k \in L_{loc}^1(\mathcal{D}\times\mathcal{D})$. Given a multi-index $\alpha$, its distributional derivative $D^{\alpha,\alpha}k$ exists and we can introduce the associated continuous bilinear form over $C_c^{\infty}(\mathcal{D)}$
\begin{align}
b_{\alpha}(\varphi,\psi) := D^{\alpha,\alpha}k(\varphi\otimes\psi) =  \int_{\mathcal{D}\times\mathcal{D}}k(x,y)\partial^{\alpha}\varphi(x)\partial^{\alpha}\psi(y)dxdy
\end{align}
Suppose that it verifies the estimate
\begin{align}
\forall \varphi,\psi \in E_2, \ \ |b_{\alpha}(\varphi,\psi)| \leq C_{\alpha}||\varphi||_{2}||\psi||_{2}
\end{align}
where $E_2$ is the set given in Lemma \ref{lemma:E_et_F_countable_line_bilin}. Then $b_{\alpha}$ can be extended to a continuous bilinear form over $L^2(\mathcal{D})$ and there exists a unique bounded, self-adjoint and positive operator $\mathcal{E}_k^{\alpha} : L^2(\mathcal{D}) \longrightarrow L^2(\mathcal{D})$ such that
\begin{align}
\forall \varphi,\psi \in C_c^{\infty}(\mathcal{D}), \ \ b_{\alpha}(\varphi,\psi) = \langle \mathcal{E}_k^{\alpha}\varphi,\psi\rangle_{L^2(\mathcal{D})}
\end{align}
\end{lemma}

\begin{lemma}\label{lemma:HSpos}
Let $k \in L^2(\mathcal{D}\times\mathcal{D})$ be a be a positive definite function and $\alpha$ a multi-index. Suppose that the weak derivative $\partial^{\alpha,\alpha}k$ exists and lies in $L^2(\mathcal{D} \times \mathcal{D})$. Then the associated Hilbert-Schmidt integral operator defined on $L^2(\mathcal{D})$
\begin{align}
(\mathcal{E}_k^{\alpha}f)(x) = \int_{\mathcal{D}} \partial^{\alpha,\alpha}k(x,y)f(y)dy
\end{align}
is self-adjoint and positive.
\end{lemma}

\begin{lemma}\label{lemma:partial_kalpha_trace}
Let $k \in L^2(\mathcal{D}\times\mathcal{D})$ be a positive definite function and $\mathcal{E}_k$ be its associated Hilbert-Schmidt operator. Let
\begin{align}
k(x,y) = \sum_{i=1}^{+\infty}\lambda_i\phi_i(x)\phi_i(y)
\end{align}
be a symmetric, positive expansion of $k$ in $L^2(\mathcal{D}\times\mathcal{D})$ where $(\lambda_i)$ is a positive sequence decreasing to $0$; it may or may not be its Mercer expansion (i.e. $(\phi_i)$ may or may not be an orthonormal basis of $L^2(\mathcal{D})$; they are still assumed to be elements of $L^2(\mathcal{D})$ though). Then
\begin{enumerate}[label=(\roman*),wide,labelwidth=0em,labelindent=0pt]
\item if the partial mixed weak derivative $\partial^{\alpha,\alpha}k$ exists and lies in $L^2(\mathcal{D}\times\mathcal{D})$, then for all $i\in \mathbb{N}$ such that $ \lambda_i \neq 0, \partial^{\alpha}\phi_i \in L^2(\mathcal{D})$.
\item if for all $i\in \mathbb{N}$ such that $ \lambda_i \neq 0, \partial^{\alpha}\phi_i \in L^2(\mathcal{D})$, then
\begin{align}\label{eq:trace_ekal}
\text{Tr}(\mathcal{E}_k^{\alpha}) = \sum_{i=1}^{+\infty}\lambda_i||\partial^{\alpha}\phi_i||_{L^2(\mathcal{D})}^2
\end{align}
whether these quantities are finite or not. If in equation \eqref{eq:trace_ekal}, either one of them is finite, then the series of functions $\sum_{i \in \mathbb{N}}\lambda_i\partial^{\alpha}\phi_i(x)\partial^{\alpha}\phi_i(y)$ is norm convergent in $L^2(\mathcal{D}\times\mathcal{D})$ (i.e. $\sum_{i \in \mathbb{N}}\lambda_i||\partial^{\alpha}\phi_i \otimes\partial^{\alpha}\phi_i||_{L^2} < +\infty$), $\partial^{\alpha,\alpha}k$ lies in $L^2(\mathcal{D}\times\mathcal{D})$ and we have the following equality:
\begin{align}\label{eq:weak_derivative_commute}
\partial^{\alpha,\alpha}k(x,y) = \sum_{i=1}^{+\infty}\lambda_i\partial^{\alpha}\phi_i(x)\partial^{\alpha}\phi_i(y) \ \ \ \text{ in }  \ \ \ L^2(\mathcal{D}\times\mathcal{D})
\end{align}
Equation \eqref{eq:weak_derivative_commute} then holds for asymmetric derivatives, as for all $|\alpha|,|\beta|\leq m$, we also have $\sum_{i \in \mathbb{N}}\lambda_i||\partial^{\beta}\phi_i \otimes\partial^{\alpha}\phi_i||_{L^2} < +\infty$.
\end{enumerate}

\end{lemma}
We can now prove Proposition \ref{prop:sobolev_reg_sto}.
\begin{proof}(Proposition \ref{prop:sobolev_reg_sto}) We successively prove $(ii) \implies (i)$, $(i) \implies (ii)$, $(ii) \iff (iii)$, $(iii) \implies (iv)$ and $(iv) \implies (iii)$.\\ Before all things, the assumptions and Lemma \ref{lemma:u_alpha_phi_GP} show that the sample paths of $U$ lie in $L_{loc}^1(\mathcal{D})$, that the random variable given by the formula
\begin{align}
U_{\varphi}^{\alpha} : \Omega \ni \omega \longmapsto (-1)^{|\alpha|}\int_{\mathcal{D}}U(x)(\omega)\partial^{\alpha}\varphi(x)dx
\end{align}
is well defined and that $(U_{\varphi}^{\alpha})_{\varphi\in F_2}$ is a Gaussian sequence. \\
\underline{$(ii) \implies (i)$ :}
From Lemma \ref{lemma:HSpos}, $\mathcal{E}_k^{\alpha}$ is a self-adjoint, positive Hilbert-Schmidt operator; it is actually trace-class by assumption. We can thus define $A_{\alpha} := \sqrt{\mathcal{E}_k^{\alpha}}$, which is a Hilbert-Schmidt, self-adjoint, positive operator. From Proposition \ref{prop:schmidt_ellips}, $A_{\alpha}(B)$ is a GB-set ($B$ is the closed unit ball of $L^2(\mathcal{D})$). Therefore, using the canonical Gaussian process of $L^2(\mathcal{D})$,
\begin{align}
\mathbb{P}(\sup_{\psi\in A_{\alpha}(B)}|V_{\psi}| < + \infty) = 1
\end{align}
which yields in particular that
\begin{align}\label{eq:V_surF_fini}
\mathbb{P}(\sup_{\varphi\in F_2}|V_{A_{\alpha}(\varphi)}| < + \infty) = 1
\end{align}
We now observe that the two Gaussian sequences $(V_{A_{\alpha}(\varphi)})_{\varphi \in F_2}$ and $(U_{\varphi}^{\alpha})_{\varphi\in F_2}$ have the same finite dimensional marginals. Indeed, they are both centered Gaussian sequences with the same covariance:
\begin{align}
\mathbb{E}[V_{A_{\alpha}(\varphi)}V_{A_{\alpha}(\psi)}] &= \langle A_{\alpha}(\varphi), A_{\alpha}(\psi) \rangle_{L^2} = \langle A_{\alpha}^2(\varphi), \psi \rangle_{L^2} = \langle \mathcal{E}_k^{\alpha}\varphi, \psi \rangle_{L^2} \\
\mathbb{E}[U_{\varphi}^{\alpha}U_{\psi}^{\alpha}] &= \mathbb{E}\bigg[\int_{\mathcal{D}}U(x)\partial^{\alpha}\varphi(x)dx\int_{\mathcal{D}}U(y)\partial^{\alpha}\psi(y)dy\bigg]\nonumber\\
&= \int_{\mathcal{D}\times\mathcal{D}}k(x,y)\partial^{\alpha}\varphi(x)\partial^{\alpha}\psi(y)dxdy \nonumber\\
&= \int_{\mathcal{D}\times\mathcal{D}}\partial^{\alpha,\alpha}k(x,y)\varphi(x)\psi(y)dxdy = \langle \mathcal{E}_k^{\alpha}\varphi, \psi \rangle_{L^2}
\end{align}
As in the proof of Proposition \ref{prop:sobolev_reg_sto_LP}, we deduce that the two random variables $\sup_{\varphi\in F_2}|U_{\varphi}^{\alpha}|$ and $\sup_{\varphi\in F_2}|V_{A_{\alpha}(\varphi)}|$ have the same law, and from equation \eqref{eq:V_surF_fini}, we obtain that
\begin{align}\label{eq:u_al_fini}
\mathbb{P}(\sup_{\varphi\in F_2}|U_{\varphi}^{\alpha}| < + \infty) = \mathbb{P}(\sup_{\varphi\in F_2}|V_{A_{\alpha}(\varphi)}| < + \infty) = 1
\end{align}
Since equation \eqref{eq:u_al_fini} holds for all $|\alpha|\leq m$, this provides a set of probability $1$ on which all the sample paths of $U$ lie in $H^m(\mathcal{D})$, which proves $(i)$. \\
\underline{$(i) \implies (ii)$ :} 
From Lemma \ref{lemma:E_et_F_countable_sobolev} and the assumption from $(i)$,
\begin{align}
\mathbb{P}(\sup_{\varphi\in F_2}|U_{\varphi}^{\alpha}| < + \infty) = 1
\end{align}
From Proposition \ref{lemma:expo_moments}, we have that
\begin{align}
C_{\alpha} := \mathbb{E}\big[\sup_{\varphi \in F_2}|U_{\varphi}^{\alpha}|^2\big] < + \infty
\end{align}
Introduce $b_{\alpha}$, the continuous bilinear form over $C_c^{\infty}(\mathcal{D})$ given by
\begin{align}
b_{\alpha}(\varphi,\psi) = \int_{\mathcal{D}\times\mathcal{D}}&k(x,y)\partial^{\alpha}\varphi(x)\partial^{\alpha}\psi(y)dxdy
\end{align}
Consider now $\varphi$ and $\psi$ in $F_2$. Then,
\begin{align}
|b_{\alpha}(\varphi,\psi)| &= \bigg|\int_{\mathcal{D}\times\mathcal{D}}k(x,y)\partial^{\alpha}\varphi(x)\partial^{\alpha}\psi(y)dxdy\bigg| = 
|\mathbb{E}[U_{\varphi}^{\alpha}U_{\psi}^{\alpha}]|\nonumber \\
&\leq\mathbb{E}[|U_{\varphi_0}^{\alpha}U_{\psi_0}^{\alpha}|] \leq \frac{1}{2}\mathbb{E}\big[(U_{\varphi_0}^{\alpha})^2 + (U_{\psi_0}^{\alpha})^2\big] \leq \mathbb{E}\big[\sup_{\varphi_0\in F}(U_{\varphi_0}^{\alpha})^2\big] = C_{\alpha}
\end{align}
From Lemma \ref{lemma:extend_balpha}, $b_{\alpha}$ can be extended to a continuous bilinear form over $L^2(\mathcal{D})$ and there exists a unique bounded, self-adjoint and positive operator $\mathcal{E}_k^{\alpha}$ which verifies
\begin{align}
\forall \varphi,\psi \in C_c^{\infty}(\mathcal{D}), \ \ \int_{\mathcal{D}\times\mathcal{D}}&k(x,y)\partial^{\alpha}\varphi(x)\partial^{\alpha}\psi(y)dxdy = b_{\alpha}(\varphi,\psi) = \langle \mathcal{E}_k^{\alpha}\varphi,\psi\rangle_{L^2}
\end{align}
Since $\mathcal{E}_k^{\alpha}$ is self-adjoint and positive, we can introduce its square root $A_{\alpha} := \sqrt{\mathcal{E}_k^{\alpha}}$, which is also a bounded, self-adjoint and positive operator. As before, we can introduce $(V_{A_{\alpha}(\varphi)})_{\varphi \in F_2}$ and observe that $(V_{A_{\alpha}(\varphi)})_{\varphi \in F_2}$ and $(U_{\varphi}^{\alpha})_{\varphi\in F_2}$ have the same law. Thus,
\begin{align}
\mathbb{P}(\sup_{\varphi\in F_2}|V_{A_{\alpha}(\varphi)}| < + \infty) = \mathbb{P}(\sup_{\varphi\in F_2}|U_{\varphi}^{\alpha}| < + \infty) = 1
\end{align}
Therefore, $A_{\alpha}(F_2)$ is a GB-set. From Proposition \ref{prop:gb_closure}(ii), $\overline{\text{Conv}(A_{\alpha}(F_2))}$ is compact. One then checks by elementary considerations that $\overline{\text{Conv}(A_{\alpha}(F_2))} = \overline{A_{\alpha}(B)}$, where $B$ is the unit ball of $L^2(\mathcal{D})$. This shows that $A_{\alpha}$ is a compact operator. 
But from Proposition \ref{prop:gb_closure}(i),  $\overline{A_{\alpha}(B)}=\overline{\text{Conv}(A_{\alpha}(F_2))}$ is also a GB-set. From Proposition \ref{prop:schmidt_ellips}, $A_{\alpha}$ is Hilbert-Schmidt and $\mathcal{E}_k^{\alpha}$ is trace-class. In particular, $\mathcal{E}_k^{\alpha}$ is a Hilbert-Schmidt operator with a kernel $k_{\alpha}$ that lies in $L^2(\mathcal{D\times\mathcal{D}})$:
\begin{align}
  \forall \varphi,\psi \in C_c^{\infty}(\mathcal{D}),\ D^{\alpha,\alpha}k(\varphi \otimes\psi) &= \int_{\mathcal{D}\times\mathcal{D}}\ \ k(x,y)\partial^{\alpha}\varphi(x)\partial^{\alpha}\psi(y)dxdy \\
  &= \int_{\mathcal{D}\times\mathcal{D}}k_{\alpha}(x,y)\varphi(x)\psi(y)dxdy = T_{k_{\alpha}}(\varphi \otimes\psi)\label{eq:bilin_schwartz}
\end{align}
Equation \eqref{eq:bilin_schwartz} shows that the distributional derivative $D^{\alpha,\alpha}k$ and the regular distribution $T_{k_{\alpha}}$ coincide on $\mathscr{D}'(\mathcal{D})\otimes\mathscr{D}'(\mathcal{D})$. From the Schwartz kernel theorem (\cite{treves2006topological}, Theorem 51.7), $D^{\alpha,\alpha}k = T_{k_{\alpha}}$ in $\mathscr{D}'(\mathcal{D}\times\mathcal{D})$, which shows that $\partial^{\alpha,\alpha}k$ exists in $L^2(\mathcal{D}\times \mathcal{D})$ and that $\partial^{\alpha,\alpha}k = k_{\alpha}$. For the existence of a representer $k_{\alpha}$ with the desired properties, we refer to the previous Proposition \ref{prop:sobolev_reg_sto_LP}. This finishes to prove $(ii)$. \\
\underline{$(ii) \iff (iii)$:} this equivalence is fully given by Lemma \ref{lemma:partial_kalpha_trace}.\\
\underline{$(iii) \implies (iv)$:} we first study how finite difference operators behave on elements of $H_k$ in order to use Lemma \ref{lemma:linear_form_sobolev}$(iii)$. First, using the reproducing formula \eqref{eq:rkhs_f}, observe that for suitable $x$ and $h\in\mathcal{D}$,
\begin{align}
\Delta_hf(x) = f(x+h)-f(x) = \langle f, k({x+h},\cdot) - k(x,\cdot)\rangle_{H_k} = \langle f, \Delta_h k(x,\cdot) \rangle_{H_k}
\end{align}
More generally, for any finite difference operator $\Delta_h$ of order $l\leq m$, $h = (h_1,...,h_l)$ and any open set $\mathcal{D}_0\Subset\mathcal{D}$ such that $\sum_i|h_i| < \text{dist}(\mathcal{D}_0,\partial\mathcal{D})$,
\begin{align}
\Delta_h f(x) = \langle f, \Delta_h k(x,\cdot) \rangle_{H_k}
\end{align}
The Cauchy-Schwarz inequality in $H_k$ yields
\begin{align}\label{eq:cs_delta}
\Delta_h f(x)^2 \leq ||f||_{H_k}^2 ||\Delta_h k(x,\cdot)||_{H_k}^2
\end{align}
Furthermore, using the bilinearity of $\langle,\cdot,\cdot\rangle_{H_k}$, we have that
\begin{align}\label{eq:bilin_delta}
||\Delta_h k(x,\cdot)||_{H_k}^2 = [(\Delta_h \otimes\Delta_h)k](x,x)
\end{align}
We then deduce that (explanation below)
\begin{align}
\forall f\in H_k,\ ||\Delta_hf||_{L^2(\mathcal{D}_0)}^2 = \int_{\mathcal{D}_0}\big(\Delta_hf\big)(x)^2dx &\leq  ||f||_{H_k}^2\int_{\mathcal{D}_0}[(\Delta_h\otimes\Delta_h)k](x,x)dx \label{eq:almost_trace}\\
&\leq ||f||_{H_k}^2 \sum_{i=1}^{+\infty}\lambda_i \int_{\mathcal{D}_0}(\Delta_h\phi_i)(x)^2dx\label{eq:delta_diff_phi} \\
&\leq ||f||_{H_k}^2 \sum_{i=1}^{+\infty}\lambda_i\Big(||\phi_i||_{H^m}^2|h_1|^2\dotsb |h_l|^2\Big)\label{eq:use_diff_sobolev_lemma} \\
&\leq ||f||_{H_k}^2\bigg(\sum_{|\alpha|\leq m}\text{Tr}(\mathcal{E}_k^{\alpha})\bigg)\big(|h_1|^2\dotsb |h_l|^2\big)\label{eq:rkhs_trace}
\end{align}
We used equations \eqref{eq:cs_delta} and \eqref{eq:bilin_delta} to obtain equation \eqref{eq:almost_trace}. In equation \eqref{eq:delta_diff_phi}, we distributed $\Delta_h\otimes\Delta_h$ on the Mercer decomposition of $k$ (which exists by the assumption $(iii)$). In equation \eqref{eq:use_diff_sobolev_lemma}, we used the fact that $\phi_i\in H^m(\mathcal{D})$ (see Lemma \ref{lemma:partial_kalpha_trace}$(i)$) conjointly with the finite difference control of Lemma \ref{lemma:linear_form_sobolev}$(iii)$. In equation \eqref{eq:rkhs_trace}, the we used the trace equality from Lemma \ref{lemma:partial_kalpha_trace}$(ii)$. From equation \eqref{eq:rkhs_trace} and Lemma \ref{lemma:linear_form_sobolev}$(iii)$ again, we obtain that $f$ lies in $H^m(\mathcal{D})$. Consider now any open set $\mathcal{D}_0\Subset \mathcal{D}$. Equation \eqref{eq:rkhs_trace} applied to $\delta_h^{\alpha}$, the finite difference approximation of $\partial^{\alpha}$ from equation \eqref{eq:def_fd_approx_alpha} with suitably chosen $h = (h_1,...,h_l)\in(\mathbb{R}_+^*)^d$, yields that
\begin{align}\label{eq:fd_control_alpha}
    \forall f\in H_k,\ ||\delta_h^{\alpha}f||_{L^2(\mathcal{D}_0)}^2 \leq ||f||_{H_k}^2\bigg(\sum_{|\alpha|\leq m}\text{Tr}(\mathcal{E}_k^{\alpha})\bigg)
\end{align}
From equation \eqref{eq:fd_cv_LP},
 $||\delta_h^{\alpha}f - \partial^{\alpha}f||_{L^2(\mathcal{D}_0)}$ converges to zero when $h\rightarrow 0$. 
Thus, the left hand side of equation \eqref{eq:fd_control_alpha} converges to $||\partial^{\alpha}f||_{L^2(\mathcal{D}_0)}^2$ when $h\rightarrow 0$. Writing $||\partial^{\alpha}f||_{L^2(\mathcal{D})} =\sup_{\mathcal{D}_0\Subset\mathcal{D}} ||\partial^{\alpha}f||_{L^2(\mathcal{D}_0)}$, equation \eqref{eq:fd_control_alpha} yields that in fact
\begin{align}\label{eq:partial_alpha_control}
    \forall f\in H_k,\ ||\partial^{\alpha} f||_{L^2(\mathcal{D})}^2 \leq ||f||_{H_k}^2\bigg(\sum_{|\alpha|\leq m}\text{Tr}(\mathcal{E}_k^{\alpha})\bigg)
\end{align}
Summing the inequality \eqref{eq:partial_alpha_control} for all $|\alpha|\leq m$, we obtain that
\begin{align}
||f||_{H^m} \leq C||f||_{H_k}
\end{align}
with $C = \big(N\sum_{|\alpha|\leq m}\text{Tr}(\mathcal{E}_k^{\alpha})\big)^{1/2}$ and $N$ is the number of indexes $\alpha$ such that $|\alpha|\leq m$. Therefore $H_k\subset H^m(\mathcal{D})$ and the corresponding imbedding $i : H_k\rightarrow H^m(\mathcal{D})$ is continuous. Using the reproducing formula \eqref{eq:rkhs_f}, its transpose $i^* : H^m(\mathcal{D}) \rightarrow H_k$ is given by
\begin{align}\label{eq:formular_transpose}
i^*(f)(x) &= \langle i^*(f),k_x\rangle_{H_k} = \langle f,i(k_x)\rangle_{H^m} = \sum_{|\alpha| \leq m} \int_{\mathcal{D}} \partial_y^{\alpha}k(x,y)\partial^{\alpha}f(y)dy
\end{align}
Above, $\partial^{\alpha}_y$ denotes differentation w.r.t. the $y$ coordinate (note that $i^*(f)$ is indeed defined pointwise, since $i^*(f) \in H_k$).
Let $(\psi_j)$ be an orthonormal basis of $H^m(\mathcal{D})$ and $k = \sum_i \lambda_i\psi_i\otimes\psi_i$ be the Mercer decomposition of $k$ provided by the assumption $(iii)$. The trace of the positive self-adjoint operator $ii^*$ is given by (explanation below)
\begin{align}
\text{Tr}(ii^*) &= \sum_{j}\langle \psi_j,ii^*(\psi_j)\rangle_{H^m} = \sum_{j}\sum_{|\beta|\leq m} \langle \partial^{\beta}\psi_j,\partial
^{\beta}ii^*(\psi_j)\rangle_{L^2}\nonumber \\
&= \sum_{j}\sum_{|\beta|\leq m} \int_{\mathcal{D}}\partial^{\beta}\psi_j(x)\partial^{\beta}ii^*(\psi_j)(x)dx\nonumber\\
&= \sum_j\sum_{|\beta| \leq m} \int_{\mathcal{D}}\partial^{\beta}\psi_j(x)\partial_x^{\beta}\sum_{|\alpha| \leq m}\int_{\mathcal{D}} \partial_y^{\alpha}k(x,y)\partial^{\alpha}\psi_j(y)dydx \label{eq:imbed_transpose}\\
&= \sum_j \sum_i \lambda_i \sum_{|\alpha| \leq m}\sum_{|\beta| \leq m} \int_{\mathcal{D}\times\mathcal{D}} \partial^{\beta}\phi_i(x)\partial^{\alpha}\phi_i(y)\partial^{\beta}\psi_j(x)\partial^{\alpha}\psi_j(y)dydx \label{eq:distribute_norm_cv}\\
&= \sum_j \sum_i \lambda_i \bigg(\sum_{|\alpha| \leq m}\int_{\mathcal{D}} \partial^{\alpha}\phi_i(x)\partial^{\alpha}\psi_j(x)dx\bigg)^2 = \sum_j \sum_i \lambda_i \bigg(\sum_{|\alpha| \leq m}\langle \partial^{\alpha}\phi_i,\partial^{\alpha}\psi_j \rangle_{L^2}\bigg)^2\nonumber\\
&= \sum_i \lambda_i \sum_j \langle \phi_i,\psi_j\rangle_{H^m}^2 =  \sum_i \lambda_i ||\phi_i||_{H^m}^2 = \sum_{|\alpha| \leq m}\sum_i \lambda_i ||\partial^{\alpha}\phi_i||_{L^2}^2 = \sum_{|\alpha| \leq m} \text{Tr}(\mathcal{E}_k^{\alpha}) \label{eq:trace_iietoile}
\end{align}
In equation \eqref{eq:imbed_transpose}, we used the fact that $i^*(\psi_j)$ given by equation \eqref{eq:formular_transpose} is a representer of $ii^*(\psi_j)$ in $H^m(\mathcal{D})$. In equation \eqref{eq:distribute_norm_cv}, we used the fact that the series of functions $\sum_i \lambda_i \partial^{\beta}\phi_i\otimes\partial^{\alpha}\phi_i$ is norm convergent (Lemma \ref{lemma:partial_kalpha_trace}$(ii)$) to distribute the partial derivatives over to the Mercer decomposition of $k$. We also used Fubini's and Tonelli's theorems ad libitum, as all the series  $\sum_i \lambda_i \partial^{\beta}\phi_i\otimes\partial^{\alpha}\phi_i$ are norm convergent.
Since $\sum_{|\alpha| \leq m} \text{Tr}(\mathcal{E}_k^{\alpha})$ is finite by assumption, equation \eqref{eq:trace_iietoile} finishes to prove $(iv)$ when $H_k$ is separable. 

When $H_k$ is not separable, observe that $\text{ker}(i)$ is closed in $H_k$ since $i$ is continuous. Therefore $H_k = \text{ker}(i) \oplus\text{ker}(i)^{\perp}$ and $\text{ker}(i)^{\perp}$ endowed with the topology of $H_k$ is a Hilbert space. Moreover, the following restriction/corestriction of $i$
\begin{align}
j := i_{|\text{ker}(i)^{\perp}}^{|\text{im}(i)}:\text{ker}(i)^{\perp}\rightarrow \text{im}(i)    
\end{align}
is continuous, linear and one to one. From the Banach bounded inverse theorem, $j^{-1}: \text{im}(i) \rightarrow\text{ker}(i)^{\perp}$ is continuous and $\text{ker}(i)^{\perp}$ is homeomorphic to $\text{im}(i)$, which is separable as a subspace of the separable Hilbert space $H^{m}(\mathcal{D})$. Thus $\text{ker}(i)^{\perp}$ is separable. Finally observe that $ii^* = jj^*$, so that equation \eqref{eq:trace_iietoile} indeed yields that $j$ is Hilbert-Schmidt.\\
\underline{$(iv) \implies (iii)$}: by assumption, $ii^*$ is a compact self-adjoint positive operator acting on the Hilbert space $H^m(\mathcal{D})$. There exists a decreasing positive sequence $(\mu_j)_{j\in\mathbb{N}}$ and a orthonormal basis of eigenvectors of $ii^*, (\psi_j)_{j\in\mathbb{N}}$ such that for all $f\in H^m(\mathcal{D})$,
\begin{align}
ii^*(f) = \sum_{j=1}^{+\infty}\mu_j \langle \psi_j,f\rangle_{H^m}\psi_j \ \ \ \text{in} \ H^m(\mathcal{D})
\end{align}
Since $ii^*$ is assumed trace class,
\begin{align}\label{eq:trace_Hm}
\sum_{|\alpha|\leq m}\sum_{j=1}^{+\infty}\mu_j||\partial^{\alpha}\psi_j||_{L^2}^2 =  \sum_{j=1}^{+\infty} \mu_j ||\psi_j||_{H^m}^2 = \sum_{j=1}^{+\infty} \mu_j < + \infty
\end{align}
We now show that the following equality holds in $L^2(\mathcal{D}\times\mathcal{D})$:
\begin{align}\label{eq:decomp_k_L2}
k(x,y) = \sum_{j=1}^{+\infty}\mu_j\psi_j(x) \psi_j(y)
\end{align}
In conjunction with equation \eqref{eq:trace_Hm}, this equation will allow us to use Lemma \ref{lemma:partial_kalpha_trace}$(ii)$, which will imply the point $(ii)$.

First, one easily shows that $\sum_{j=1}^{+\infty}\mu_j\psi_j\otimes\psi_j$, the right-hand side of equation \eqref{eq:decomp_k_L2}, is indeed in $L^2(\mathcal{D}\times\mathcal{D})$ (e.g. use that $\sum_j \mu_j <+\infty$). Equation \eqref{eq:int_0} will then show that $k$ is indeed in $L^2(\mathcal{D}\times\mathcal{D})$.
Now, decompose $i(k_x)\in H^m(\mathcal{D})$ on the basis $(\psi_j)_{j\in\mathbb{N}}$ given any $x\in \mathcal{D}$: 
\begin{align}\label{eq:decomp_ikx}
i(k_x) = \sum_{j=1}^{+\infty} \langle \psi_j,i(k_x)\rangle_{H^m} \psi_j \ \ \ \text{in} \ H^m(\mathcal{D})
\end{align}
In equation \eqref{eq:decomp_ikx}, the scalar $\langle\psi_j,i(k_x)\rangle_{H^m}$ is obtained through the reproducing formula \eqref{eq:rkhs_f}:
\begin{align}\label{eq:scalar_ikx}
\langle\psi_j,i(k_x)\rangle_{H^m} = \langle i^*(\psi_j),k_x\rangle_{H_k} = 
i^*(\psi_j)(x)
\end{align}
Moreover, $\psi_j$ is an eigenvector of $ii^*$: $\mu_j\psi_j = ii^*(\psi_j)$ in $H^m(\mathcal{D})$. In particular,
\begin{align}
||\mu_j\psi_j - ii^*(\psi_j)||_{L^2(\mathcal{D})} = 0 \label{eq:eigen_ii}
\end{align}
But the pointwise defined function $i^*(\psi_j)$ is a representer of $ii^*(\psi_j)$ in $H^m(\mathcal{D})$, since $i$ is the imbedding of $H_k$ in $H^m(\mathcal{D})$. 
Setting $S = \sum_j\mu_j = \text{Tr}(ii^*)$, one has (explanation below)
\begin{align}
\bigg|\bigg|k-\sum_{j=1}^{+\infty}\mu_j\psi_j\otimes \psi_j\bigg|\bigg|_{L^2(\mathcal{D}\times\mathcal{D})}^2 &= \int_{\mathcal{D}\times\mathcal{D}}\Big(k(x,y) - \sum_{j=1}^{+\infty}\mu_j\psi_j(x) \psi_j(y)\Big)^2dxdy\\
&= \int_{\mathcal{D}}\int_{\mathcal{D}} \Big(k_x(y) - \sum_{j=1}^{+\infty}\mu_j\psi_j(x) \psi_j(y)\Big)^2dydx \label{eq:tonelli}\\
&= \int_{\mathcal{D}}\int_{\mathcal{D}} \Big(i(k_x)(y) - \sum_{j=1}^{+\infty}\mu_j\psi_j(x) \psi_j(y)\Big)^2dydx\label{eq:imbed_kx}\\
&= \int_{\mathcal{D}}\int_{\mathcal{D}} \bigg( \sum_{j=1}^{+\infty}\mu_j\psi_j(y)\big({\mu_j}^{-1}i^*(\psi_j)-\psi_j(x)\big)\bigg)^2dydx\label{eq:facto} \allowdisplaybreaks\\
&\leq \int_{\mathcal{D}}\int_{\mathcal{D}}S \sum_{j=1}^{+\infty}\mu_j\psi_j(y)^2\big({\mu_j}^{-1}i^*(\psi_j)-\psi_j(x)\big)^2dydx\label{eq:jensen_square}\\
&\leq S\sum_{j=1}^{+\infty}\mu_j\int_{\mathcal{D}}\psi_j(y)^2dy \int_{\mathcal{D}}\big({\mu_j}^{-1}ii^*(\psi_j)-\psi_j(x)\big)^2dx\label{eq:reorder}\\
&\leq S\sum_{j=1}^{+\infty}\mu_j ||\psi_j||_{L^2(\mathcal{D})}^2||{\mu_j}^{-1}ii^*(\psi_j)-\psi_j||_{L^2(\mathcal{D})}^2 = 0 \label{eq:int_0}
\end{align}
Above, we used Tonelli's theorem in equation \eqref{eq:tonelli}. We imbedded $k_x$ in $H^m(\mathcal{D})$ in equation \eqref{eq:imbed_kx}. We used equations \eqref{eq:decomp_ikx} and \eqref{eq:scalar_ikx} in equation \eqref{eq:facto}. We used Jensen's discrete inequality on the squaring function $(\cdot)^2$  with the weights $\mu_j/S$ ($\mu_j/S\geq 0, \sum_j \mu_j/S = 1$) in equation \eqref{eq:jensen_square}. We imbedded $i^*(\psi_j)$ in $H^m(\mathcal{D})$ and used Tonelli's theorem in equation \eqref{eq:reorder}. We used equation \eqref{eq:eigen_ii} in equation \eqref{eq:int_0}.

Therefore we have proved that equation \eqref{eq:decomp_k_L2} holds. By the assumption that $ii^*$ is trace class and using Lemma \ref{lemma:partial_kalpha_trace}$(ii)$,
\begin{align}
\sum_{|\alpha| \leq m} \text{Tr}(\mathcal{E}_k^{\alpha}) = \sum_{|\alpha|\leq m}\sum_{j=1}^{+\infty}\mu_j||\partial^{\alpha}\psi_j||_{L^2}^2 =  \sum_{j=1}^{+\infty} \mu_j ||\psi_j||_{H^m}^2 = \sum_{i} \mu_j = \text{Tr}(ii^*) < + \infty
\end{align}
Therefore, Lemma \ref{lemma:partial_kalpha_trace}$(ii)$ implies that every $\mathcal{E}_k^{\alpha}$ is indeed trace-class, which shows $(ii)$. 
\end{proof}

\section{Concluding remarks and perspectives}\label{sec:ccl} 
Given $p\in(1,+\infty)$ and $m\in\mathbb{N}_0$, we showed that the $W^{m,p}$-Sobolev regularity of integer order of a measurable Gaussian process $((U(x))_{x\in\mathcal{D}}\sim GP(0,k)$ is fully equivalent to the fact that  ${\partial^{\alpha,\alpha}k}$ lies in ${L^p(\mathcal{D}\times\mathcal{D})}$ combined with the integrability in $L^p(\mathcal{D})$ of the associated standard deviation. Using  general results on Gaussian measures over Banach spaces of type $2$ and cotype  $2$, we translated this criteria as the existence of suitable nuclear decompositions of the covariance. These can be understood as generalizations to Banach spaces of the eigenfunction expansion of symmetric, positive and trace class operators. In the Hilbert space case $p=2$, we linked this property with the Hilbert-Schmidt nature of the imbedding of the RKHS in $H^m(\mathcal{D})$, and gave explicit formulas for the traces of the involved integral operators in terms of the Mercer decomposition of the kernel.

The results presented in this article provide a theoretical background w.r.t. the use of Gaussian processes for solving physics-related machine learning problems, in particular when modelling solutions of PDEs as sample paths of some Gaussian process. These results also come along with certain key quantities for controlling the Sobolev norm of the corresponding sample paths (see Remark \ref{rk:control_sobolev_norm}). The application of the Gaussian process principles identified here to PDE-related machine learning problems is certainly an interesting continuation of the results of this article.

The following directions are interesting for generalizing the results presented here. First, similar spectral/integral criteria should be obtained for fractional Sobolev and Besov spaces. Second, similar results should be sought to tackle the cases $p=1$ and $p=+\infty$. Linked to the case $p=1$, results should be sought for the space of functions of bounded variations (\cite{brezis2011functional}, p. 269), which are important in many problems related to physics. The following open questions are also relevant for Gaussian process theory: $(i)$ Can the small ball problem for Gaussian
processes whose sample paths lie in a Sobolev spaces be tackled only using spectral properties of the covariance operator, e.g. its nuclear norm?
$(ii)$ Are all Gaussian measures over $W^{m,p}(\mathcal{D})$ induced by some Gaussian process? Proposition \ref{prop_bogachev_LP} shows that it is the case for $m=0$, i.e. $L^p(\mathcal{D})$.

\section{Appendix : proof of intermediary results and lemmas}\label{sec:proofs}
\begin{proof} (Lemma \ref{lemma:linear_form_sobolev}) This proof follows exactly the lines of the proof of Proposition 9.3 from \cite{brezis2011functional}.\\
\underline{$(i) \iff (ii)$:} suppose that $u \in W^{m,p}(\mathcal{D})$, use the fact that the distributional derivative $D^{\alpha}u$ is a regular distribution represented by a function that lies in $L^p(\mathcal{D})$, denoted by $\partial^{\alpha}u$ :
\begin{align}
\forall \varphi \in C_c^{\infty}(\mathcal{D}),\ \ \ \int_{\mathcal{D}}u(x)\partial^{\alpha}\varphi(x)dx = (-1)^{|\alpha|}\int_{\mathcal{D}}\partial^{\alpha}u(x)\varphi(x)dx
\end{align}
Hölder's inequality yields \eqref{eq:LP_control} with $C_{\alpha} = ||\partial^{\alpha}u||_{L^p}$.
Conversely, suppose that \eqref{eq:LP_control} holds and consider any $|\alpha|\leq m$. Since $C_c^{\infty}(\mathcal{D})$ is dense in $L^q(\mathcal{D})$ (whatever the open set $\mathcal{D}$, \cite{fournier_adams_sobolev}, section 2.30), equation \eqref{eq:LP_control} shows that the linear form $L_{\alpha}:\varphi \longmapsto (-1)^{|\alpha|}\int_{\mathcal{D}}u(x)\partial^{\alpha}\varphi(x)dx, \varphi\in C_c^{\infty}(\mathcal{D}),$ can be extended to a continuous linear form over $L^q(\mathcal{D})$. From Riesz' representation lemma, there exists $v_{\alpha} \in L^p(\mathcal{D})$ such that $L_{\alpha}(\varphi) = \langle v_{\alpha},\varphi\rangle_{L^p,L^q}$ for all $\varphi \in L^q(\mathcal{D})$. In particular, this is valid for all $\varphi \in C_c^{\infty}(\mathcal{D})$, which shows that for all $|\alpha|\leq m, \partial^{\alpha}u$ exists and is equal to $v_{\alpha}$. Thus $u\in W^{m,p}(\mathcal{D})$. Finally, Hölder's inequality and the density of $C_c^{\infty}(\mathcal{D})$ in $L^q(\mathcal{D})$ yield
\begin{align*}
||\partial^{\alpha}u||_{L^p(\mathcal{D})} = \sup_{\varphi \in C_c^{\infty}(\mathcal{D})\setminus\{0\}} \Bigg|\int_{\mathcal{D}}u(x)\frac{\partial^{\alpha}\varphi(x)}{||\varphi||_{L^q(\mathcal{D})}}dx \Bigg|
\end{align*}
\underline{$(iii) \implies (ii)$:} suppose $(iii)$, let us show $(ii)$. Let $|\alpha|\leq  m$ and let $\varphi \in C_c^{\infty}(\mathcal{D})$. Note $K:=\text{Supp}(\varphi)$ its compact support and consider an open set $\mathcal{D}_0$ such that $K \subset\mathcal{D}_0 \Subset \mathcal{D}$. Let $h = (h_1,...,h_d)$ be such that $\sum_i |h_i| < \text{dist}(\mathcal{D}_0,\partial \mathcal{D})$. Recall that $\delta_{h}^{\alpha}$ from equation \eqref{eq:def_fd_approx_alpha} is a finite difference approximation of $\partial^{\alpha}$ and from $(iii)$,
\begin{align}\label{eq:delta_alpha_control}
\bigg|\int_{\mathcal{D}}\delta_{h}^{\alpha}u(x) \varphi(x)dx\bigg| \leq ||\varphi||_{L^q(\mathcal{D}_0)}||\delta_{h}^{\alpha}u||_{L^p(\mathcal{D}_0)} \leq C ||\varphi||_{L^q(\mathcal{D})}
\end{align}

Note also that we have the discrete integration by parts formula since $h$ is suitably chosen:
\begin{align}\label{eq:discrete_ibp}
\int_{\mathcal{D}}\delta_{h}^{\alpha}u(x) \varphi(x)dx = \int_{\mathcal{D}}u(x) (\delta_{h}^{\alpha})^*\varphi(x)dx
\end{align}

Therefore,
\begin{align}\label{eq:discrete_ibp_control}
\bigg|\int_{\mathcal{D}}u(x)(\delta_{h}^{\alpha})^*\varphi(x) dx\bigg| \leq C ||\varphi||_{L^q(\mathcal{D})}
\end{align}
The Lebesgue dominated convergence theorem yields that the left hand side converges to $\big|\int_{\mathcal{D}}u(x)\partial^{\alpha}\varphi(x)dx\big|$. We therefore have $(ii)$. \\
\underline{$(i) \implies (iii)$:}
We will use recursively the fact that if $f\in W^{1,p}(\mathcal{D})$, then for all  $\mathcal{D}_0 \Subset \mathcal{D}$ and $h \in \mathbb{R}^d$ such that $|h|<\text{dist}(\mathcal{D}_0,\partial\mathcal{D})$, there exists an open set $\mathcal{D}_1\Subset\mathcal{D}$ which verifies $\mathcal{D}_0+th \subset\mathcal{D}_1$ for all $t\in[0,1]$ and
\begin{align}\label{eq:d_0_d_1}
||\Delta_hf||_{L^p(\mathcal{D}_0)}^p = ||\tau_h f - f||_{L^p(\mathcal{D}_0)}^p\leq |h|^p ||\nabla f ||_{L^p(\mathcal{D}_1)}^p = |h|^p \sum_{j=1}^d||\partial_{x_j} f ||_{L^p(\mathcal{D}_1)}^p
\end{align}
(this is equation 4 p. 268 in \cite{brezis2011functional}, found in the proof of Proposition 9.3 in \cite{brezis2011functional}). 
First, one easily checks that weak partial derivatives and finite difference operators all commute together. 
Let $l\leq m$, $\mathcal{D}_0 \Subset \mathcal{D}$ and $h = (h_1,...,h_{l}) \in (\mathbb{R}^d)^{l}$ such that $\sum_i |h_i| < \text{dist}(\mathcal{D}_0,\partial \mathcal{D})$. Recall that
\begin{align}
\Delta_h = \prod_{i=1}^{l}\Delta_{h_i}
\end{align}
Note now that $\prod_{i=2}^{l}\Delta_{h_i}u$ lies in $W^{1,p}(\mathcal{D})$. Since $|h_1| \leq \sum_i |h_i| < \text{dist}(\mathcal{D}_0,\partial \mathcal{D})$, from equation \eqref{eq:d_0_d_1} there exists an open set $\mathcal{D}_1\Subset \mathcal{D}$ such that $\mathcal{D}_0 +th_1\subset \mathcal{D}_1$ for all $t\in[0,1]$. Moreover, one can choose $\mathcal{D}_1$ small enough so that $\text{dist}(\mathcal{D}_1,\partial\mathcal{D}) < \sum_{i=2}^l|h_i|$.
\begin{align}
||\Delta_h u||_{L^p(\mathcal{D}_0)}^p &= \Big|\Big| \Delta_{h_1}\prod_{i=2}^{l}\Delta_{h_i}u \Big|\Big|_{L^p(\mathcal{D}_0)}^p \leq |h_1|^p \Big|\Big|\nabla \prod_{i=2}^{l}\Delta_{h_i} u \Big|\Big|_{L^p(\mathcal{D}_1)}^p \label{eq:brezis_control}\\
&\leq |h_1|^p \sum_{j=1}^d \Big|\Big|\partial_{x_j}\prod_{i=2}^{l}\Delta_{h_i}u \Big|\Big|_{L^p(\mathcal{D}_1)}^p \label{eq:apply_brezis_diff} \leq |h_1|^p \sum_{j=1}^d \Big|\Big|\prod_{i=2}^{l}\Delta_{h_i}\big( \partial_{x_j}u\big) \Big|\Big|_{L^p(\mathcal{D}_1)}^p
\end{align}
We used equation \eqref{eq:d_0_d_1} in equation \eqref{eq:brezis_control} which then yields equation \eqref{eq:apply_brezis_diff}. 
But note that for all $j$, $ \partial_{x_j}u \in W^{1,p}(\mathcal{D})$. One can then proceed by induction and perform the above step sequentially over $i \in\{2,...,l\}$, which yields a sequence of open sets $\mathcal{D}_0 \subset \mathcal{D}_1 \subset ... \subset \mathcal{D}_{l} \Subset \mathcal{D}$ such that
\begin{align}
||\Delta_h u||_{L^p(\mathcal{D}_0)}^p &\leq |h_1|^p \times ... \times |h_{l}|^p \sum_{|\beta|\leq l} ||\partial^{\beta}u||_{L^p(\mathcal{D}_l)}^p \\
&\leq |h_1|^p \times ... \times |h_{l}|^p ||u||_{W^{l,p}(\mathcal{D})}^p \leq |h_1|^p \times ... \times |h_{l}|^p ||u||_{W^{m,p}(\mathcal{D})}^p
\end{align}
which shows equation \eqref{eq:big_delta_control_sobolev} with $C = ||u||_{W^{m,p}(\mathcal{D})}$. We finally show equation \eqref{eq:fd_cv_LP}. First,
for all $u\in C_c^{\infty}(\mathcal{D})$ and suitably chosen $x$ and $h_i$, recall that
\begin{align}
    u(x+h_ie_i) - u(x) = \int_0^1 h_ie_i\cdot\nabla u(x+th_ie_i)dt
\end{align}
which yields that
\begin{align}
    \frac{u(x+h_ie_i) - u(x)}{h_i} - \partial_{x_i}u(x) = \int_0^1 \partial_{x_i}u(x+th_ie_i)-\partial_{x_i}u(x)dt
\end{align}
Generalizing this formula to $\partial^{\alpha}$ with $\alpha = (\alpha_1,...,\alpha_d)\in\mathbb{N}_0^d$,
\begin{align}
 \delta_h^{\alpha}u(x) -\partial^{\alpha}u(x) &= \int_{[0,1]^{\alpha_1}}\dotsb\int_{[0,1]^{\alpha_d}} \partial^{\alpha}u\Big(x + \sum_{i=1}^d h_i\Big(\sum_{j=1}^{\alpha_i}t_{ij}\Big)e_i\Big)\prod_{i=1}^{d}\prod_{j=1}^{\alpha_i}dt_{ij} -\partial^{\alpha}u(x) \\
 = \int_{[0,1]^{\alpha_1}}&\dotsb\int_{[0,1]^{\alpha_d}} \bigg(\partial^{\alpha}u\Big(x + \sum_{i=1}^d h_i\Big(\sum_{j=1}^{\alpha_i}t_{ij}\Big)e_i\Big)-\partial^{\alpha}u(x)\bigg)\prod_{i=1}^{d}\prod_{j=1}^{\alpha_i}dt_{ij}\label{eq:integral_rep_partial_alpha}
\end{align}
Integrating equation \eqref{eq:integral_rep_partial_alpha} over $\mathcal{D}_0$ and using Jensen's inequality for the function $(\cdot)^p$ and the Lebesgue measure $\prod_{i=1}^{d}\prod_{j=1}^{\alpha_i}dt_{ij}$ over the unit cube $[0,1]^{\alpha_1}\times \dotsb \times [0,1]^{\alpha_d} = [0,1]^{|\alpha|}$ yields
\begin{align}
 ||\delta_h^{\alpha}u -&\partial^{\alpha}u||_{L^p(\mathcal{D}_0)} = \int_{\mathcal{D}_0}|\delta_h^{\alpha}u(x) -\partial^{\alpha}u(x)|^pdx \\
 &\leq \int_{[0,1]^{\alpha_1}}\dotsb\int_{[0,1]^{\alpha_d}}\int_{\mathcal{D}_0}\bigg| \partial^{\alpha}u\Big(x + \sum_{i=1}^d h_i\Big(\sum_{j=1}^{\alpha_i}t_{ij}\Big)e_i\Big)-\partial^{\alpha}u(x)\bigg|^pdx\prod_{i=1}^{d}\prod_{j=1}^{\alpha_i}dt_{ij}\\
 &\leq \int_{[0,1]^{\alpha_1}}\dotsb\int_{[0,1]^{\alpha_d}}\bigg|\bigg| \tau_{\sum_{i=1}^d h_i\big(\sum_{j=1}^{\alpha_i}t_{ij}\big)e_i}\partial^{\alpha}u-\partial^{\alpha}u\bigg|\bigg|_{L^p(\mathcal{D}_0)}^p\prod_{i=1}^{d}\prod_{j=1}^{\alpha_i}dt_{ij}\label{eq:lebesgue_dom}
\end{align}
It is standard that for all $v\in L^p(\mathcal{D})$ and $\mathcal{D}_0\Subset\mathcal{D}$, $||\tau_hv-v||_{L^p(\mathcal{D}_0)} \rightarrow 0$ when $h\rightarrow 0$. The following control also holds uniformly in $h$
\begin{align}
    \bigg|\bigg| \tau_{\sum_{i=1}^d h_i\big(\sum_{j=1}^{\alpha_i}t_{ij}\big)e_i}\partial^{\alpha}u-\partial^{\alpha}u\bigg|\bigg|_{L^p(\mathcal{D}_0)} &\leq \Big|\Big| \tau_{\sum_{i=1}^d h_i\big(\sum_{j=1}^{\alpha_i}t_{ij}\big)e_i}\partial^{\alpha}u\Big|\Big|_{L^p(\mathcal{D}_0)} + || \partial^{\alpha}u||_{L^p(\mathcal{D}_0)}\nonumber \\
    &\leq 2|| \partial^{\alpha}u||_{L^p(\mathcal{D})}
\end{align}
Thus the Lebesgue dominated convergence theorem applied to equation \eqref{eq:lebesgue_dom} yields that $||\partial^{\alpha}u - \partial^{\alpha}u||_{L^p(\mathcal{D}_0)}$ converges to zero when $h\rightarrow 0$. One then readily generalizes the density result from \cite{brezis2011functional}, Theorem 9.2 to multiple derivatives, to obtain that $||\delta_h^{\alpha}u - \partial^{\alpha}u||_{L^p(\mathcal{D}_0)} \rightarrow 0$ holds for any $f\in W^{m,p}(\mathcal{D})$.
\end{proof}
\begin{proof}(Lemma \ref{lemma:E_et_F_countable_sobolev}) 
We begin by explicitly constructing the family $(\Phi_n^q)$. First, use the fact that $L^q(\mathcal{D})$ is a separable Banach space (\cite{fournier_adams_sobolev}, Theorem 2.21) : let $(f_n)_{n \in\mathbb{N}} \subset L^q(\mathcal{D})$ be a dense countable subset of $L^q(\mathcal{D})$. For all $n \in \mathbb{N}$, let $(\phi_{nm})_{m \in \mathbb{N}} \subset C_c^{\infty}(\mathcal{D})$ be such that $\phi_{nm} \longrightarrow f_n$ for the $L^q(\mathcal{D})$ topology (recall that $C_c^{\infty}(\mathcal{D})$ is dense in $L^q(\mathcal{D})$, \cite{fournier_adams_sobolev}, Corollary 2.30). We relabel the countable family $(\phi_{nm})_{n,m \in \mathbb{N}}$ as $(\varphi_n)_{n \in \mathbb{N}}$, which is thus dense in $L^q(\mathcal{D})$. Second, let $(h_n)_{n\in \mathbb{N}} \subset C_c^{\infty}(\mathcal{D})$ be a dense subset of $C_c^{\infty}(\mathcal{D})$ for its LF-space topology (see Lemma \ref{lemma:d_separable}). We then define $E_q$ to be the set of all finite linear combinations of elements  of $(\varphi_n)$ and $(h_n)$ with rational coefficients :
\begin{align}
E_q &= \text{Span}{_\mathbb{Q}}\{\varphi_n, n \in \mathbb{N}\} + \text{Span}{_\mathbb{Q}}\{h_m, m \in \mathbb{N}\} \\
&= \bigcup_{n,m \in \mathbb{N}}\Big\{\sum_{i=1}^nq_i\varphi_i + \sum_{j=1}^m r_j h_j, (q_1,...,q_n,r_1,...,r_m)\in \mathbb{Q}^{n+m}\Big\} \label{eq:grand_phi}
\end{align}
Note that $E_q$ is countable, as a countable union of countable sets. We then define the family $(\Phi_n^q)$ to be an enumeration of $E_q$ : $E_q = \{\Phi_n^q, n \in \mathbb{N}\}$.  \\
\underline{Proof of $(i)$:} Suppose that $T = T_v$ for some $v \in L^p(\mathcal{D})$. Then the control \eqref{eq:LP_control_distrib} is obviously true. Now, suppose that this countable control holds : let us show that $T = T_v$ for some $v\in L^p(\mathcal{D})$.

We begin by showing that the map ${T}_{|E_q}$,the restriction of $T$ to the set $E_q$, can be uniquely extended to a continuous linear form $\tilde{T}$ over $L^q(\mathcal{D})$.
Begin with the fact that for all $f,g \in E_q$, then $f-g \in E_q$ and from equation \eqref{eq:LP_control_den},
\begin{align}\label{eq:control_E}
|T(f)-T(g)| = |T(f-g)| \leq C ||f-g||_{q}
\end{align}
Equation \eqref{eq:control_E} shows that ${T}_{|E_q}$ is Lipschitz over $E_q$ and therefore uniformly continuous on $E_q$. Since $\mathbb{R}$ is complete and $E_q$ is dense in $L^q(\mathcal{D})$, ${T}_{|E_q}$ can be uniquely extended by a map $\tilde{T}$ defined over $L^q(\mathcal{D})$, which is itself uniformly continuous (\cite{royden1988real}, Problem 44, p. 196). We briefly recall the construction procedure of $\tilde{T}$ over $L^q(\mathcal{D})$. Given $f \in L^q(\mathcal{D})$ and $(f_n) \subset E_q$ \textit{any} sequence such that $||f_n-f||_{L^q} \rightarrow 0$, one shows that the sequence $(T(f_n))_{n\in\mathbb{N}}$ is Cauchy, thus convergent and one sets $\tilde{T}(f):=\lim_n T(f_n)$. One proves that the value $\tilde{T}(f)$ does not depend on the sequence $(f_n)$, which implies that $\tilde{T}$ is well defined and coincides with ${T}$ on $E_q$.
 
We now check that $\tilde{T}$ remains linear. Let $f,g \in L^q(\mathcal{D})$ and $\lambda \in \mathbb{R}$. Let $(f_n),(g_n) \subset E_q$ and $(\lambda_n) \subset \mathbb{Q}$ be sequences such that $f_n \rightarrow f, g_n \rightarrow g$ both in $L^q(\mathcal{D})$ and $\lambda_n \rightarrow \lambda$. Then $\lambda_n f_n + g_n \rightarrow \lambda f+g$ in $L^q(\mathcal{D})$, and the sequence $(\lambda_nf_n + g_n)$ is contained in $E_q$. Since $\tilde{T}$ is well defined, we have that
\begin{align}
\tilde{T}(\lambda f + g) = \lim_{n\rightarrow \infty} T(\lambda_n f_n + g_n) = \lim_{n\rightarrow \infty} \lambda_n T(f_n) + T(g_n) = \lambda \tilde{T}(f) + \tilde{T}(g)
\end{align}
Thus, $\tilde{T}$ is a (uniformly) continuous linear form over $L^q(\mathcal{D})$. Riesz' representation lemma yields a function $v \in L^p(\mathcal{D})$ such that
\begin{align}\label{eq:riesz}
\forall f \in L^q(\mathcal{D}), \ \ \tilde{T}(f) = \int_{\mathcal{D}}f(x)v(x)dx
\end{align}
We now need to check that in fact $\tilde{T}(\varphi) = T(\varphi)$ if $\varphi \in C_c^{\infty}(\mathcal{D})$, to show that $\tilde{T}$ is indeed an extension of $T$. For this, notice that $T$ and $\tilde{T}$ both define continuous linear forms over $C_c^{\infty}(\mathcal{D})$, w.r.t. its LF-topology ($v$ lies in $L_{loc}^1(\mathcal{D})$). Note also that $T$ and $\tilde{T}$ coincide on $E_q$, by construction of $\tilde{T}$ :
\begin{align}
\forall n \in \mathbb{N},\ \ T(\Phi_n) - \tilde{T}(\Phi_n) = 0
\end{align}
But $E_q$ is chosen so that it contains $(h_n)$, which is a dense subset of $C_c^{\infty}(\mathcal{D})$. Given $\varphi \in C_c^{\infty}(\mathcal{D})$, consider $(j_n)$ a subsequence of $(h_n)$ such that $j_n \longrightarrow \varphi$ for the topology of $C_c^{\infty}(\mathcal{D})$. Then,
\begin{align}
(T-\tilde{T})(\varphi) = \lim_{n\rightarrow \infty}(T-\tilde{T})(j_n) = \lim_{n\rightarrow \infty} 0 = 0
\end{align}
which shows that in fact, $\tilde{T}(\varphi) = T(\varphi)$.\\
\underline{Proof of $(ii)$:} if $b$ can be extended to a continuous linear form over $L^q(\mathcal{D})$, then the estimate \eqref{eq:bilin_Lq} is obviously true, by continuity over $L^q(\mathcal{D})$ of the said extension. Suppose now that \eqref{eq:bilin_Lq} holds. Let $\varphi \in E_q$. Then $L_{\varphi}$, the continuous linear form over $C_c^{\infty}(\mathcal{D})$ defined by
\begin{align}
\forall \psi \in C_c^{\infty}(\mathcal{D}), \ \ L_{\varphi}(\psi) = b(\varphi,\psi)
\end{align}
verifies
\begin{align}
\forall \psi \in E_q, \ \ |L_{\varphi}(\psi)| \leq C||\varphi||_{q}||\psi||_{q}
\end{align}
From the point $(i)$, $L_{\varphi}$ is a regular distribution with a representer $v_{\varphi} \in L^p(\mathcal{D})$ which is unique in $L^p(\mathcal{D})$. Define the map $B : E_q \rightarrow L^p(\mathcal{D})$ by $B\varphi = v_{\varphi}$. Then $B$ verifies
\begin{align}
\forall \varphi \in E_q, \forall \psi \in L^q(\mathcal{D}), \ \ |\langle B\varphi,\psi \rangle_{L^p,L^q}| = |L_{\varphi}(\psi)| \leq C||\varphi||_{q}||\psi||_{q}
\end{align}
Taking the supremum w.r.t. $\psi\in L^q(\mathcal{D})$ yields
\begin{align}\label{eq:B_lpz}
\forall \varphi \in E_q, \ \ ||B\varphi||_{p} \leq C||\varphi||_{q}
\end{align}
Observe now that the bilinearity of $b$ yields $B(\varphi + \lambda\psi) = B\varphi + \lambda B\psi$ if $\varphi,\psi \in E_q$ and $\lambda \in \mathbb{Q}$. Taking the exact same steps as for the proof of point $(i)$ and using equation \eqref{eq:B_lpz}, $B : E_q \rightarrow L^p(\mathcal{D})$ is Lipschitz continuous over $E_q$, and can thus be uniquely extended as a uniformly continuous map $\tilde{B} : L^q(\mathcal{D}) \rightarrow L^p(\mathcal{D})$. This relies on the fact that $E_q$ is dense in $L^q(\mathcal{D})$ and that $L^q(\mathcal{D})$ is complete. As above, one checks that $\tilde{B}$ is linear. Being uniformly continuous, it is then a bounded operator from $L^q(\mathcal{D})$ to $L^p(\mathcal{D})$ (its adjoint $\tilde{B}^*$ is then automatically bounded). Denote by $\tilde{b}$ the continuous bilinear form over $L^q(\mathcal{D})$ defined by
\begin{align}
\tilde{b}(f,g) = \langle \tilde{B}f,g\rangle_{L^p,L^q} \ \ \ \forall f,g \in L^q(\mathcal{D})
\end{align}
We now need to check that $\tilde{b}$ indeed coincides with $b$ over $C_c^{\infty}(\mathcal{D})$, so that it is indeed an extension of $b$. For this, let $\varphi,\psi \in C_c^{\infty}(\mathcal{D})$ and $(\varphi_n), (\psi_n)$ two sequences of elements of $E_q$ that converge to $\varphi$ and $\psi$ respectively. Then $b$ and $\tilde{b}$ coincide on $E_q$:
\begin{align}\label{eq:b_btilde_coincide_E}
b(\varphi_n,\psi_m) = \tilde{b}(\varphi_n,\psi_m)
\end{align}
Observe the following chain of equalities, which rely on the sequential continuity (for the LF topology of $C_c^{\infty}(\mathcal{D})$) of the linear forms $\varphi \mapsto b(\varphi,\psi), \psi \mapsto b(\varphi,\psi)$ and $T_v : \varphi \mapsto T_v(\varphi) = \langle v,\varphi \rangle_{L^q,L^p}$ for any $v\in L^q(\mathcal{D})$, as well equation \eqref{eq:b_btilde_coincide_E}.
\begin{align}
b(\varphi,\psi) &= \lim_{n \rightarrow \infty} b(\varphi_n,\psi)  = \lim_{n \rightarrow \infty}\lim_{m \rightarrow \infty} b(\varphi_n,\psi_m) = \lim_{n \rightarrow \infty}\lim_{m \rightarrow \infty} \tilde{b}(\varphi_n,\psi_m) \nonumber \\
&= \lim_{n \rightarrow \infty}\lim_{m \rightarrow \infty}\langle \tilde{B}\varphi_n,\psi_m \rangle_{L^p,L^q} = \lim_{n \rightarrow \infty}\lim_{m \rightarrow \infty} T_{\tilde{B}\varphi_n}(\psi_m) = \lim_{n \rightarrow \infty} T_{\tilde{B}\varphi_n}(\psi) \nonumber\\
&= \lim_{n \rightarrow \infty} \langle \tilde{B}\varphi_n, \psi\rangle_{L^p,L^q} = \lim_{n \rightarrow \infty} \langle \varphi_n,\tilde{B}^* \psi\rangle_{L^q,L^p} = \lim_{n \rightarrow \infty}T_{\tilde{B}^*\psi}(\varphi_n) = T_{\tilde{B}^*\psi}(\varphi) \nonumber\\
&= \langle \varphi,\tilde{B}^* \psi\rangle_{L^q,L^p} = \langle \tilde{B}\varphi,\psi\rangle_{L^p,L^q} = \tilde{b}(\varphi,\psi)
\end{align}
The uniqueness of $b$ follows from the uniqueness of $\tilde{B}$ as an extension of $B$.

\end{proof}
\begin{proof}(Lemma \ref{lemma:u_alpha_phi_GP})
Let $(K_n)$ be an increasing sequence of compact subsets of $\mathcal{D}$ such that $\bigcup_nK_n = \mathcal{D}$. From the measurability of $U$ and Tonelli's theorem, $\omega \mapsto \int_{K_n}|U_{\omega}(x)|dx$ is measurable and we have that
\begin{align}\label{eq:traj_L1loc}
     \mathbb{E}\bigg[\int_{K_n}|U(x)|dx\bigg] = \int_{K_n}\mathbb{E}[|U(x)|]dx = \sqrt{\frac{2}{\pi}}\int_{K_n}\sigma(x)dx < +\infty
\end{align}
From equation \eqref{eq:traj_L1loc}, $\omega \mapsto \int_{K_n}|U_{\omega}(x)|dx$ is finite almost surely.
Since the family $(K_n)$ is countable, one obtains a set $\Omega_0 \subset \Omega$ of probability one such that for all $\omega\in\Omega_0$ and for all $n\in\mathbb{N}, \int_{K_n}|U_{\omega}(x)|dx <+\infty$.
Given now any compact subset $K$ of $\mathcal{D}$, there exists $N\in\mathbb{N}$ such that $K\subset K_N$ and thus for all $\omega\in\Omega_0,\ \int_K|U_{\omega}(x)|dx < +\infty$. Therefore, the sample paths of $U$ lie in $L_{loc}^1(\mathcal{D})$ almost surely. From this fact and Fubini's theorem, we next obtain that given any $\varphi \in C_c^{\infty}$ and $|\alpha|\leq m$, the following map
\begin{align}
U_{\varphi}^{\alpha} : \Omega \ni \omega \longmapsto \int_{\mathcal{D}}U_{\omega}(x)\partial^{\alpha}\varphi(x)dx
\end{align}
is a well defined random variable (i.e. it is measurable; see e.g. \cite{doob_sto_pro}, Theorem 2.7, p. 62). Moreover, one can show that it is a limit in probability of suitably chosen Riemann sums of the integrand (\cite{doob_sto_pro}, Theorem 2.8, p. 65). But here, those Riemann sums are all Gaussian random variables because $U$ is a Gaussian process. Thus $U_{\varphi}^{\alpha}$ is a Gaussian random variable. a a limit in probability of Gaussian random variables. This also shows that $\{U_{\varphi}^{\alpha}, \varphi \in C_c^{\infty}(\mathcal{D})\}$ is in fact a Gaussian process, since the linearity of $\partial^{\alpha}$ yields
\begin{align}
\sum_{i=1}^n a_i U_{\varphi_i}^{\alpha} = U_{(\sum_{i=1}^n a_i\varphi_i)}^{\alpha}
\end{align}
and thus $\sum_{i=1}^n a_i U_{\varphi_i}^{\alpha}$ is a Gaussian random variable. An alternative proof is found in \cite{bogachev1998gaussian}, Example 2.3.16. p. 58-59.
\end{proof}
\begin{proof}(Lemma \ref{lemma:extend_balpha})
First, the map $k$ is measurable over $\mathcal{D} \times \mathcal{D}$. Then, given a compact set $K \subset \mathcal{D} \times \mathcal{D}$, there exists a compact set $K_0 \subset \mathcal{D}$ such that $K \subset K_0\times K_0$ (see e.g. the text before equation \eqref{eq:K_0_D_0}). Then, using the Cauchy-Schwarz inequality for $k$,
\begin{align}
\int_K |k(x,y)|dxdy \leq \int_{K_0\times K_0}\sigma(x)\sigma(y)dxdy = \bigg(\int_{K_0}\sigma(x)dx\bigg)^2 < +\infty
\end{align}
Therefore, $k\in L_{loc}^1(\mathcal{D}\times \mathcal{D})$ and for all mutli-index $\alpha$, $b_{\alpha}$ is a bilinear continuous form over $C_c^{\infty}(\mathcal{D})$. From Lemma \ref{lemma:E_et_F_countable_line_bilin}, $b_{\alpha}$ can be uniquely extended to a continuous bilinear form over $L^2(\mathcal{D})$. Denote by $\mathcal{E}_k^{\alpha}$ the associated bounded operator over $L^2(\mathcal{D})$. We now need to show that $\mathcal{E}_k^{\alpha}$ is self-adjoint and positive. First note that for all $\varphi, \psi \in C_c^{\infty}(\mathcal{D})$,
\begin{align}\label{eq:ek_alpha_selfad}
\langle \mathcal{E}_k^{\alpha}\varphi, \psi \rangle_{L^2}
 &= \int_{\mathcal{D}\times \mathcal{D}}k(x,y)\partial^{\alpha}\varphi(x)\partial^{\alpha}\psi(y)dydx = \langle \varphi, \mathcal{E}_k^{\alpha}\psi \rangle_{L^2}
\end{align}
Equation \eqref{eq:ek_alpha_selfad}, conjoined with the density of $C_c^{\infty}(\mathcal{D})$ in $L^2(\mathcal{D})$ and the continuity of the bilinear form $(f,g) \mapsto \langle \mathcal{E}_k^{\alpha}f, g \rangle_{L^2}$ yields that $\langle \mathcal{E}_k^{\alpha}f, g \rangle_{L^2} = \langle f, \mathcal{E}_k^{\alpha}g \rangle_{L^2}$ for all $f,g \in L^2(\mathcal{D})$. Therefore $\mathcal{E}_k^{\alpha}$ is self-adjoint. For the positivity, consider again $\varphi \in C_c^{\infty}(\mathcal{D})$. Then from Fubini's theorem (justified below),
\begin{align}
\langle \mathcal{E}_k^{\alpha}\varphi, \varphi \rangle
&= \int_{\mathcal{D}\times \mathcal{D}}k(x,y)\partial^{\alpha}\varphi(x)\partial^{\alpha}\varphi(y)dydx = \int_{\mathcal{D}\times \mathcal{D}}\mathbb{E}[U(x)U(y)]\partial^{\alpha}\varphi(x)\partial^{\alpha}\varphi(y)dydx \nonumber \\
&= \mathbb{E}\bigg[\bigg(\int_{\mathcal{D}}U(x)\partial^{\alpha}\varphi(x)dx\bigg)^2\bigg] \geq 0
\end{align}
Indeed the following integrability condition holds, setting $K = \text{Supp}(\varphi)$ :
\begin{align}
\mathbb{E}\bigg[&\int_{\mathcal{D}\times\mathcal{D}}|\partial^{\alpha}\varphi(x)\partial^{\alpha}\varphi(y)U(x)U(y)|dxdy\bigg] = \int_{K\times K}|\partial^{\alpha}\varphi(x)\partial^{\alpha}\varphi(y)|\mathbb{E}[|U(x)U(y)|]dxdy \nonumber \\
&\leq \int_{K\times K}|\partial^{\alpha}\varphi(x)\partial^{\alpha}\varphi(y)|\sigma(x)\sigma(y)dxdy 
= \bigg(\int_{K}|\partial^{\alpha}\varphi(x)|\sigma(x)dx\bigg)^2 \nonumber\\
&\leq \sup_{x\in K}|\partial^{\alpha}\varphi(x)|^2 \bigg(\int_{K}\sigma(x)dx\bigg)^2 < +\infty \label{eq:e_k_al_pos}
\end{align}
Equation \eqref{eq:e_k_al_pos}, conjoined with the density of $C_c^{\infty}(\mathcal{D})$ in $L^2(\mathcal{D})$ and the continuity of the quadratic form $f \mapsto \langle \mathcal{E}_k^{\alpha}f, f \rangle_{L^2}$ yields that $\langle \mathcal{E}_k^{\alpha}f, f \rangle_{L^2} \geq 0$ for all $f \in L^2(\mathcal{D})$. Therefore $\mathcal{E}_k^{\alpha}$ is positive.
\end{proof}

\begin{proof} (Lemma \ref{lemma:HSpos}) 
Introduce $b_{\alpha}$ the continuous bilinear map over $C_c^{\infty}(\mathcal{D})$ defined by
\begin{align}
    b_{\alpha}(\varphi,\psi) &= \int_{\mathcal{D}\times\mathcal{D}}k(x,y)\partial^{\alpha}\varphi(x)\partial^{\alpha}\psi(y)dxdy = \int_{\mathcal{D}\times\mathcal{D}}\partial^{\alpha,\alpha}k(x,y)\varphi(x)\psi(y)dxdy \nonumber\\
    &= \langle \mathcal{E}_k^{\alpha}\varphi,\psi\rangle_{L^2} \label{eq:b_alpha_lemma_HSpos}
\end{align}
From Cauchy-Schwarz's inequality, it verifies
\begin{align}
   \forall \varphi,\psi \in C_c^{\infty}(\mathcal{D}), \  |b_{\alpha}(\varphi,\psi)| \leq ||\partial^{\alpha,\alpha}k||_2||\varphi||_2||\psi||_2
\end{align}
From Lemma \ref{lemma:extend_balpha}, there exists a unique bounded, self-adjoint and positive operator $B_{\alpha}$ over $L^2(\mathcal{D})$ such that $b_{\alpha}(\varphi,\psi) = \langle B_{\alpha}\varphi,\psi\rangle_{L^2}$ for all $\varphi, \psi \in C_c^{\infty}(\mathcal{D})$. The uniqueness of $B_{\alpha}$ and equation \eqref{eq:b_alpha_lemma_HSpos} yield $B_{\alpha} = \mathcal{E}_k^{\alpha}$, and thus $\mathcal{E}_k^{\alpha}$ is self-adjoint and positive.
\end{proof}

\begin{proof}(Lemma \ref{lemma:partial_kalpha_trace})
 \underline{$(i):$} first, let $j$ be such that $\lambda_j \neq 0$. Let $\varphi \in C_c^{\infty}(\mathcal{D})$. Then
\begin{align}
\lambda_j \bigg(\int_{\mathcal{D}}\phi_j(x)\partial^{\alpha}\varphi(x)dx\bigg)^2 &\leq \sum_{i=1}^{+\infty}\lambda_i\bigg(\int_{\mathcal{D}}\phi_i(x)\partial^{\alpha}\varphi(x)dx\bigg)^2 \nonumber \\
&\leq \sum_{i=1}^{+\infty}\lambda_i \int_{\mathcal{D}\times \mathcal{D}}\phi_i(x)\phi_i(y)\partial^{\alpha}\varphi(x)\partial^{\alpha}\varphi(y)dxdy \nonumber \\
&\leq \int_{\mathcal{D}\times \mathcal{D}}k(x,y)\partial^{\alpha}\varphi(x)\partial^{\alpha}\varphi(y)dxdy \nonumber \\
&\leq \int_{\mathcal{D}\times \mathcal{D}}\partial^{\alpha,\alpha}k(x,y)\varphi(x)\varphi(y)dxdy \nonumber \\
&\leq ||\partial^{\alpha,\alpha}k||_{L^2(\mathcal{D}\times\mathcal{D})}||\varphi||_{L^2(\mathcal{D})}^2
\end{align}
Therefore, from Lemma \ref{lemma:linear_form_sobolev}, $\partial^{\alpha}\phi_j \in L^2(\mathcal{D})$. \\
\underline{$(ii):$} introduce the finite rank kernel $k_n$ defined by
\begin{align}
k_n(x,y) = \sum_{i=1}^{n}\lambda_i\phi_i(x)\phi_i(y)
\end{align}
Then its mixed derivative $\partial^{\alpha,\alpha}k_n(x,y)$ is equal to $\sum_{i=1}^{n}\lambda_i\partial^{\alpha}\phi_i(x)\partial^{\alpha}\phi_i(y)$ in $ L^2(\mathcal{D}\times\mathcal{D})$ and the associated operator $\mathcal{E}_{k_n}^{\alpha}$ is trace class, with
\begin{align}\label{eq:trace_finite_rank}
\text{Tr}(\mathcal{E}_{k_n}^{\alpha}) &= \sum_{j=1}^{+\infty} \langle \mathcal{E}_{k_n}^{\alpha}\phi_j,\phi_j\rangle_{L^2} = \sum_{j=1}^{+\infty} \sum_{i=1}^n \lambda_i \langle \partial^{\alpha}\phi_i,\phi_j\rangle_{L^2}^2 \\
&= \sum_{i=1}^n \lambda_i \sum_{j=1}^{+\infty}\langle \partial^{\alpha}\phi_i,\phi_j\rangle_{L^2}^2 = \sum_{i=1}^n \lambda_i ||\partial^{\alpha}\phi_i||_{L^2}^2
\end{align}
Now, observe that $\mathcal{E}_{k_n}^{\alpha} \leq \mathcal{E}_{k}^{\alpha}$ in the sense of the Loewner order. Indeed, let first $\varphi \in C_c^{\infty}(\mathcal{D})$:
\begin{align}
\langle (\mathcal{E}_{k}^{\alpha} - \mathcal{E}_{k_n}^{\alpha})\varphi,\varphi\rangle_{L^2} = \langle (\mathcal{E}_{k} - \mathcal{E}_{k_n})\partial^{\alpha}\varphi,\partial^{\alpha}\varphi\rangle_{L^2}
= \sum_{i=n+1}^{+\infty} \lambda_i \langle \phi_i,\partial^{\alpha}\varphi\rangle_{L^2}^2 \geq 0
\end{align}
The density of $C_c^{\infty}(\mathcal{D})$ in $L^2(\mathcal{D})$ and the continuity of the quadratic form $f \mapsto \langle (\mathcal{E}_{k}^{\alpha} - \mathcal{E}_{k_n}^{\alpha})f,f\rangle_{L^2}$ over $L^2(\mathcal{D})$ yields indeed that $\mathcal{E}_{k_n}^{\alpha} \leq \mathcal{E}_{k}^{\alpha}$. Taking the trace :
\begin{align}\label{eq:tr_comp}
\sum_{i=1}^n \lambda_i ||\partial^{\alpha}\phi_i||_{L^2}^2 = \text{Tr}(\mathcal{E}_{k_n}^{\alpha}) = \sum_{j=1}^{+\infty} \langle \mathcal{E}_{k_n}^{\alpha}\phi_j,\phi_j\rangle_{L^2} \leq \sum_{j=1}^{+\infty} \langle \mathcal{E}_{k}^{\alpha}\phi_j,\phi_j\rangle_{L^2} = \text{Tr}(\mathcal{E}_{k}^{\alpha}) 
\end{align}
Taking the limit when $n$ goes to infinity yields $\sum_{i=1}^{+\infty} \lambda_i ||\partial^{\alpha}\phi_i||_{L^2}^2 \leq \text{Tr}(\mathcal{E}_{k}^{\alpha})$. Suppose now that $\text{Tr}(\mathcal{E}_{k}^{\alpha}) < + \infty$. Equation \eqref{eq:tr_comp} shows that the series of functions $\sum_i \lambda_i \partial^{\alpha}\phi_i \otimes\partial^{\alpha}\phi_i$ converges in norm in $L^2(\mathcal{D}\times\mathcal{D})$. Moreover, we check that it is equal to $\partial^{\alpha,\alpha}k$ : taking $\varphi \in C_c^{\infty}(\mathcal{D}\times\mathcal{D})$, then
\begin{align}
\int_{\mathcal{D}\times\mathcal{D}} k(x,y)\partial^{\alpha,\alpha}\varphi(x,y)dxdy
 &= \sum_{i}\lambda_i \int_{\mathcal{D}\times\mathcal{D}} \phi_i(x)\phi_i(y)\partial^{\alpha,\alpha}\varphi(x,y)dxdy \\
 &= \sum_{i}\lambda_i\int_{\mathcal{D}\times\mathcal{D}} \partial^{\alpha}\phi_i(x)\partial^{\alpha}\phi_i(y)\varphi(x,y)dxdy \\
 &= \int_{\mathcal{D}\times\mathcal{D}} \bigg(\sum_{i}\lambda_i\partial^{\alpha}\phi_i(x)\partial^{\alpha}\phi_i(y)\bigg)\varphi(x,y)dxdy 
\end{align}
We can then write,, following the steps of equation \eqref{eq:trace_finite_rank}
\begin{align}
\text{Tr}(\mathcal{E}_{k}^{\alpha})  = \sum_{j=1}^{+\infty} \langle \mathcal{E}_{k}^{\alpha}\phi_j,\phi_j\rangle_{L^2} =  \sum_{j=1}^{+\infty} \sum_{i}\lambda_i \langle\partial^{\alpha}\phi_i,\phi_j\rangle_{L^2}^2 = \sum_{i=1}^{+\infty} \lambda_i ||\partial^{\alpha}\phi_i||_{L^2}^2 
\end{align}
Suppose now that $\sum_{i=1}^{+\infty} \lambda_i ||\partial^{\alpha}\phi_i||_{L^2}^2 < +\infty$. Then as observed before, the series of functions $\sum_i \lambda_i \partial^{\alpha}\phi_i \otimes\partial^{\alpha}\phi_i$ converges in norm in $L^2(\mathcal{D}\times\mathcal{D})$, one verifies that $\partial^{\alpha,\alpha}k$ exists in $L^2(\mathcal{D})$ and is in fact given by
\begin{align}
\partial^{\alpha,\alpha}k = \sum_i \lambda_i \partial^{\alpha}\phi_i \otimes\partial^{\alpha}\phi_i \ \ \ \text{in} \ \ \ L^2(\mathcal{D}\times\mathcal{D})
\end{align}
Finally,
\begin{align}
\sum_{i=1}^{+\infty} \lambda_i ||\partial^{\alpha}\phi_i||_{L^2}^2  &= \sum_{i=1}^{+\infty} \lambda_i \sum_{j} \langle\partial^{\alpha}\phi_i,\phi_j\rangle_{L^2}^2 \\
 &= \sum_{j}\sum_{i=1}^{+\infty} \lambda_i\bigg(\int_{\mathcal{D}}\partial^{\alpha}\phi_i(x)\phi_j(x)dx\bigg)^2 \\
 &= \sum_{j} \int_{\mathcal{D}\times\mathcal{D}}\sum_{i} \lambda_i\partial^{\alpha}\phi_i(x)\partial^{\alpha}\phi_i(y)\phi_j(x)\phi_j(y)dxdy \\
 &= \sum_{j}\langle \mathcal{E}_k^{\alpha}\phi_j,\phi_j\rangle_{L^2} = \text{Tr}(\mathcal{E}_k^{\alpha})
\end{align}
Therefore $\mathcal{E}_k^{\alpha}$ is trace class and $\text{Tr}(\mathcal{E}_k^{\alpha}) = \sum_{i=1}^{+\infty} \lambda_i ||\partial^{\alpha}\phi_i||_{L^2}^2 $. For asymmetric derivatives, simply observe that for all $|\alpha|,|\beta|\leq m$,
\begin{align}
    ||\partial^{\alpha}\phi_i \otimes \partial^{\beta}\phi_i||_2 = ||\partial^{\alpha}\phi_i||_2|| \partial^{\beta}\phi_i||_2 \leq \frac{||\partial^{\alpha}\phi_i||_2^2 + ||\partial^{\beta}\phi_i||_2^2}{2}
\end{align}
Therefore the norm convergence of the series $\sum_{i\in\mathbb{N}} \lambda_i ||\partial^{\alpha}\phi_i\otimes\partial^{\alpha}\phi_i||_{L^2}$ for all $|\alpha|\leq m$ implies that of all the series of the form $\sum_{i\in\mathbb{N}} \lambda_i ||\partial^{\alpha}\phi_i\otimes\partial^{\beta}\phi_i||_{L^2}$ converge, provided that $|\alpha|\leq m$ and $|\beta|\leq m$. As previously, one readily checks that $\partial^{\alpha,\beta}k = \sum_{i=0}^{\infty}\lambda_i\partial^{\alpha}\phi_i\otimes\partial^{\beta}\phi_i$.
\end{proof}
\bibliographystyle{abbrv}
\bibliography{bibliography}

\begin{thebibliography}{10}

\bibitem{adler2007}
R.~J. Adler and J.~E. Taylor.
\newblock {\em Random Fields and Geometry}.
\newblock Springer-Verlag New York, NY, 2007.

\bibitem{azais_level_2009}
J.-M. Aza{\"i}s and M.~Wschebor.
\newblock {\em {Level sets and extrema of random processes and fields}}.
\newblock {Wiley \& Sons}, 2009.

\bibitem{agnan2004}
A.~Berlinet and C.~Thomas-Agnan.
\newblock {\em Reproducing Kernel {H}ilbert Spaces in Probability and
  Statistics}.
\newblock Springer US, 2004.

\bibitem{bogachev1998gaussian}
V.~I. Bogachev.
\newblock {\em {G}aussian measures}.
\newblock Number~62 in Mathematical Surveys and Monographs. American
  Mathematical Soc., 1998.

\bibitem{bogachev2007measure}
V.~I. Bogachev and M.~A.~S. Ruas.
\newblock {\em Measure theory}, volume~1.
\newblock Springer, 2007.

\bibitem{brenner2008mathematical}
S.~C. Brenner, L.~R. Scott, and L.~R. Scott.
\newblock {\em The mathematical theory of finite element methods}, volume~3.
\newblock Springer, 2008.

\bibitem{brezis2011functional}
H.~Br{\'e}zis.
\newblock {\em Functional analysis, {S}obolev spaces and partial differential
  equations}.
\newblock Springer New York, NY, 2011.

\bibitem{brislawn1988kernels}
C.~Brislawn.
\newblock Kernels of trace class operators.
\newblock {\em Proceedings of the American Mathematical Society},
  104(4):1181--1190, 1988.

\bibitem{CHEN_owhadi_2021}
Y.~Chen, B.~Hosseini, H.~Owhadi, and A.~M. Stuart.
\newblock Solving and learning nonlinear {PDE}s with {G}aussian processes.
\newblock {\em Journal of Computational Physics}, 447:110668, 2021.

\bibitem{chobanjan1977gaussian}
S.~Chobanjan and V.~Tarieladze.
\newblock {G}aussian characterizations of certain {B}anach spaces.
\newblock {\em Journal of Multivariate Analysis}, 7(1):183--203, 1977.

\bibitem{ciesielski1991modulus}
Z.~Ciesielski.
\newblock Modulus of smoothness of the {B}rownian paths in the {$L^p$} norm.
\newblock {\em Constructive theory of functions (Varna, Bulgaria, 1991)}, pages
  71--75, 1991.

\bibitem{ciesielski1993quelques}
Z.~Ciesielski.
\newblock Quelques espaces fonctionnels associ{\'e}s {\`a} des processus
  gaussiens.
\newblock {\em Studia Mathematica}, 107:171--204, 1993.

\bibitem{cirel1976norms}
B.~S. Cirel'son, I.~A. Ibragimov, and V.~N. Sudakov.
\newblock Norms of {G}aussian sample functions.
\newblock In {\em Proceedings of the Third Japan—USSR Symposium on
  Probability Theory}, pages 20--41. Springer, 1976.

\bibitem{di2012hitchhiker}
E.~Di~Nezza, G.~Palatucci, and E.~Valdinoci.
\newblock Hitchhiker{'}s guide to the fractional {S}obolev spaces.
\newblock {\em Bulletin des sciences math{\'e}matiques}, 136(5):521--573, 2012.

\bibitem{doob1937stochastic}
J.~L. Doob.
\newblock Stochastic processes depending on a continuous parameter.
\newblock {\em Transactions of the American Mathematical Society},
  42(1):107--140, 1937.

\bibitem{doob_sto_pro}
J.~L. Doob.
\newblock {\em Stochastic processes}.
\newblock Wiley Classics Library. Wiley-Interscience, 1990.

\bibitem{dudley1967sizes}
R.~M. Dudley.
\newblock The sizes of compact subsets of {H}ilbert space and continuity of
  {G}aussian processes.
\newblock {\em Journal of Functional Analysis}, 1(3):290--330, 1967.

\bibitem{evans1998}
L.~Evans.
\newblock {\em Partial Differential Equations}.
\newblock Graduate studies in mathematics. American Mathematical Society, 1998.

\bibitem{gapaillard8processus}
J.~Gapaillard and J.~Michaux.
\newblock Sur les processus lin{\'e}aires d{\'e}finis sur un espace
  nucl{\'e}aire.
\newblock In {\em Annales de la Facult{\'e} des sciences de Toulouse:
  Math{\'e}matiques}, volume~8, pages 75--92, 1986.

\bibitem{gelfand1964generalized}
I.~M. Gel'fand and G.~E. Shilov.
\newblock {\em Generalized functions, Vol. 4: applications of harmonic
  analysis}.
\newblock Academic Press, 1964.

\bibitem{henderson_arxiv}
I.~Henderson, P.~Noble, and O.~Roustant.
\newblock Stochastic processes under linear differential constraints :
  Application to {G}aussian process regression for the 3 dimensional free space
  wave equation, 2021.

\bibitem{ibragimov1994conditions}
I.~Ibragimov.
\newblock Conditions for {G}aussian homogeneous fields to belong to classes
  {$H^r_p$}.
\newblock {\em Journal of Mathematical Sciences}, 68(4):484--497, 1994.

\bibitem{hoffman_jorgensen_book}
J.~J{\o}rgensen.
\newblock {\em Probability With a View Towards Statistics, Volume {II}}.
\newblock Chapman and Hall/CRC, 1994.

\bibitem{kanagawa2018gaussian}
M.~Kanagawa, P.~Hennig, D.~Sejdinovic, and B.~K. Sriperumbudur.
\newblock {G}aussian processes and kernel methods: A review on connections and
  equivalences.
\newblock {\em arXiv preprint arXiv:1807.02582}, 2018.

\bibitem{kerkyacharian2018regularity}
G.~Kerkyacharian, S.~Ogawa, P.~Petrushev, and D.~Picard.
\newblock Regularity of {G}aussian processes on {D}irichlet spaces.
\newblock {\em Constructive Approximation}, 47(2):277--320, 2018.

\bibitem{PhamTheLai1973-1974}
P.~T. Lai.
\newblock Noyaux d'{A}gmon.
\newblock {\em Séminaire Jean Leray}, 277(4):1--37, 1973-1974.

\bibitem{linde1980characterization}
V.~Linde, V.~I. Tarieladze, and S.~A. Chobanyan.
\newblock Characterization of certain classes of {B}anach spaces by properties
  of {G}aussian measures.
\newblock {\em Theory of Probability \& Its Applications}, 25(1):159--164,
  1980.

\bibitem{nikol2012approximation}
S.~M. Nikol'skii.
\newblock {\em Approximation of functions of several variables and imbedding
  theorems}, volume 205.
\newblock Springer Science \& Business Media, 2012.

\bibitem{owhadi_bayes_homog}
H.~Owhadi.
\newblock {B}ayesian numerical homogenization.
\newblock {\em Multiscale Modeling \& Simulation}, 13(3):812--828, 2015.

\bibitem{owhadi2017separability}
H.~Owhadi and C.~Scovel.
\newblock Separability of reproducing kernel spaces.
\newblock {\em Proceedings of the American Mathematical Society},
  145(5):2131--2138, 2017.

\bibitem{pinns}
M.~Raissi, P.~Perdikaris, and G.~Karniadakis.
\newblock Physics-informed neural networks: A deep learning framework for
  solving forward and inverse problems involving nonlinear partial differential
  equations.
\newblock {\em Journal of Computational Physics}, 378:686--707, 2019.

\bibitem{raissi2017}
M.~Raissi, P.~Perdikaris, and G.~E. Karniadakis.
\newblock Machine learning of linear differential equations using {G}aussian
  processes.
\newblock {\em Journal of Computational Physics}, 348:683--693, 2017.

\bibitem{fournier_adams_sobolev}
J.~J. F.~F. Robert A.~Adams.
\newblock {\em {S}obolev spaces}.
\newblock Pure and Applied Mathematics. Academic Press, 2 edition, 2003.

\bibitem{royden1988real}
H.~L. Royden and P.~Fitzpatrick.
\newblock {\em Real analysis}, volume~32.
\newblock Macmillan New York, 1988.

\bibitem{roynette1993mouvement}
B.~Roynette.
\newblock Mouvement {B}rownien et espaces de {B}esov.
\newblock {\em Stochastics: An International Journal of Probability and
  Stochastic Processes}, 43(3-4):221--260, 1993.

\bibitem{SCHEUERER2010}
M.~Scheuerer.
\newblock Regularity of the sample paths of a general second order random
  field.
\newblock {\em Stochastic Processes and their Applications},
  120(10):1879--1897, 2010.

\bibitem{stein1970singular}
E.~M. Stein.
\newblock {\em Singular integrals and differentiability properties of
  functions}, volume~2.
\newblock Princeton university press, 1970.

\bibitem{steinwart2019convergence}
I.~Steinwart.
\newblock Convergence types and rates in generic {K}arhunen-{L}oeve expansions
  with applications to sample path properties.
\newblock {\em Potential Analysis}, 51(3):361--395, 2019.

\bibitem{talagrand2014upper}
M.~Talagrand.
\newblock {\em Upper and lower bounds for stochastic processes}, volume~60.
\newblock Springer, 2014.

\bibitem{treves2006topological}
F.~Treves.
\newblock {\em Topological Vector Spaces, Distributions and Kernels}.
\newblock Dover books on mathematics. Dover Publications, 2006.

\bibitem{triebel1968approximationszahlen}
H.~Triebel.
\newblock {\"U}ber die approximationszahlen der einbettungsoperatoren
  305-01305-01305-01.
\newblock {\em Archiv der Mathematik}, 19(3):305--312, 1968.

\bibitem{van2011information}
A.~Van Der~Vaart and H.~Van~Zanten.
\newblock Information rates of nonparametric {G}aussian process methods.
\newblock {\em Journal of Machine Learning Research}, 12(6), 2011.

\bibitem{wahlstrom2013}
N.~Wahlstrom, M.~Kok, T.~B. Sch{\"o}n, and F.~Gustafsson.
\newblock Modeling magnetic fields using {G}aussian processes.
\newblock {\em 2013 IEEE International Conference on Acoustics, Speech and
  Signal Processing}, pages 3522--3526, 2013.

\end{thebibliography}

\end{document}